%% file: main.tex
\providecommand{\main}{.}
\documentclass{article}

\usepackage{amsmath,amsfonts,amsthm}
\usepackage{tikz}
\usepackage{enumitem}
\usepackage{bm}
\usepackage{graphicx}
\usepackage{geometry}
\usepackage{\main/commands}
\usepackage{subfiles}
\usepackage[]{algorithm2e}
\usepackage{subcaption}
\usepackage{setspace} 
\newtheorem{theorem}{Theorem}

\newtheorem{definition}{Definition}
\newtheorem{remark}{Remark}
\newtheorem{lemma}{Lemma}
\newtheorem{proposition}{Proposition}
\newcommand{\br}{\bar{r}}
\newcommand{\ba}{\bar{a}}
\newcommand{\bb}{\bar{b}}
\newcommand{\bc}{\bar{c}}
\newcommand{\by}{\bar{y}}
\newcommand{\be}{\bar{\ve}}

\def\biblio{\bibliographystyle{alpha}\bibliography{\main/bib}}
\begin{document} 
\def\biblio{}

\title{On fast-slow consensus networks with a dynamic weight}
\author{Hildeberto Jard\'on-Kojakhmetov\thanks{Technical University of Munich (TUM), Faculty of Mathematics, Boltzmannstr. 3, 85748 Garching b. M{\"u}nchen, Germany; e-mail: \tt{h.jardon.kojakhmetov@tum.de} } \, and Christian Kuehn\thanks{Technical University of Munich (TUM), Faculty of Mathematics, Boltzmannstr. 3, 85748 Garching b. M{\"u}nchen, Germany; e-mail: \tt{ckuehn@ma.tum.de}}}

\maketitle

\begin{abstract}
We study dynamic networks under an undirected consensus communication protocol and with one state-dependent weighted edge. We assume that the aforementioned dynamic edge can take values over the whole real numbers, and that its behaviour depends on the nodes it connects and on an extrinsic slow variable. We show that, under mild conditions on the weight, there exists a reduction such that the dynamics of the network are organized by a transcritical singularity. As such, we detail a slow passage through a transcritical singularity for a simple network, and we observe that an exchange between consensus and clustering of the nodes is possible. In contrast to the classical planar fast-slow transcritical singularity, the network structure of the system under consideration induces the presence of a maximal canard. Our main tool of analysis is the blow-up method. Thus, we also focus on tracking the effects of the blow-up transformation on the network's structure. We show that on each blow-up chart one recovers a particular dynamic network related to the original one. We further indicate a numerical issue produced by the slow passage through the transcritical singularity. 
\end{abstract}

\newpage

\tableofcontents

\newpage

\section{Introduction}\label{sec:intro}

\subfile{subfiles/intro.tex}
\section{Preliminaries}\label{sec:prel}
\subfile{subfiles/prel.tex}
\section{A triangle motif}\label{sec:triangle}

\subfile{subfiles/triangle.tex}
\subsection{Blow-up analysis}\label{sec:blowup}

\subfile{subfiles/blowup.tex}
\section{Some generalizations}\label{sec:gen}

\subfile{subfiles/generalizations.tex}
\section{Summary and outlook}\label{sec:con}

\subfile{subfiles/conclusions.tex}

\section*{Acknowledgements}
HJK is supported by a Technical University Foundation Fellowship and by an Alexander von Humboldt Research Fellowship. CK is partially supported by a Lichtenberg Professorship of the VolkswagenStiftung. CK is also partially supported by project B10 of the SFB/TR109 ``Discretization in Geometry and Dynamics'' funded by the German Science Foundation (DFG).  

\appendix  
\section{Numerical analysis and simulations}\label{sec:num}

\subfile{subfiles/num.tex}

%\textcolor{red}{TODO: (a) 'uniformize' references; (b) change the format so that we have 
%less pages; maybe even use directly the JNLS style; (c) change 'constraint equation' and 'layer equation' to small letters everywhere (in case I have overlooked sth); (d) change order notation from non-caligraphic to caligraphic;}

\clearpage

%\doublespacing

\bibliographystyle{plain}
\bibliography{\main/bib}

\end{document}

%% file: subfiles/intro.tex
A wide range of scientific disciplines, such as biochemistry \cite{aral2009distinguishing,xia2004analyzing}, economics \cite{schweitzer2009economic}, social sciences \cite{moreno1934shall,proskurnikov2017tutorial,proskurnikov2018tutorial}, epidemiology \cite{barabasi2016network,pastor2001epidemic}, among many others, benefit from the progress in network theory. Network science represents an important paradigm for modeling and analysis of complex systems. In network theory, individuals are represented as vertices on a graph. These individuals, or agents, can be persons in a community, robots working in an assembly line, computers in an office building, particles of a chemical substance, etc\'etera. The interactions between such individuals are then represented as edges or links of the graph. For example, if two persons talk to each other, and one influences the other, we then associate an edge on the graph to such an activity. Similarly, interactions can be identified with edges in all other aforementioned examples. Thus, the individual's own dynamics (the nodes' dynamics) and the interaction it has with other individuals of the network (the topology of the network) will not only determine its own fate, but that of the entire group of individuals. This rather convenient way of describing complicated dynamic behavior is quite powerful and has attracted an enormous scientific interest \cite{albert2002statistical,barrat2004architecture,boccaletti2006complex,strogatz2001exploring}. 

From an applied mathematical perspective, topics such as stability, convergence rates, synchronization, connectivity, robustness, and many others can all be formally described and have important implications in other sciences. In a large part of the mathematical studies of networks, one considers that the interactions between the agents are fixed~\cite{BarratBarthelemyVespignani}; see also Section \ref{sec:prel}. In another large part of the theory one frequently assumes that the network structure evolves without dynamics at nodes~\cite{vanderHofstad}. In most cases, these assumptions are a simplification since it may be expected that there is coupled dynamics \emph{of} and \emph{on} the network, i.e., one has to deal with adaptive (or co-evolutionary) networks~\cite{GrossSayama}. A crucial assumption to approximate an adaptive network by a partially static one with either just dynamics on or dynamics of the network is time scale separation~\cite{kuehn2015multiple}. Yet, if one assumes that either the dynamics on the nodes or the dynamics of the edges are infinitely slow, hence static, leads to a singular limit description. This limit is known to miss adaptive network dynamics effects induced by the interaction of dynamical variables for finite time scale separation~\cite{KuehnNetworks}. Also from the viewpoint of applications, a finite but large time scale separation is far more reasonable. As an example, consider a group of people that communicate with each other daily but whose mutual influences shape the way they handle elections. One of such activities occurs in time scales of minutes or hours, while the other in time scales of years, yet both are interrelated in a complex manner. Similar examples where different sorts of relations occur at distinct time scales can be found in population dynamics, telecommunication networks, power grids, etc\'etera. So, although time scales add an extra level of difficulty to the analysis of networks, they may be useful for a more accurate representation of certain phenomenona. On the other hand, dynamical systems with two or more time scales have also been of interest from many perspectives, particularly in applied mathematics. The overall idea is to distinguish slow from fast subprocesses, analyze them separately, and then come up with an appropriate description of the problem~\cite{jones1995geometric,Kaper,kuehn2015multiple,OMalley1991,verhulst2005methods}. This basic idea can be made rigorous, and has proven to be powerful. However, there are generic complex systems in which the time scale separation can no longer be clearly distinguished. Thus more advanced mathematical techniques are required to analyze multi-scale adaptive networks. 

In this article we bring together network and multi-scale theories to study a class of adaptive networks. We are interested in networks whose agents communicate in a rather simple way, known as linear average consensus protocol (see the details in Section \ref{sec:consensus}). This type of communication has been largely studied due to its relevance in all kinds of sciences \cite{mesbahi2010graph,ren2005survey}. On this class of networks we assume that there is one interaction or communication link that slowly changes over time and investigate the implications of it. We shall see that the aforementioned setting leads to a nontrivial problem from both, network and multi-scale, contexts. As a result we describe the overall behavior of the network by adequately incorporating techniques from consensus dynamics and geometric singular perturbation theory. 

The forthcoming parts of this work are arranged as follows: in Section \ref{sec:prel} we provide a short technical introduction to the main topics of this paper, namely fast-slow systems and consensus networks. In Section \ref{sec:triangle} we present our main contribution, which consists in the analysis of a simple network that has a dynamic weight and whose overall dynamics evolve in two time scales. Next, in Section \ref{sec:gen} we show that, in qualitative terms, the analysis performed for the aforementioned simple network can be extended to arbitrary networks with one dynamic edge. We finish in Section \ref{sec:con} with concluding remarks and an outlook on future research.

\biblio

%% file: subfiles/prel.tex
	In this section we provide a brief recollection of the two mathematical areas that come together in this paper. We first state what a fast-slow system formally is, the concept of normal hyperbolicity, and two relevant geometric techniques of analysis. Afterwards, to place our work into context, we recall and provide appropriate references to some of the relevant results on dynamic networks.

\subsection{Fast-slow systems}\label{sec:sfs}

A fast-slow system is a singularly perturbed Ordinary Differential Equation (ODE) of the form
\begin{equation}\label{eq:pre-sfs1}
	\begin{split}
		\ve \dot x &= f(x,y,\ve)\\
		\dot y &= g(x,y,\ve),
	\end{split}
\end{equation}
where $x\in\bbR^m$ and $y\in\bbR^n$ are, respectively, the fast and slow variables, and where $0<\ve\ll1$ is a small parameter accounting for the time scale difference between the variables. The overdot denotes derivative with respect to the slow time $\tau$. By defining the fast time $t=\tau/\ve$, one can rewrite \eqref{eq:pre-sfs1} as
\begin{equation}\label{eq:pre-sfs2}
	\begin{split}
		x' &= f(x,y,\ve)\\
		y' &= \ve g(x,y,\ve),
	\end{split}
\end{equation}
where the prime denotes the derivative with respect to the fast time $t$. Observe that, for $\ve>0$, the only difference between \eqref{eq:pre-sfs1} and \eqref{eq:pre-sfs2} is their time-parametrization. Therefore, we say that \eqref{eq:pre-sfs1} and \eqref{eq:pre-sfs2} are equivalent.

Although there are several approaches~\cite{kuehn2015multiple} to the analysis of fast-slow systems, e.g.~classical asymptotics~\cite{eckhaus2011asymptotic,eckhaus2011matched,OMalley1991,verhulst2005methods}, here we take a geometric approach \cite{fenichel1979geometric,jones1995geometric}, which is called Geometric Singular Perturbation Theory. The overall idea is to consider \eqref{eq:pre-sfs1} and \eqref{eq:pre-sfs2} restricted to $\ve=0$, understand the resulting systems, and then use perturbation results to obtain a description of \eqref{eq:pre-sfs1} and \eqref{eq:pre-sfs2} for $\ve>0$ sufficiently small. Therefore, two important subsystems to be considered are
\begin{equation}
\label{eq:singlimits}
	\begin{split}
		0 &= f(x,y,0)\\
		\dot y &= g(x,y,0),
	\end{split}\hspace{5cm}
	\begin{split}
		x' &= f(x,y,0)\\
		y' &= 0,
	\end{split}
\end{equation}
which are called \emph{the constraint equation}~\cite{takens1976constrained} (or slow subsystem or reduced system) and \emph{the layer equation} (or fast subsystem) respectively. It is important to note that the constraint and layer equations are not equivalent any more, there are even different classes of differential equations as the constraint equation is a differential-algebraic equation~\cite{kunkel2006differential}, while the layer equation is an ODE, where the slow variables $y$ can be viewed as parameters. In some sense the time scale separation is infinitely large between two singular limit systems~\eqref{eq:singlimits}. However, a geometric object that relates the two is \emph{the critical manifold}.

\begin{definition} The critical manifold of a fast-slow system is defined by
\begin{equation}
	\cC_0=\left\{ (x,y)\in\bbR^m\times\bbR^n\,|\, f(x,y,0)=0 \right\}.
\end{equation}
\end{definition}

The critical manifold is, on the one hand, the set of solutions of the algebraic equation in the constraint equation, and on the other hand, the set of equilibrium points of the layer equation. There is an important property that critical manifolds may have, called \emph{normal hyperbolicity}. 
\begin{definition} A point $p\in\cC_0$ is called hyperbolic if the eigenvalues of the matrix $\textnormal D_xf(p,0)$, where $\textnormal D_x$ denotes the total derivative with respect to $x$, have non-zero real part. The critical manifold $\cC_0$ is called \emph{normally hyperbolic} if every point $p\in\cC_0$ is hyperbolic. On the contrary, if for a point $p\in\cC_0$ we have that $\textnormal D_xf(p,0)$ has at least one eigenvalue on the imaginary axis, we then call $p$ \emph{non-hyperbolic}.
\end{definition}

In a general sense, whether a critical manifold has non-hyperbolic points or not, dictates the type of mathematical techniques that are suitable for analysis. For the case when the critical manifold is normally hyperbolic, Fenichel's theory \cite{fenichel1979geometric} (see also \cite{tikhonov1952systems} and \cite[Chapter 3]{kuehn2015multiscale}) asserts that, under compactness of the critical manifold, the constraint and the layer equations give a good enough approximation of the dynamics near $\cC_0$ of the fast-slow system for $\ve>0$ sufficiently small. In the normally hyperbolic case for $0<\ve\ll1$, there exists a slow manifold $\cC_\ve$, which can be viewed as a perturbation of $\cC_0$; see also~\cite{fenichel1979geometric,kuehn2015multiscale}.

The case when the critical manifold has non-hyperbolic points is considerably more difficult. One mathematical technique that has proven highly useful for the analysis in such a scenario is the \emph{blow-up method} \cite{dumortier1996canard}. Briefly speaking, the blow-up method consists on a well-suited generalized polar change of coordinates. What one aims to gain with such a coordinate transformation is enough hyperbolicity so that the dynamics can be analyzed using standard techniques of dynamical systems. Nowadays, the blow-up method is widely used to analyze the dynamics of fast-slow systems having non-hyperbolic points in a broad range of theoretical contexts and applications. For detailed information on the blow-up technique the reader may refer to \cite{dumortier1996canard,jardon2019survey,krupa2001extending}, \cite[Chapter 7]{kuehn2015multiscale} and references therein.

\subsection{Consensus networks}\label{sec:consensus}

In this section we formally introduce the type of consensus problems on an adaptive network which we are concerned with in this work. Let us start by introducing some notation: we denote by $\cG=\left\{\cV,\cE,\cW \right\}$ an undirected weighted graph where $\cV=\left\{ 1,\ldots,m \right\}$ denotes the set of vertices, $\cE=\left\{ e_{ij}\right\}$ the set of edges and $\cW=\left\{ w_{ij}\right\}$ the set of weights. We assume that the graph is undirected, that there are only simple edges, and that there are no self-loops, that is $e_{ij}=e_{ji}$ and $e_{ii}\notin\cE$. To each edge $e_{ij}$ we assign a weight $w_{ij}\in\bbR$ and thus we identify  the presence (resp.~absence) of an edge with a nonzero (resp.~zero) weight. Moreover, we shall say that a graph is \emph{unweighted} if all the nonzero weights are equal to one. The Laplacian \cite{merris1994laplacian} of the graph $\cG$ is denoted by $L=[l_{ij}]$ and is defined by
\begin{equation}\label{eq:l}
	l_{ij}=\begin{cases}
		-w_{ij}, & i\neq j\\
		\sum_{j=1}^m w_{ij}, & i=j.
	\end{cases}
\end{equation}

\begin{remark}
The majority of the scientific work regarding adaptive/dynamic networks considers non-negative weights. One of the reasons for such considerations is that the spectrum of the Laplacian matrix is well identified \cite{barrat2004architecture,mohar1991laplacian,olfati2007consensus}, which simplifies the analysis. When the weights are allowed to be positive and negative one usually refers to $L$ as a \emph{signed Laplacian}. Difficulties arise due to the fact that many of the convenient properties of non-negatively weighted Laplacians do not hold for signed Laplacians. In some part of the literature, see for example~\cite{altafini2013consensus,proskurnikov2016opinion}, the diagonal entries of the Laplacian matrix are rather defined by $\sum_{j=1}^m |w_{ij}|$. In this case, however, the Laplacian matrix is positive semi-definite and the potential loss of stability due to dynamic weights (the main topic of this paper) is not possible. One the other hand, Laplacian matrices defined by \eqref{eq:l} are relevant in many applications. For example, in~\cite{bronski2014spectral,knyazev2017signed,pan2016laplacian} problems like agent clustering are studied while the stability of networks under uncertain perturbations is considered in~\cite{chen2016semidefiniteness,zelazo2017robustness}.
\end{remark}

We identify each vertex $i$ of the graph $\cG$ with the state of an agent $x_i$. Here we are interested on scalar agents, that is $x_i\in\bbR$ for all $i=1,\ldots,m$. We now have a couple of important definitions.

\begin{definition}\leavevmode

\begin{itemize}[leftmargin=*]
	\item We say that the agents $x_i$ and $x_j$ \emph{agree} if and only if $x_i=x_j$.
	\item Consider a continuous time dynamical system defined by
	\begin{equation}\label{eq:cons}
		\dot x =f(x),
	\end{equation}
	where $x=(x_1,\ldots,x_m)\in\bbR^m$ is the vector of agents' states. Let $x(0)$ denote initial conditions and $\chi:\bbR^m\to\bbR$ be a smooth function. We say that the graph $\cG$ reaches \emph{consensus} with respect to $\chi$ if and only if all the agents agree and $x_i=\chi(x(0))$ for all $i\in\cV$.
	\item We say that $f(x)$ defines a \emph{consensus communication protocol} over $\cG$ if the solutions of \eqref{eq:cons} reach consensus.
\end{itemize}	
\end{definition}

We note that the above definition of consensus is rather general, in the sense that there can be ``discrete consensus'' if all agents only agree at discrete time points; ``finite time consensus'' if $x_i(T)=\chi(x(0))$ for all $i\in\cV$ and $t>T$ with $0\leq T<\infty$; ``asymptotic consensus'' if $\lim_{t\to\infty} x_i(t)=\chi(x(0))$ for all $i\in\cV$; and so on. Similarly, several consensus protocols can be classified with respect to the function $\chi$, see e.g. \cite{olfati2007consensus,olfati2004consensus,saber2003consensus}.

The appeal in studying consensus problems and protocols is due to their wide range of applications in, for example, computer science \cite{thomas1977majority}, formation control of autonomous vehicles \cite{fax2003information,jadbabaie2003coordination,ren2007information}, biochemistry \cite{chen2013programmable,holland2004using}, sensor networks \cite{olfati2005distributed}, social networks \cite{alves2007unveiling,xie2011social}, among many others. A simple example of consensus would be a group of people in which all agree to vote for the same candidate in an election. Another example would be a group of autonomous vehicles that are set to move with the same velocity.

In this paper we are interested in one of the simplest consensus protocols that leads to \emph{average consensus}, that is $\chi(x(0))=\frac{1}{m}\sum_{i=1}^m x_i(0)$ with the protocol defined by
\begin{equation}
	f_i(x)=\sum_{j=1}^m w_{ij} (x_j-x_i).
\end{equation}
This communication protocol is particularly interesting since it is an instance of a distributed protocol. In other words, the time evolution of $x_i$ is solely determined by its interaction with other agents directly connected to it. This type of protocols are widely investigated in engineering applications, for example to design controllers that only require local information in order to achieve their tasks \cite{lynch1996distributed,moreau2004stability,ren2008distributed,xiao2007distributed}. Alternatively, this linear average consensus protocol can be written as
\begin{equation}\label{eq:cons2}
	\dot x =-Lx,
\end{equation}
where $L$ denotes the Laplacian of $\cG$ as defined by \eqref{eq:l}. It is then clear that the behaviour of the agents is determined by the spectral properties of the Laplacian matrix \cite{merris1994laplacian,mohar1991laplacian,zelazo2014definiteness}. One of the most relevant results for systems defined by \eqref{eq:cons2} is that, if the graph is connected and all the weights are positive, then \eqref{eq:cons2} reaches average consensus asymptotically \cite{olfati2007consensus}. Although most of the scientific work has been focused on consensus protocols over unweighted graphs and with fixed topology, there is an increased interest in investigating dynamical systems defined on weighted graphs with varying and/or switching topologies \cite{casteigts2012time,mesbahi2005state,moreau2005stability,olfati2004consensus,proskurnikov2013average,tanner2007flocking}. 

In the main part of this article, Sections \ref{sec:triangle} and \ref{sec:gen}, we are going to consider linear average consensus protocols with a dynamic weight. This dynamic weight is assumed to have a slower time scale that that of the nodes. Therefore, it makes sense to approach the problem from a singular perturbation perspective. We will see that under generic conditions on the weight, the fact that the dynamics are defined on a network, induces the presence of a non-hyperbolic point. As we have described in Section \ref{sec:sfs}, one suitable technique of analysis to describe the system is then the blow-up method. Since the blow-up method is a coordinate transformation, one should check whether such a transformation preserves the network structure or not. For general networks, this is a classical problem and it is known that for certain coordinate changes, network structure is not preserved~\cite{Field,GolubitskyStewart2}. Yet, sometimes symmetries help to gain a better understanding for certain classes such as coupled cell network dynamics~\cite{NijholtRinkSanders}. As we will show, the blow-up method not only preserves the network structure for our consensus problem but the blown-up networks in different coordinate charts also have natural dynamical and network interpretations. In qualitative terms this tells us that the blow-up method is a suitable technique for the analysis of adaptive networks with multiple time scales.

Before proceeding to our main contribution, in the next section we present a first interconnection between the topics discussed above. We show that Fenichel's theory suffices to analyze state-dependent linear consensus networks with two time scales, and for which the Laplacian matrix has just a simple zero eigenvalue.

\subsection{State-dependent fast-slow consensus networks with a simple zero eigenvalue}\label{sec:app1}

\subfile{\main/subfiles/app.tex}

Next we are going to consider a case study in which Fenichel's theory is not enough to describe the dynamics of a fast-slow network.

%% file: subfiles/app.tex
In this section we show that Fenichel's theorem is enough to describe the dynamics of arbitrary fast-slow consensus networks with state dependent Laplacian as long as $\lambda_1=0$ is a simple eigenvalue, i.e., we are going to show that the zero eigenvalue corresponds to a trivial parametrized direction and that for each parameter we have a normally hyperbolic structure. The result presented below is motivated by a similar claim that appears in \cite[Section B]{awad2018time}. However, here we are not concerned with the stability of the fast nor the slow dynamics, and the use of Fenichel's theorem appears more aligned to the contents of this paper. Let us then consider the fast-slow system
\begin{equation}\label{eq:app1}
	\begin{split}
		x' &= -L(x,y,\ve)x\\
		y' &= \ve g(x,y,\ve),
	\end{split}
\end{equation}
where $x\in\bbR^m$,  $y\in\bbR^n$, $\ve>0$ is a small parameter, and $L(x,y,\ve)$ is a state dependent Laplacian matrix.

\begin{theorem} Consider \eqref{eq:app1} and a compact region $U_x\times U_y\subseteq\bbR^m\times\bbR^n$. Let $\bm 1_m:=(1,1,\ldots,1)^\top\in U_x$. If for all $(x,y)\in U_x\times U_y$ one has that $\ker L(x,y,0)=\Span\left\{ \bm 1_m\right\}$, then the set
	\begin{equation}\label{eq:app2}
		\cS_0=\left\{ (x,y)\in U_x\times U_y \, | \, x_i=\frac{1}{m}\bm 1_m^\top x(0) ~\forall i\right\}
	\end{equation}
	is a normally hyperbolic family of critical manifolds of \eqref{eq:app1}.

\end{theorem}
\begin{proof}

	Let $X=(\bar X,\hat X)\in\bbR\times\bbR^{m-1}$ be new coordinates defined by 
	\begin{equation}
	\label{eq:cchangenhyp}
	X=\begin{bmatrix}
			\bar X \\ \hat X
		\end{bmatrix}=Px=\begin{bmatrix}
			\frac{1}{m}\bm 1_m^\top \\
			Q
		\end{bmatrix}x,
	\end{equation}	
where the matrix $Q$ is found via the Gram-Schmidt process after selecting the first component
as indicated in~\eqref{eq:cchangenhyp}. Although $L(x,y,0)$ cannot really be regarded as a fixed linear operator acting on $\bbR^m$ as it depends upon $(x,y)$, the choice of the eigenvector $\bm 1_m$ is justified due to the fact that $\lambda_1=0$ is a simple zero eigenvalue of the Laplacian matrix $L(x,y,0)$ if and only if $L(x,y,0)\bm 1_m=0$ for all $(x,y)\in U_x\times U_y$. Note then that $\bar X$ denotes the average of the nodes' states. It now follows that from the equation of $x'$ we have
	\begin{equation}
		\begin{split}
			\begin{bmatrix}
				\bar X' \\
				\hat X'
			\end{bmatrix} = \begin{bmatrix}
			\frac{1}{m}\bm 1_m^\top \\
			Q
		\end{bmatrix}x' = - \begin{bmatrix}
			\frac{1}{m}\bm 1_m^\top L \bm 1_m \bar X + \frac{1}{m}\bm 1_m^\top LQ^\top \hat X\\
			QL\bm 1_m \bar X +  QLQ^\top \hat X
		\end{bmatrix}=\begin{bmatrix}
			0\\
			QLQ^\top\hat X
		\end{bmatrix}.
		\end{split}
	\end{equation}

Therefore we have that \eqref{eq:app1} is conjugate to
\begin{equation}
	\begin{split}
		\bar X' &=0\\
		\hat X' &= -\hat L(\bar X,\hat X,y,\ve) \hat X\\
		y' &= \ve \hat g(\bar X,\hat X,y,\ve),
	\end{split}
\end{equation}
where $\hat L(\bar X,\hat X,y,\ve) = QL(P^{-1}X,y,\ve)Q^\top$ and $\hat g(\bar X,\hat X,y,\ve)=g(P^{-1}X,y,\ve)$. One observes that, as expected, $\bar X$ has the role of a parameter. Furthermore, due to our hypothesis and definition of $\hat L$, we have that the matrix $\hat L(\bar X,\hat X,y,0)$ is invertible within the compact region of interest. Therefore, the corresponding critical manifold is given by $\hat \cS_0 =\left\{ \hat X=0 \right\}$. Denoting $f(\bar X,\hat X,y,\ve) =-\hat L(\bar X,\hat X,y,\ve) \hat X$ we have that $\frac{\partial f}{\partial \hat X}(\bar X,0,y,0)=-L(\bar X,0,y,0)$, which is invertible, implying that $\hat \cS_0$ is normally hyperbolic. The proof is finalized by returning to the original coordinates leading to \eqref{eq:app2}.
\end{proof}

\biblio

%% file: subfiles/triangle.tex
In this section we study a motif \cite{milo2002network}. Motifs can be seen as building blocks of more general and complex networks. Indeed, as we describe throughout this article, all the dynamic traits and properties that the triangle motif exhibits can be extended to arbitrary networks, see Section \ref{sec:gen}.

Let us consider the following network

\begin{figure}[htbp]\centering
\begin{tikzpicture}
		\node[shape=circle,draw=black, inner sep=1pt] (n1) at (0,0) {\small $1$};
		\node[shape=circle,draw=black, inner sep=1pt] (n2) at (2,0) {\small $2$};
		\node[shape=circle,draw=black, inner sep=1pt] (n3) at (1,2) {\small $3$};
		\draw[thick] (n1) -- (n2) node[midway, below] { $w$};
		\draw[thick] (n2) -- (n3) node[midway, right] { $1$};
		\draw[thick] (n1) -- (n3) node[midway, left] { $1$};
	\end{tikzpicture}
	\caption{Triangle motif: a network of three nodes connected on a cycle.}
	\label{fig:triangle}
\end{figure}
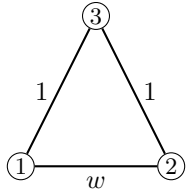

To each node $i=1,2,3$ we assign a state $x_i=x_i(t)\in\bbR$. We assume that the
dynamics of each node are defined only by diffusive coupling. Moreover, we assume
that $w\in\bbR$ is a dynamic weight depending on the vertices it connects
and on an external state $y\in\bbR$, which is assumed to have much slower time evolution
than that of the nodes. Hence, we study the fast-slow system
\begin{equation}\label{eq:sfs1}
	\begin{split}
		x' &= -L(x,y,\ve)x,\\
		y' &= \ve g(x,y,\ve),
	\end{split}\qquad\qquad
	L=\begin{bmatrix}
		w + 1 & -w & -1\\
		-w & w+1 & -1\\
		-1 & -1 & 2
	\end{bmatrix},
\end{equation}
where $w=w(x_1,x_2,y,\ve)$ is a smooth function of its arguments and
$0<\ve\ll 1$ is a small parameter. In this section we shall consider the simple
case in which $w$ is affine in the state variables, that is
\begin{equation}\label{eq:weight0}
	w = \alpha_0 + \alpha_1x_1 + \alpha_2x_2 + \alpha_3y,\\
\end{equation}
with $\alpha_0,\alpha_1,\alpha_2,\alpha_3$ real constants. We further assume
the non-degeneracy condition $\alpha_3\neq0$ to ensure coupling between the
slow and fast variables. By shifting and rescaling $y\mapsto \alpha_0+\alpha_3y$, 
 and a possible a change of signs of the variables, we may also
 assume
\begin{equation}\label{eq:weight}
	w = y + \alpha_1x_1 + \alpha_2x_2,
\end{equation}
with $\alpha_1\geq0$ and $\alpha_2\geq0$.

\subsection{Preliminary analysis (the singular limit)}\label{sec:prel_a}

The following transformation, which is simple to obtain, will be useful
throughout this work.

\begin{lemma}\label{lemma:T}
 	Consider the symmetric matrix $L$ defined in \eqref{eq:sfs1}. Then the
	orthogonal matrix
 	\begin{equation}\label{eq:T}
 		T = \begin{bmatrix}
 			\frac{\sqrt 3}{3} & -\frac{\sqrt 6}{6} & -\frac{\sqrt 2}{2}\\[1ex]
 			\frac{\sqrt 3}{3} & -\frac{\sqrt 6}{6} & \frac{\sqrt 2}{2}\\[1ex]
 			\frac{\sqrt 3}{3} & \frac{\sqrt 6}{3}  & 0
 		\end{bmatrix}
 	\end{equation}
 	diagonalizes $L$ as $D=T^\top L T=\diag\left\{ 0,3,2w+1 \right\}$.
 \end{lemma}

Thus, applying the coordinate transformation defined by $(X,Y)=(T^\top x,y)$ one obtains the conjugate diagonalized system
\begin{equation}\label{eq:diag}
	\begin{split}
		X' &= - D(X,Y)X\\
		Y' &= \ve G(X,Y,\ve),
	\end{split}
\end{equation}
where $D(X,Y)=\diag\left\{0,3,2W+1\right\}$ and
\begin{equation}
	\begin{split}
 		 W &= Y + \underbrace{\frac{\sqrt 3}{3}(\alpha_1+\alpha_2)}_{=:\beta_1}X_1 - \underbrace{\frac{\sqrt 6}{6}(\alpha_1+\alpha_2)}_{=:\beta_2}X_2 + \underbrace{\frac{\sqrt 2}{2}(\alpha_2-\alpha_1)}_{=:\beta_3}X_3\\
 		 &=Y+\beta_1 X_1 + \beta_2 X_2 + \beta_3 X_3\\
 		G(X,Y,\ve)&=g(TX,Y,\ve).
	\end{split}
\end{equation}

Observe that the fast-slow system \eqref{eq:diag} has a conserved quantity given by $X_1'=0$, which arises due the zero eigenvalue of the Laplacian matrix $L$ of
\eqref{eq:sfs1}. Since this is a trivial eigenvalue, that is, independent of the dynamics, we shall assume that $X_1$ is a coordinate on the critical
manifold and not in the fast foliation, see also Section \ref{sec:app1}.

\begin{remark}\label{rem:A}
	Due to $X_2'=-3X_2$, the set
	$\cA=\left\{ (X_1,X_2,X_3,Y)=(X_1,0,X_3,Y) \right\}$ is uniformly globally
	exponentially stable. On the other hand, the local stability properties of
	$\left\{ X_3=0 \right\}$ are dictated by the sign of $2W+1$.
\end{remark}

The previous  observations allow us to reduce the analysis of \eqref{eq:diag}
 to that of the planar fast-slow system
\begin{equation}
\begin{split}
	X_3' &= -(2\tilde W+1)X_3\\
	Y' &= \ve \tilde G(X_1,0,X_3,Y,\ve)
\end{split}, \qquad \tilde W = W|_{\cA}=Y + \beta_1 X_1 + \beta_3 X_3,
\end{equation}
where $X_1$ is regarded as a parameter. It follows that the corresponding critical manifold
is
\begin{equation}
	\tilde \cC_0 = \left\{ (X_3,Y)\in\bbR^2\, |\,  (2\tilde W+1)X_3=0  \right\}.
\end{equation}
It is now straightforward to see that, for fixed
$X_1$, $\tilde p = \left\{ (X_3,Y)=\left(0,-\frac{1}{2}-\beta_1X_1\right) \right\}$ is a
non-hyperbolic point of the critical manifold.

\begin{remark}
	Our goal will be to describe the dynamics of the network shown in Figure
	\ref{fig:triangle} as trajectories pass through the non-hyperbolic point
	$\tilde p$. The reason to
	consider this will become clear below when we give an interpretation of the
	singular dynamics in terms of the network. Thus, we assume that
	$\tilde G(X_1,0,0,-\frac{1}{2}-\beta_1X_1,0)<0$.
\end{remark}

\begin{description}[leftmargin=*]
	\item[Singular dynamics in terms of the network:] from the definition
	$X=T^\top x$ we have \[(X_1,X_2,X_3)=\left(\frac{\sqrt{3}}{3}(x_1+x_2+x_3),\, \frac{\sqrt{6}}{6}(2x_3-x_1-x_2),\, \frac{\sqrt{2}}{2}(x_2-x_1) \right).\]
	So, first of all, we have that the uniformly globally exponentially stable set
 $\cA$, previously defined by
 $\cA=\left\{ (X_1,X_2,X_3,Y)\in\bbR^4\,|\, X_2=0 \right\}$ (see Remark \ref{rem:A}), is equivalently
 given by
\begin{equation}
	\cA = \left\{ (x_1,x_2,x_3,y)\in\bbR^4\,|\, x_3=\frac{x_1+x_2}{2} \right\}.
\end{equation}
Naturally, the uniform global stability of $\cA$ is still valid. Next,
if we restrict \eqref{eq:sfs1} to $\cA$ we obtain
\begin{equation}\label{eq:tr-red1}
	\begin{split}
		x_1' &= \left(w+\half\right)(x_2-x_1)\\
		x_2' &= \left(w+\half\right)(x_1-x_2)\\
		y' &= \ve g\left(x_1,x_2,\frac{x_1+x_2}{2},y,\ve\right),
	\end{split}
\end{equation}
which is the model of a $2$-node $1$-edge fast-slow network as show in Figure
\ref{fig:reducedgraph}.

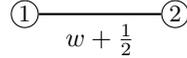
\begin{figure}[htbp]\centering
	\begin{tikzpicture}
		\node[shape=circle,draw=black, inner sep=1pt] (n1) at (0,0) {\small $1$};
		\node[shape=circle,draw=black, inner sep=1pt] (n2) at (2,0) {\small $2$};
		\draw[thick] (n1) -- (n2) node[midway, below] { $w+\frac{1}{2}$};
	\end{tikzpicture}
	\caption{Reduced graph corresponding to \eqref{eq:tr-red1}. The dynamics of the triangle motif converge
	exponentially to the dynamics of this simpler graph.}
	\label{fig:reducedgraph}
\end{figure}

Next, we note in \eqref{eq:tr-red1} that $x_1'+x_2'=0$, which implies that
$x_1(t)+x_2(t)=x_1(0)+x_2(0)=:\sigma_0$ for all $t\geq 0$. Therefore, just as in the
diagonalized system above, we can reduce the analysis of the triangle motif to
the analysis of the planar fast-slow system
\begin{equation}\label{eq:red2}
	\begin{split}
		x_1' &= \left(\half + \underbrace{y +(\alpha_1-\alpha_2)x_1+\alpha_2\sigma_0}_{=:\tilde w}\right)(-2x_1+\sigma_0)\\
		y' &= \ve g(x_1,-x_1+\sigma_0,\frac{\sigma_0}{2},y,\ve),
	\end{split}
\end{equation}

Now, it is straightforward to see that the critical manifold is given by
\begin{equation}
	\cC_0 = \left\{ (x_1,y)\in\bbR^2\, | \,  \left(\tilde w+\half\right)(-2x_1+\sigma_0) = 0\right\}.
\end{equation}
Let us consider the lines
\begin{equation}
	\begin{split}
		\cM_0 &= \left\{ (x_1,y)\in\bbR^2\, | \,  \tilde w+\half = 0\right\}\\
		\cN_0 &=  \left\{ (x_1,y)\in\bbR^2\, | \, -2x_1+\sigma_0 = 0\right\},
	\end{split}
\end{equation}
which are subsets of the critical manifold since $\cC_0=\cM_0\cup\cN_0$. It is
 clear that the intersection
 $p=\cM_0\cap\cN_0=\left\{ (x_1,y)=\left(\frac{\sigma_0}{2},-\frac{1+\sigma_0(\alpha_1+\alpha_2)}{2}\right) \right\}$
 is the only non-hyperbolic point of the layer equation of \eqref{eq:red2},
 and that the stability properties of $\cC_0$ are as shown in Figure
 \ref{fig:crit0}. For brevity let $q=-\frac{1+\sigma_0(\alpha_1+\alpha_2)}{2}$.
\begin{figure}[htbp]\centering
	\begin{tikzpicture}
	\pgftext{
	\includegraphics{\main/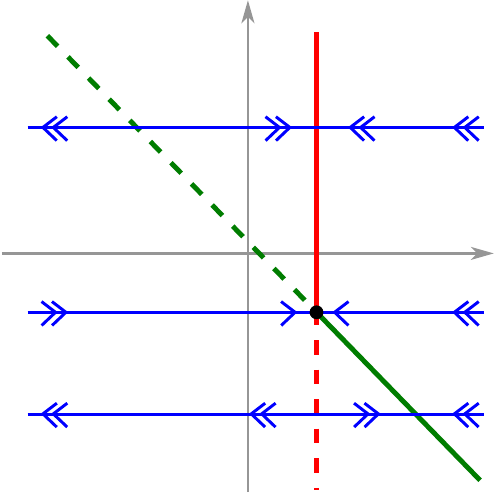}
	}
	\node at (2.65,-.075) {$x$};
	\node at (0,2.65) {$y$};
	\node[red] at (.75, 2.5) {$\cN_0^\txta$};
	\node[red] at (.75,-2.7) {$\cN_0^\txtr$};
	\node[green!50!black] at (2.4,-2.6) {$\cM_0^\txta$};
	\node[green!50!black] at (-2.2,2.4) {$\cM_0^\txtr$};
	\node[] at (.5,-.95) {$p$};
	\end{tikzpicture}\hspace{2cm}
	\begin{tikzpicture}
	\pgftext{
	\includegraphics{\main/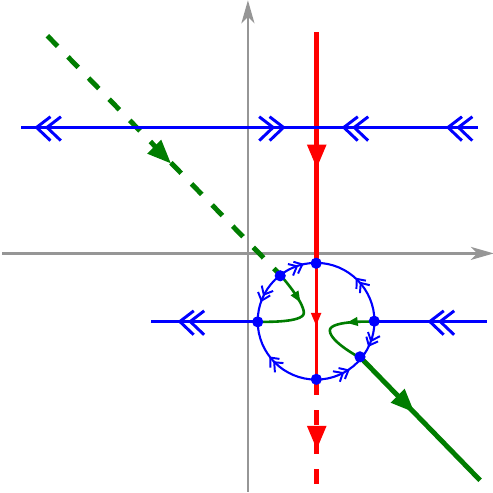}
	}
	\node at (2.65,-.075) {$x$};
	\node at (0,2.65) {$y$};
	\node[red] at (.75, 2.5) {$\cN_0^\txta$};
	\node[red] at (.75,-2.7) {$\cN_0^\txtr$};
	\node[green!50!black] at (2.4,-2.6) {$\cM_0^\txta$};
	\node[green!50!black] at (-2.2,2.4) {$\cM_0^\txtr$};
	\node[red] at (0.5,-1) {$\gamma_c$};
	\end{tikzpicture}
	\caption{Left: stability properties of the critical manifold $\cC_0$, where we
	partition the sets $\cM_0$ and $\cN_0$ into their attracting and repelling parts, and where $p=\left(\frac{\sigma_0}{2},-\frac{1+\sigma_0(\alpha_1+\alpha_2)}{2}\right)$ is a non-hyperbolic point of the fast dynamics. The case $\alpha_1-\alpha_2=0$ is degenerate and corresponds to the case where $\cM_0$ is aligned with the fast foliation. Right: blow-up of the non-hyperbolic point $p$, where $\gamma_c$ is a (singular) maximal canard. The details of the blow-up analysis are given in Section \ref{sec:blowup}.
	}
	\label{fig:crit0}
\end{figure}
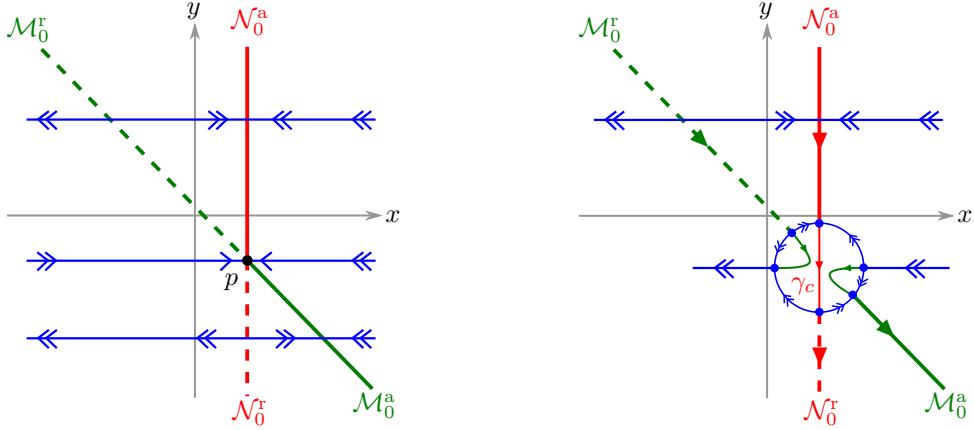

 Next, suppose trajectories converge to $\cN_0^\txta$. This means that $(x_1(t),x_2(t),x_3(t))\to \frac{\sigma_0}{2}(1,1,1)$ as $t\to\infty$. That is, the agents reach consensus, hence we call $\cN_0$ the consensus manifold. On the other hand, assume trajectories converge to $\cM_0^\txta$. For this it is necessary that $\alpha_1-\alpha_2\neq0$,
 otherwise $\cM_0$ is tangent to the fast foliation. Then $(x_1(t),x_2(t),x_3(t))\to\left( -\frac{\half+\alpha_2\sigma_0+y}{\alpha_1-\alpha_2}, \frac{\half+\alpha_1\sigma_0+y}{\alpha_1-\alpha_2},\frac{\sigma_0}{2} \right)$ as $t\to\infty$.
That is, for fixed values of $y$, agents converge to different values depending on their initial conditions. Therefore, we call $\cM_0$ the clustering manifold. Our goal will be to describe the dynamics of the network as agents transition from consensus into clustering. Thus, we also assume that $g(p,0)<0$.
\begin{remark}
	Note that the sign of $(\alpha_1-\alpha_2)$ only changes the orientation of $\cM_0$. In fact, if we denote \eqref{eq:red2} by $X(x_1,y,\ve,\rho)$ with $\rho=\alpha_1-\alpha_2\geq0$, one can show that $X(x_1,y,\ve,-\rho)=-X(-x_1,y,\ve,\rho)$. From this we shall further assume that $\alpha_1-\alpha_2\geq0$. For completeness we show the singular limit for the case $(\alpha_1-\alpha_2)<0$ in Figure \ref{fig:crit1}, but shall not be further discussed.

\begin{figure}[htbp]\centering
	\begin{tikzpicture}
	\pgftext{
	\includegraphics{\main/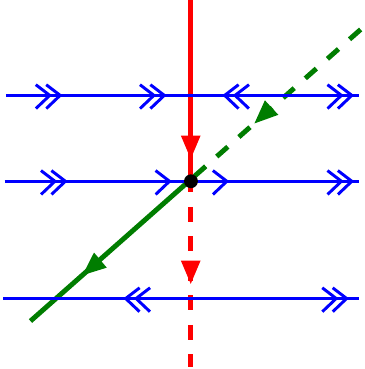}
	}
	\end{tikzpicture}\hspace{2cm}
	\begin{tikzpicture}
	\pgftext{
	\includegraphics{\main/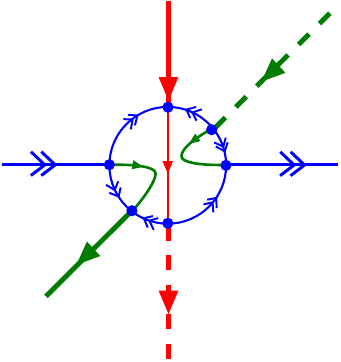}
	}
	\end{tikzpicture}
	\caption{Singular limit for the case $(\alpha_1-\alpha_2)<0$, compare with Figure \ref{fig:crit0}.}
	\label{fig:crit1}
\end{figure}

\end{remark}

\end{description}

It should be clear up to this point that the main difficulty for the analysis of the fast-slow system \eqref{eq:red2} is given by the transition across a transcritical singularity \cite{de2015planar,krupa2001extendingtrans}. Our goal is not to present a new analysis of this phenomenon but rather to study the effects of the blow-up transformation in a network. We shall show below that on each chart of the blow-up space, the resulting blow-up system can also be interpreted as a particular adaptive network. More importantly, it turns out that via the blow-up transformation one gains a clear distinction between the dynamics occurring at the different time scales. On a more technical matter, we will also show that the fact that the problem under study is defined on a graph results on a maximal canard, which in \cite{de2015planar,krupa2001extendingtrans} is non-generic.

\subsection{Main result}\label{sec:main}

Since $p$ is non-hyperbolic, the classical Fenichel theorem is not enough to conclude that for $\ve>0$ sufficiently small we have a qualitatively equivalent behaviour to the one in the limit $\ve=0$ described above. Therefore, a more detailed analysis is needed for our purposes. To state our main result, and for the analysis to be performed later, it will be convenient to move the origin of the coordinate system to the non-hyperbolic point $p$ and to relabel the coordinates of the nodes. So, let us perform the following steps

\begin{enumerate}[leftmargin=*]
	\item Relabel the fast coordinates as $x=(x_1,x_2,x_3)=(a,b,c)$. 
	This will make our notations across the blow-up charts simpler.
	\item Translate coordinates according to $(a,b,c,y)\mapsto ((a,b,c)+\frac{\sigma_0}{2}\bm 1_3,y-q)$, where $\bm 1_m=(1,\ldots,1)\in\bbR^m$ and we recall that $\sigma_0=\frac{2}{3}(a(0)+b(0)+c(0))$ and $q=-\frac{1+\sigma_0(\alpha_1+\alpha_2)}{2}$. This translation has the convenient implication $a(t)+b(t)+c(t)=0$ for all $t\geq0$.
	\item Rescale the parameter $\ve$ by $\ve\mapsto\frac{\ve}{|g(0)|}$. Thus,  we may assume that $g(0)=-1$.
\end{enumerate}

With the above we now consider
\begin{equation}\label{eq:orig2}
	\begin{split}
		\begin{bmatrix}
			a'\\ b'\\ c'
		\end{bmatrix} &=-\begin{bmatrix}
			w + 1 & -w & -1\\
			-w & w+1 & -1\\
			-1 & -1 & 2
		\end{bmatrix}\begin{bmatrix}
			a \\ b\\ c\\
		\end{bmatrix}\\
		y' &=\ve (-1+\cO(a,b,c,y,\ve))
	\end{split}
\end{equation}
where $w=-\frac{1}{2} + y + \alpha_1 a + \alpha_2 b$. Next, let us define the sections
\begin{equation}\label{eq:sections}
	\begin{split}
		\Sigma^{\txten} &= \left\{ (a,b,c,y)\in\bbR^4\,|\, y=\delta  \right\},\\[1ex]
		\Sigma^{\txtex} &= \left\{ (a,b,c,y)\in\bbR^4\,|\, y=-\delta  \right\},\\[1ex]
	\end{split}
\end{equation}
where $\delta>0$ is of order $\cO(1)$. We further define the map
\begin{equation}
	\begin{split}
		\Pi &:\Sigma^{\txten} \to	\Sigma^{\txtex},\\
	\end{split}
\end{equation}
which is induced by the flow of \eqref{eq:sfs1}. We prove the following.

\begin{theorem}\label{thm:main1} Consider the fast-slow system \eqref{eq:orig2}, where $\alpha_1-\alpha_2\geq0$. Then
\begin{itemize}[leftmargin=.5cm]
	\item[(T1)] The set $\cA=\left\{ (a,b,c,y)\in\bbR^4\,|\, c=\frac{a+b}{2}\right\}$ is globally attracting.
	\item[(T2)] The critical manifold of \eqref{eq:orig2} is contained in $\cA$ and is given by the union
	\begin{equation}
		\cC_0=\cN_0^\txta\cup\cN_0^\txtr\cup\cM_0^\txta\cup\cM_0^\txtr\cup\left\{ 0 \right\},
	\end{equation}
	where 
	\begin{equation}
		\begin{split}
			\cN_0^\txta &= \left\{ (a,b,c,y)\in\bbR^4 \, | \, a=b=c=0, \, y>0 \right\},\\
			\cN_0^\txtr &= \left\{ (a,b,c,y)\in\bbR^4 \, | \, a=b=c=0, \, y<0 \right\},\\
			\cM_0^\txta &= \left\{ (a,b,c,y)\in\bbR^4 \, | \, y+\alpha_1 a + \alpha_2 b = 0, \, y<0, \, \alpha_1-\alpha_2>0 \right\}\\
			\cM_0^\txtr &= \left\{ (a,b,c,y)\in\bbR^4 \, | \, y+\alpha_1 a + \alpha_2 b = 0, \, y>0, \, \alpha_1-\alpha_2>0 \right\}.
		\end{split}
	\end{equation}

	\item[(T3)] Restriction to $\cA$ is equivalent to the restriction to $\left\{ b=-a, \; c=0 \right\}$.

\end{itemize}

Restricted to $\cA$ and for $\ve>0$ sufficiently small:

\begin{itemize}[leftmargin=.5cm]
	\item[(T4)] There exists a slow manifold $\cN_\ve=\left\{  (a,y)\in\bbR^2\, | \, a=0 \right\}$ that is a maximal canard. Moreover, $\cN_\ve$ is attracting for $y>0$ and repelling for $y<0$.

	\item[(T5)] If $\alpha_1=\alpha_2$ then $\Pi|_{\cA}(a,b,c,y)=\Pi(a,-a,0,y)=(a,-a,0,-y)$. Moreover, every trajectory with initial condition in $\Sigma^\txten$ with $a\neq0$ diverges from $\cN_\ve$ exponentially fast as $t\to\infty$.

	\item[(T6)] If $\alpha_1-\alpha_2>0$, there exist slow manifolds $\cM_\ve^\txta$ and $\cM_\ve^\txtr$ given by
	\begin{equation}
		\begin{split}
			\cM_\ve^\txta &= \left\{ (a,y)\in\bbR^2\, | \, a=H(y,\ve)+\cO(\ve^{1/2}), \, y<0 \right\}\\
			\cM_\ve^\txtr &= \left\{ (a,y)\in\bbR^2\, | \, a=H(y,\ve)+\cO(\ve^{1/2}), \, y>0 \right\},
		\end{split}
	\end{equation}
	where 
	\begin{equation}
		H(y,\ve)=-\frac{\ve^{1/2}}{2(\alpha_1-\alpha_2)D_+(\ve^{-1/2}y)},
	\end{equation}
	with $D_+$ denoting the Dawson function \cite[pp. 219 and 235]{abramowitz1972handbook}.
	In this case, if $(a-b)|_{\Sigma^\txten}>0$ then the map $\Pi$ is well-defined and the corresponding trajectories converge towards $\cM_\ve^\txta$ as $t\to\infty$. On the contrary, if $(a-b)|_{\Sigma^\txten}<0$, then the corresponding trajectories diverge exponentially fast as $t\to\infty$.
\end{itemize}

\end{theorem}
\begin{proof}
	Items (T1) and (T2) have already been proven in our preliminary analysis of Section \ref{sec:prel_a}. Item (T3) readily follows from the relations $a+b+c=0$ and $c=\frac{a+b}{2}$, which are simultaneoulsy satisfied on $\cA$. The proof of items (T3)-(T6) is given in Section \ref{sec:trans}.
\end{proof}

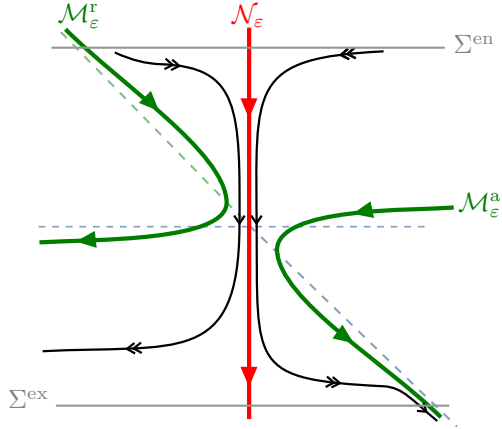
\begin{figure}[htbp]\centering
	\begin{tikzpicture}[scale=.75]
		\pgftext{
		\includegraphics[scale=1.5]{\main/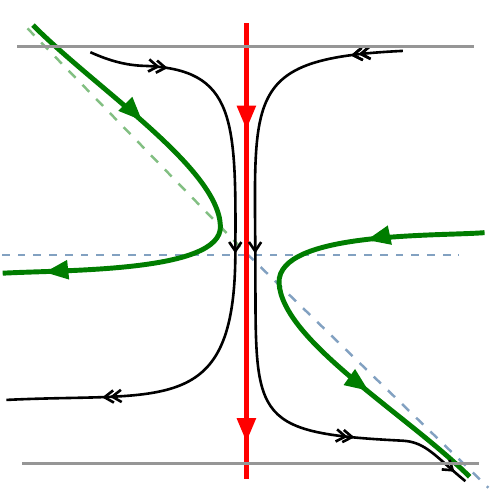}
		}
		\node[red] at (0,3.6) {$\cN_\ve$};
		\node[green!50!black] at (4.1,.25) {$\cM_\ve^\txta$};
		\node[green!50!black] at (-3,3.6) {$\cM_\ve^\txtr$};
		\node[gray] at (4.,3.1) {$\Sigma^\txten$};
		\node[gray] at (-3.9,-3.2) {$\Sigma^\txtex$};
	\end{tikzpicture}
	\caption{Schematic representation of the dynamics of \eqref{eq:orig2} in $\cA$ and for $\ve>0$ sufficiently small, compare with Figure \ref{fig:crit0}. }
	\label{fig:main1}
\end{figure}

The claims of Theorem \ref{thm:main1} are sketched in Figure \ref{fig:main1}.

\begin{description}[leftmargin=*]
	\item[Interpretation:] In terms of the network, Theorem \ref{thm:main1} tells us that: 
	\begin{itemize}[leftmargin=*]
		\item The time evolution of the node $c$ (the node that is not connected by the dynamic weight) can always be described as a combination of the dynamics of the nodes $(a,b)$ (those connected to the dynamic weight).
		\item The parameters $\alpha_1,\alpha_2$ in the definition of the weight $w=-\frac{1}{2}+y+\alpha_1 a + \alpha_2 b$, play an essential role: i) if $\alpha_1=\alpha_2$ then there is no ``clustering manifold''. Another way to interpret this degenerate case is that the nodes (or agents) have an equal contribution towards the value of the weight. This results in a zero net contribution of the nodes towards the dynamics of the weight. This is already noticeable in \eqref{eq:red2}, where $\alpha_1=\alpha_2$ results in $\tilde w$ being independent on the nodes' state. In this case the dynamics are rather simple, trajectories are attracted towards consensus for $y>0$ and repelled from consensus for $y<0$; ii) if $\alpha_1\neq\alpha_2$, then the clustering manifold exists. For suitable initial conditions, the nodes first approach consensus, but then, when $y<0$, the nodes tend towards a clustered state in which $b=-a$ and $c=0$. 
		\item The consensus manifold $\cN_\ve$ is a maximal canard, which implies that one observes a delayed loss of stability of $\cN_\ve$. In other words, one expects that trajectories exponentially near $\cN_\ve$ stay close to it for time of order $\cO(1)$ after they cross the transcritical singularity before being repelled from it. See also Appendix~\ref{sec:num}.
	\end{itemize}
\end{description}

\biblio

%% file: subfiles/blowup.tex
\renewcommand{\br}{\bar{r}}
\renewcommand{\ba}{\bar{a}}
\renewcommand{\bb}{\bar{b}}
\renewcommand{\bc}{\bar{c}}
\renewcommand{\by}{\bar{y}}
\renewcommand{\be}{\bar{\ve}}

In this section we are going to study the trajectories of \eqref{eq:orig2} in a small neighbourhood of the origin. To do this we employ the blow-up method \cite{dumortier1996canard,krupa2001extending,kuehn2015multiple,jardon2019survey}. 

\begin{remark}
	We could naturally perform the blow-up analysis restricted to the invariant and attracting subset $\cA$. However, since one of our goals is to investigate the effects of the blow-up on network dynamics, we shall proceed by blowing up \eqref{eq:orig2} and track, on each chart, the resulting ``blown-up network dynamics''.
\end{remark}

Let the blow-up map be defined by

\begin{equation}\label{eq:bu}
	a=\br\ba, \; b=\br\bb, \; c=\br\bc, \; y=\br\by, \; \ve=\br^2\be,
\end{equation}
where $\ba^2+\bb^2+\bc^2+\by^2+\be^2=1$ and $\br\geq0$.

We define the charts
\begin{equation}
\begin{split}
	K_1 = \left\{ \by = 1 \right\}, \; K_2 = \left\{ \be=1  \right\}, K_3 = \left\{ \by = -1 \right\}.
\end{split}
\end{equation}

Accordingly we define local coordinates on each chart by%

\begin{equation}
	\begin{array}{llllll}
		K_1 :  &a=r_1a_1,& b=r_1b_1,& c=r_1c_1,& y=r_1,& \ve=r_1^2\ve_1,\\
		K_2 :  &a=r_2a_2,& b=r_2b_2,& c=r_2c_2,& y=r_2y_2,& \ve=r_2^2,\\
		K_3 :  &a=r_3a_3,& b=r_3b_3,& c=r_3c_3,& y=-r_3,& \ve=r_3^2\ve_3.\\
	\end{array}
\end{equation}

The following relationship between the local blow-up coordinates will be used throughout our analysis.

\begin{lemma}\label{lemma:matching} Let $\kappa_{ij}$ denote the transformation map between charts $K_i$ and $K_j$. Then
	\begin{equation}\label{eq:matching}
		\begin{array}{llllllll}
			\kappa_{12} : &r_2 = r_1\ve_1^{1/2},\, &a_2 = \ve_1^{-1/2}a_1,\, &b_2 = \ve_1^{-1/2}b_1,\, &c_2 = \ve_1^{-1/2}c_1,\, &y_2 = \ve_1^{-1/2},\\[1ex]
			\kappa_{21} : &r_1=r_2y_2, \, &a_1=y_2^{-1}a_2,\,  &b_1=y_2^{-1}b_2,\,  &c_1=y_2^{-1}c_2,\, &\ve_1=y_2^{-2}, \\[1ex]
			\kappa_{32} : &r_2 = r_3\ve_3^{1/2},\, &a_2 = \ve_3^{-1/2}a_3,\, &b_2 = \ve_3^{-1/2}b_3,\, &c_2 = \ve_3^{-1/2}c_3,\, &y_2 = -\ve_3^{-1/2},\\[1ex]
			\kappa_{23} : &r_3=-r_2y_2, \, &a_3=-y_2^{-1}a_2,\,  &b_3=-y_2^{-1}b_2,\,  &c_3=-y_2^{-1}c_2,\, &\ve_3=y_2^{-2}. 
		\end{array}
	\end{equation}

Note that $\kappa_{ij}^{-1}=\kappa_{ji}$.

\end{lemma}

Let us now proceed with the blow-up analysis on each of the charts. We recall that on $K_1$ one studies orbits of \eqref{eq:orig2} as they approach the origin, on $K_2$ orbits within a small neighborhood of the origin, and finally on $K_3$ orbits as they leave a small neighborhood of the origin. 

\subsubsection{Analysis in the entry chart $K_1$}\label{sec:K1}
\subfile{\main/subfiles/chartK1.tex}
\subsubsection{Analysis in the rescaling chart $K_2$}\label{sec:K2}

\subfile{\main/subfiles/chartK2.tex}
\subsubsection{Analysis in the exit chart $K_3$}\label{sec:K3}
\subfile{\main/subfiles/chartK3.tex}
\subsubsection{Full transition and proof of main result}\label{sec:trans}
\subfile{\main/subfiles/transition.tex}

\biblio

%% file: subfiles/chartK1.tex
\renewcommand{\br}{r_1}
\renewcommand{\by}{y_1}
\renewcommand{\ba}{a_1}
\renewcommand{\bb}{b_1}
\renewcommand{\bc}{c_1}
\renewcommand{\be}{\ve_1}
\newcommand{\bw}{\bar w}
\newcommand{\bA}{A_1}
\newcommand{\bB}{B_1}
\newcommand{\bC}{C_1}
\newcommand{\bW}{\bar W}

	In this chart the blow up map is given by
\begin{equation}\label{eq:bu1}
	a=\br\ba, \; b=\br\bb, \; c=\br\bc,\; y = \br,\; \ve=\br^2\be.
\end{equation}

We then obtain the blown up vector field
\begin{equation}\label{eq:K1_0}
	\begin{split}
		\begin{bmatrix}
			\ba' \\ \bb' \\ \bc'
		\end{bmatrix} &=
		-\begin{bmatrix}
			\half & \half & -1 \\
			\half & \half & -1\\
			-1 & -1 & 2
		\end{bmatrix}\begin{bmatrix}
			\ba\\ \bb \\ \bc
		\end{bmatrix}+r_1f_1(\ba,\bb,\bc,\br,\be)\\
		\br' &= \br^2\be(-1+\cO(\br))\\
		\be' &= -2\br\be^2(-1+\cO(\br)),
	\end{split}
\end{equation}
where $f_1(\ba,\bb,\bc,\br,\be)$ reads as
\begin{equation}\label{eq:f1}
	f_1 = \left( (1+\alpha_1\ba+\alpha_2\bb) \begin{bmatrix}
		-1 & 1 & 0 \\
		1 & -1 & 0\\
		0 & 0 & 0
	\end{bmatrix}  - \be(-1+\cO(\br)) \begin{bmatrix}
		1 & 0 & 0\\
		0 & 1 & 0\\
		0 & 0 & 1
	\end{bmatrix} \right)\begin{bmatrix}
		\ba\\ \bb \\ \bc
	\end{bmatrix}.
\end{equation}

\begin{remark}
	In \eqref{eq:K1_0} and \eqref{eq:f1} the term $-1+\cO(\br)$ stands for exactly the same function.
\end{remark}

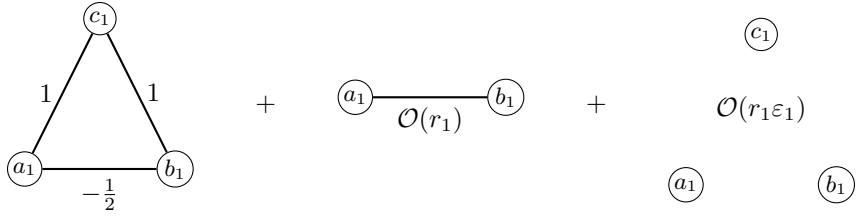
\begin{figure}[htbp]\centering
	\begin{tikzpicture}[baseline=(current bounding box.center)]
		\node[shape=circle,draw=black, inner sep=1pt] (n1) at (0,0) {\small $\ba$};
		\node[shape=circle,draw=black, inner sep=1pt] (n2) at (2,0) {\small $\bb$};
		\node[shape=circle,draw=black, inner sep=1pt] (n3) at (1,2) {\small $\bc$};
		\draw[thick] (n1) -- (n2) node[midway, below] { $-\half$};
		\draw[thick] (n2) -- (n3) node[midway, right] { $1$};
		\draw[thick] (n1) -- (n3) node[midway, left] { $1$};
	\end{tikzpicture}
	$ \qquad+ \qquad$
	\begin{tikzpicture}[baseline=(current bounding box.center)]
		\node[shape=circle,draw=black, inner sep=1pt] (n1) at (0,0) {\small $\ba$};
		\node[shape=circle,draw=black, inner sep=1pt] (n2) at (2,0) {\small $\bb$};
		\draw[thick] (n1) -- (n2) node[midway, below] { $\cO(\br)$};
	\end{tikzpicture}
	$ \qquad+ \qquad$
	\begin{tikzpicture}[baseline=(current bounding box.center)]
		\node[shape=circle,draw=black, inner sep=1pt] (n1) at (0,0) {\small $\ba$};
		\node[shape=circle,draw=black, inner sep=1pt] (n2) at (2,0) {\small $\bb$};
		\node[shape=circle,draw=black, inner sep=1pt] (n3) at (1,2) {\small $\bc$};
		\node at (1,1) {$\cO(\br\be)$};
	\end{tikzpicture}
\caption{Network interpretation corresponding to \eqref{eq:K1_0}. The order $\cO(1)$ terms in \eqref{eq:K1_0} correspond to a triangle motif with fixed weights. The particular values of the weights make such a network degenerate in the sense that the corresponding Laplacian has a kernel of dimension two. Next, the order $\cO(\br)$ terms in \eqref{eq:K1_0} correspond to two nodes connected by a dynamic weight. Finally, the $\cO(\br,\be)$ correspond to internal node dynamics.}
	\label{fig:K10}
\end{figure}

Let us interpret the equations in the first chart from a network dynamics perspective. We are interested in the dynamics of \eqref{eq:K1_0} for $\be$ small and with $\br\to0$. This is because $\br\to0$ is equivalent to $y(t)$ approaching the origin in \eqref{eq:orig2}. Thus, we may regard \eqref{eq:K1_0} as a perturbation of a network with fixed weights as shown in Figure \ref{fig:K10}.

The order $\cO(\br)$ terms can be seen as a smaller network, only involving the nodes $(\ba,\bb)$ and with dynamic edge with weight $\br(1+\alpha_1\ba+\alpha_2\bb)$. The  order $\cO(\br\be)$ can be interpreted as internal dynamics on each node.

Continuing with the analysis, it is straightforward to check (with the help of \eqref{eq:T}) that for $\br=0$ we have $\bc(t_1)\to\frac{\ba(t_1)+\bb(t_1)}{2}$ as $t_1\to\infty$, where $t_1$ denotes the time parameter of \eqref{eq:K1_0}.  We now proceed with a more detailed analysis of \eqref{eq:K1_0} as follows.

\begin{proposition}\label{prop:K1eqs}
	System \eqref{eq:K1_0} has the following sets of equilibrium points.
	\begin{equation}
		\begin{split}
			S_{1,0} &= \left\{ (\ba,\bb,\bc,\br,\be)\in\bbR^5\,|\, \br = 0, \, \bc = \frac{\ba+\bb}{2} \right\},\\
			N_{1,0} &= \left\{ (\ba,\bb,\bc,\br,\be)\in\bbR^5\,|\, \be = 0, \, \bc = \frac{\ba+\bb}{2}, \, \ba=\bb \right\},\\
			M_{1,0} &= \left\{ (\ba,\bb,\bc,\br,\be)\in\bbR^5\,|\, \be = 0, \, \bc = \frac{\ba+\bb}{2}, \, 1+\alpha_1\ba+\alpha_2\bb=0 \right\}.\\
		\end{split}
	\end{equation}
\end{proposition}
\begin{proof}
	Straightforward computations.
\end{proof}

Next, we show that the set defined by $\bc = \frac{\ba+\bb}{2}$ is an attracting centre manifold.

\begin{proposition}\label{prop1K1}
	The system given by \eqref{eq:K1_0} has a local $4$-dimensional centre manifold $\cW_1^\txtc$ and a local $1$-dimensional stable manifold $\cW_1^\txts$. The centre manifold $\cW_1^\txtc$ contains the sets of Proposition \ref{prop:K1eqs}. Furthermore, $\cW_1^\txtc$ is given by $\bc=\frac{\ba+\bb}{2}$, and the flow along it reads as
	\begin{equation}\label{eq:K1p1}
		\begin{split}
			\ba' &= \br(1+\alpha_1\ba+\alpha_2\bb)(\bb-\ba) - \br\be(-1+\cO(\br))\ba\\
			\bb' &= \br(1+\alpha_1\ba+\alpha_2\bb)(\ba-\bb) - \br\be(-1+\cO(\br))\bb\\
			\br' &= \br^2\be(-1+\cO(\br))\\
			\be' &= -2\br\be^2(-1+\cO(\br)).
		\end{split}
	\end{equation}

\end{proposition}
\begin{proof}
We start by using the similarity transformation~$\begin{bmatrix}
	\bA & \bB & \bC
\end{bmatrix}^\top=T^\top \begin{bmatrix}
	\ba & \bb & \bc
\end{bmatrix}^\top$, where $T$ is defined in \eqref{eq:T}. Under such a transformation one rewrites  \eqref{eq:K1_0} as
	 \begin{equation}\label{eq:K1_diag}
		\begin{split}
		\begin{bmatrix}
			\bA' \\ \bB' \\ \bC'
		\end{bmatrix} &= - \begin{bmatrix}
			0 & 0 & 0 \\
			0 & 3 & 0 \\
			0 & 0 & 0
		\end{bmatrix} \begin{bmatrix}
			\bA\\ \bB \\ \bC
		\end{bmatrix} + \br F_1(\bA,\bB,\bC,\br,\be)\\
		\br' &= \br^2\be(-1+\cO(\br))\\
		\be' &= -2\br\be^2(-1+\cO(\br)),
	\end{split}
	 \end{equation}
	 where
	 \begin{equation}
		 \begin{split}
			 F_1 = &\left(\underbrace{\left( 1+\frac{\sqrt 3}{3}(\alpha_1+\alpha_2)\bA-\frac{\sqrt 6}{6}(\alpha_1+\alpha_2)\bB+\frac{\sqrt 2}{2}(\alpha_2-\alpha_1)\bC \right)}_{=:\bW}\begin{bmatrix}
 				0 & 0 & 0 \\
 				0 & 0 & 0\\
 				0 & 0 & -2
 			\end{bmatrix} - \right.\\
			&\left. \be(-1+\cO(\br)) \begin{bmatrix}
 				1 & 0 & 0\\
 				0 & 1 & 0\\
 				0 & 0 & 1
 			\end{bmatrix}\right)\begin{bmatrix}
 				\bA\\ \bB \\ \bC
 			\end{bmatrix}.
		 \end{split}
	 \end{equation}

It is now straightforward to see that there is a $1$-dimensional stable manifold $\cW_1^s$ tangent to the $B_1$-axis and a $4$-dimensional centre manifold $\cW_1^c$ containing the set of equilibrium points $\left\{ (\bA,\bB,\bC,\br,\be)\in\bbR^5\,|\, \br=\bB=0 \right\}$.

	 \begin{remark}
	 	Observe that, due to the term $(-1+\cO(\br))$, the vector field corresponding to $\bB'$ is not decoupled from the center directions. However, we show below that $\cW_1^c$ is indeed given by $\bB=0$.
	 \end{remark}

The centre manifold $\cW_1^c$ can be expressed by $\bB=h_1(\br,\bA,\bC,\be)$ satisfying $h_1(0)=0$, $\textnormal{D}h_1(0)=0$, where $\textnormal{D}h_1$ denotes the Jacobian of $h_1$. Let $h_1$ be given as
	 \begin{equation}\label{eq:h1}
	 	h_1 = \sum_{\substack{i,j,k,l\geq0\\i+j+k+l\geq2}}\sigma_{ijkl}\bA^i\bC^j\br^k\be^l,
	 \end{equation}
	 where $\sigma_{ijkl}$ denotes scalar coefficients. Substituting \eqref{eq:h1} into the equation for $\bB'$ we get
	 \begin{equation}\label{eq:K1-cm1}
	 	-3h_1  = \br\be(-1+\cO(\br))\left( -h_1 + \bA\frac{\partial h_1}{\partial\bA}+\bC\frac{\partial h_1}{\partial\bC}+\br\frac{\partial h_1}{\partial\br}-2\be\frac{\partial h_1}{\partial\be} \right)-2\br\bW\bC\frac{\partial h_1}{\partial\bC}.
	 \end{equation}

We now have the following observations:
\begin{enumerate}[leftmargin=*]
	\item All the monomials in the right hand side of \eqref{eq:K1-cm1} are of degree at least $3$, therefore, all coefficients $\sigma_{ijkl}$ with $i+j+k+l=2$ are zero.
	\item Since the right hand side of \eqref{eq:K1-cm1} is of order $\cO(\br)$ we have that all coefficients $\sigma_{ij0l}$ are zero for all $i+j+l\geq3$. Naturally, we then have that $h_1\in \cO(\br)$ and thus $k\geq1$.
	\item The coefficients $\sigma_{ijk0}$, $k\geq1$, are computed from the equality
	\begin{equation}\label{eq:K1-ser1}
		\begin{split}
			-3h_1 &= -2\br\bW\bC\frac{\partial h_1}{\partial\bC}\\
			&=-2\br\left( 1+\frac{\sqrt 3}{3}(\alpha_1+\alpha_2)\bA-\frac{\sqrt 6}{6}(\alpha_1+\alpha_2)h_1+\frac{\sqrt 2}{2}(\alpha_2-\alpha_1)\bC \right)\bC\frac{\partial h_1}{\partial\bC}\\
			&=-2\br(1+\eta_1\bA - \eta_2 h_1 + \eta_3\bC)\bC\frac{\partial h_1}{\partial\bC},
		\end{split}
	\end{equation}
	where the last equation is introduced for simplicity. We readily see that all coefficients $\sigma_{i0k0}$ with $i+k\geq3$ are zero. Next, for $i+j+k=3$, the term $h_1\bC\frac{\partial h_1}{\partial\bC}$ does not play a role because its degree is at least $4$. It follows from the first item that $\sigma_{ijk0}=0$ for $i+j+k=3$. Next, let us write \eqref{eq:K1-ser1} in a simplified form by i) expanding it, ii) writing all monomials in the exact same form $\bA^i\bC^j\br^k$, iii) by omitting the monomial, and iv) omitting the $0$ of the superscript $\sigma_{ijk0}$. We get
	\begin{equation}
		\begin{split}
			-3\sum\sigma_{ijk} &= -2\sum j\sigma_{ij(k-1)} -2\eta_1\sum j \sigma_{(i-1)j(k-1)}-2\eta_3\sum j\sigma_{i(j-1)(k-1)}\\
			&+2\eta_2\sum\sigma_{ijk}\sum j\sigma_{ij(k-1)}
		\end{split}
	\end{equation}
Now, it suffices to note that \emph{for each monomial}, the coefficient $\sigma_{ijk}$, with $i+j+k=n$ and $n>3$, of the left-hand side depends exclusively on coefficients $\sigma_{ijk}$ with $i+j+k<n$. From the previous items, and by progressing at each degree $n$, it follows that $\sigma_{ijk0}=0$ for all $i+j+k\geq2$.
\item The exact same argument as in item 3. applies for $l\geq1$.
\end{enumerate}

The expression of the centre manifold in the original coordinates is obtained by noting that $\bB=\frac{\sqrt 6}{6}(2\bc-\ba-\bb)$, implying that $\bc = \frac{\ba+\bb}{2}$ as stated. Finally, the flow along the centre manifold is obtained by taking into account the restriction $\bc = \frac{\ba+\bb}{2}$.
\end{proof}

\begin{remark}
	$\cW_1^\txtc$ is the blow-up of $\cA$.
\end{remark}

Before proceeding with the analysis on $\cW_1^\txtc$, we have the next observation.

\begin{proposition}\label{prop:K1sign}
	Let $(\ba(0),\bb(0),\bc(0))$ denote initial conditions of \eqref{eq:K1_0} and let $\be=0$. Then $\sign(\ba)\to\sign(\ba(0)-\bb(0))$ as $t\to\infty$.
\end{proposition}
\begin{proof}
	It is easier to see the claim in \eqref{eq:K1_diag} with $\be=0$, and where $\bA=\frac{\sqrt{3}}{3}(\ba+\bb+\bc)=0$, $\sign(\bB)=\sign(2\bc-\ba-\bb)$ and $\sign(\bC)=\sign(\bb-\ba)$. We note in \eqref{eq:K1_diag} that $\sign(\bB)$ and  $\sign(\bC)$ are invariant. Therefore as $\bB\to0$ (equivalently $2\bc\to\ba+\bb$) we have $\ba+\bb\to0$ and therefore $\sign(\bb-\ba)\to\sign(-2\ba)$ from which the claim immediately follows.
\end{proof}

The previous observation is important since, as we will see, in $\cW_1^\txtc$ the set $\left\{ \ba=0\right\}$ is invariant.  Note that we can now desingularize the dynamics restricted to $\cW_1^\txtc$ by dividing by $\br$ in \eqref{eq:K1p1}, as is usually the case when blowing up, to obtain
\begin{equation}\label{eq:K1_1}
	\begin{split}
		\ba' &= (1+\alpha_1\ba+\alpha_2\bb)(\bb-\ba) - \be(-1+\cO(\br))\ba\\
		\bb' &= (1+\alpha_1\ba+\alpha_2\bb)(\ba-\bb) - \be(-1+\cO(\br))\bb\\
		\br' &= \br\be(-1+\cO(\br))\\
		\be' &= -2\be^2(-1+\cO(\br)).
	\end{split}
\end{equation}

\begin{remark}\label{rem:ic}
	Recall that $\ba(t_1)+\bb(t_1)+\bc(t_1)=0$ for all $t_1\geq0$. Moreover, since in $\cW_1^c$ we have $\bc=\frac{\ba+\bb}{2}$ we further have $\ba(t_1)+\bb(t_1)=0$ for all $t_1\geq0$. Therefore, we can consider instead of \eqref{eq:K1_1} the $3$-dimensional system
	\begin{equation}\label{eq:K1-3d}
		\begin{split}
			\ba' &= -2(1+(\alpha_1-\alpha_2)\ba)\ba - \be(-1+\cO(\br))\ba\\
			\br' &= \br\be(-1+\cO(\br))\\
			\be' &= -2\be^2(-1+\cO(\br)).
	\end{split}
	\end{equation}

	Naturally, solutions of \eqref{eq:K1-3d} give solutions of \eqref{eq:K1_1} by adding $\bb(t_1)=-\ba(t_1)$. Therefore we proceed by studying \eqref{eq:K1-3d}. It is worth noting that on $\cW_1^\txtc$, the set $\left\{ (\ba,\br,\be)\in\bbR^3\,|\,\ba=0\right\}$ is invariant. Therefore, it is important to keep track of the sign of $\ba$ as it approaches  $\cW_1^\txtc$. Such sign is given by Proposition \ref{prop:K1sign}. That is, if $\ba(0)-\bb(0)>0$ (resp. $\ba(0)-\bb(0)<0$), then $\ba>0$ (resp. $\ba<0$) on $\cW_1^\txtc$. Similarly, if $\ba(0)-\bb(0)=0$, then $\ba=0$ on $\cW_1^\txtc$. Finally, we recall that $\cW_1^\txtc$ coincides precisely with the invariant set $\cA$ written in the coordinates of this chart (see the statement of Theorem \ref{thm:main1}). 

\end{remark}

To study the dynamics in this chart, we are going to be interested in the properties of the flow between the sections
\begin{equation}
	\begin{split}
		\Delta_1^\txten &= \left\{ (\ba,\br,\be)\in\bbR^3\, | \, \br=\delta_1,\, \be<\mu_1 \right\}\\
		\Delta_1^\txtex &= \left\{ (\ba,\br,\be)\in\bbR^3\, | \, \br<\delta_1,\, \be=\mu_1 \right\},
	\end{split}
\end{equation}
where $\delta_1>0$, and $\mu_1>0$ is sufficienlty small. The precise meaning of these sections becomes clear in Section \ref{sec:trans} where we compute a transition map through al whole neighborhood of the origin of\eqref{eq:orig2}. For now it shall be enough to mention that the definition $\Delta_1^\txten$ is motivated by the entry section $\Sigma^\txten$ (recall \eqref{eq:sections}), while $\Delta_1^\txtex$ is a convenient section allowing us to transition towards the central chart $K_2$.

We observe that the subspaces $\left\{ \ba=0\right\}$, $\left\{ \br=\be=0\right\}$, $\left\{ \br=0\right\}$, and $\left\{ \be=0\right\}$ are all invariant and thus are helpful to describe the overall dynamics \eqref{eq:K1-3d}. So, we proceed as follows.

\begin{description}[leftmargin=*]

\item[In $\left\{ \ba=0 \right\}$] we have the planar system
		\begin{equation}\label{eq:K1a0}
		\begin{split}
				\br' &= \br\be(-1+\cO(\br))\\
				\be' &= -2\be^2(-1+\cO(\br)).
	\end{split}
	\end{equation}
which has a line of zeros $(\br,\be)=(\br,0)$ and an unstable invariant manifold $(\br,\be)=(0,\be)$. Note that away from $\left\{\be=0\right\}$ the flow of \eqref{eq:K1a0} is equivalent to that of a planar saddle. Next, we want to compute the time it takes to travel from $\Delta_1^\txten$ to $\Delta_1^\txtex$. Therefore, assume initial conditions $(\br,\be)=(\delta_1,\be^*)$ and boundary conditions $(\br,\be)=(\br(T_1),\mu_1)$. From \eqref{eq:K1a0} we find that
\begin{equation}
	\br(T_1) = \delta_1\left( \frac{\be^*}{\mu_1}\right)^{1/2}.
\end{equation}

Then, one can estimate the transition time $T_1$ by integrating the equation for $\be'$, so that we get
\begin{equation}\label{eq:T1}
	T_1=\frac{1}{2}\left( \frac{1}{\be^*}-\frac{1}{\mu_1} \right)(1+\cO(\delta_1)), \qquad 0<\be^*\leq \mu_1.
\end{equation}

	\item[In $\left\{ \br=\be=0 \right\}$] we have the $1$-dimensional system
	\begin{equation}\label{eq:K1_00}
		\ba' = -2(1+(\alpha_1-\alpha_2)\ba)\ba,
	\end{equation}
	where we recall that $\alpha_1-\alpha_2\geq0$. In this case we have, generically, two hyperbolic equilibrium points: one stable at $\ba=0$ and one unstable at $\ba=\frac{1}{\alpha_2-\alpha_1}$. If $\alpha_1=\alpha_2$ then only the stable equilibrium at the origin exists. It will be useful to integrate \eqref{eq:K1_00}, that is
	\begin{equation}\label{eq:K1-sola}
		\ba(t_1)=-\frac{\ba^*}{(\alpha_1-\alpha_2)\ba^* - ((\alpha_1-\alpha_2)\ba^*+1)\exp(2t_1)},
	\end{equation}
	where $\ba^*$ denotes an initial condition for \eqref{eq:K1_00}. 

\item[In $\left\{ \br=0 \right\}$] we have
\begin{equation}\label{eq:K1_1p}
		\begin{split}
		\ba' &= -2(1+(\alpha_1-\alpha_2)\ba)\ba + \be\ba\\
		\be' &= 2\be^2.
	\end{split}
	\end{equation}
	Then, we have two $1$-dimensional centre manifolds: $\cE_1^\txta$ is a centre manifold to the equilibrium point $(\ba,\be)=(0,0)$ and  $\cE_1^\txtr$ to $(\ba,\be)=\left(\frac{1}{\alpha_2-\alpha_1},0\right)$. The flow on both centre manifolds is given by $\be' = 2\be^2$, and we have that $\cE_1^\txta$ is tangent to the $\be$-axis while $\cE_1^\txtr$ is tangent to the vector $\begin{bmatrix} 1 & -2(\alpha_2-\alpha_1)\end{bmatrix}^\top$. In fact, one can show that $\cE_1^\txta$ is actually given by $\ba=0$ and that it is unique. On the other hand $\cE_1^\txtr$ is not unique and has the expansion $\ba=\frac{1}{\alpha_2-\alpha_1}+\frac{1}{2(\alpha_1-\alpha_2)}\be+\cO(\be^2)$. Since $\be\geq 0$, we have that in a small neighbourhood of  $(\ba,\be)=(0,0)$ the flow is equivalent to that of a saddle, while in a small neighbourhood of $(\ba,\be)=(\frac{1}{\alpha_2-\alpha_1},0)$ the flow is equivalent to that of a source. From this analysis we conclude that the flow of \eqref{eq:K1_1p} is as sketched in Figure \ref{fig:K1r0}.
	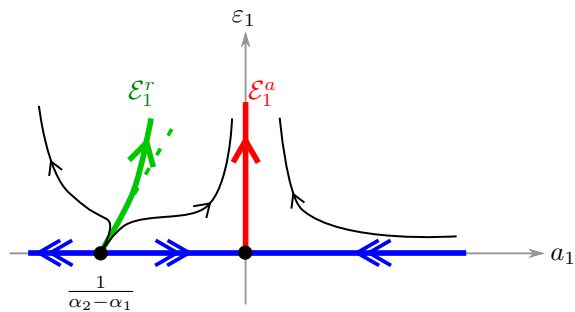
\begin{figure}[htbp]\centering
		\begin{tikzpicture}
			\node at (0,0){
			\includegraphics{\main/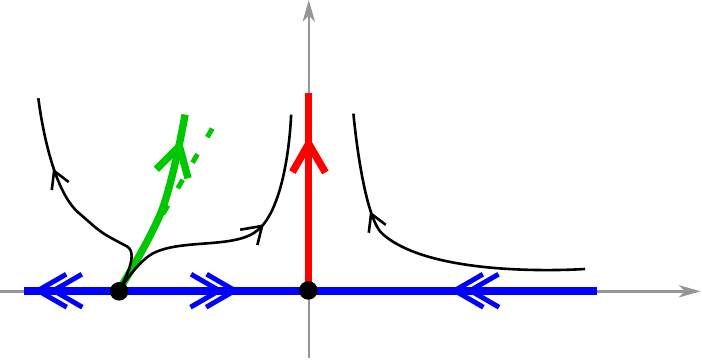}
			};
			\node at (3.8,-1.17) {$\ba$};
			\node at (-.45,2) {$\be$};
			\node[green!50!black] at (-1.8,1) {$\cE_1^r$};
			\node[red]  at (-.2,1) {$\cE_1^a$};
			\node at (-2.35,-1.65) {$\frac{1}{\alpha_2-\alpha_1}$};
		\end{tikzpicture}
		\caption{Schematic of the flow of \eqref{eq:K1_1p} for $\alpha_2-\alpha_1>0$ and $\be$ small.  }
		\label{fig:K1r0}
	\end{figure}

\begin{remark}
	The orbit $\cE_1^\txtr$ can be identified with the critical manifold $\cM_0^\txtr$ as it goes up on the blow-up sphere. The same correspondence holds for $\cE_1^\txta$ and $\cN_0$. Compare Figures \ref{fig:K1r0} and \ref{fig:crit0}.
\end{remark}

\item[In $\left\{ \be=0 \right\}$] we have
\begin{equation}\label{eq:K1e0}
	\begin{split}
		\ba' &= -2(1+(\alpha_1-\alpha_2)\ba)\ba\\ 
		\br' &= 0.
	\end{split}
\end{equation}
Therefore, the $(\ba,\br)$-plane is foliated by lines parallel to the $\br$-axis. Along each leaf the flow is given by \eqref{eq:K1_00}.

\end{description}

We can now summarize the previous analysis in the following Proposition, which completely characterizes the dynamics of \eqref{eq:K1_0}.

\begin{proposition}The following statements hold for \eqref{eq:K1_0}.
\begin{enumerate}[leftmargin=*]
	\item There exist a $1$-dimensional local stable manifold $\cW_1^\txts$ and a $4$-dimensional local centre-stable manifold  $\cW_1^\txtc$, which is given by the graph of $\bc=\frac{\ba+\bb}{2}$.
	
	\end{enumerate}
	Restricted to $\cW_1^\txtc$ one has $\bb=-\ba$, which implies $\bc=0$, and:
	\begin{enumerate}[leftmargin=*]
		\setcounter{enumi}{1}
		\item There is an attracting $2$-dimensional centre manifold $\cC_1^\txta$. The manifold $\cC_1^\txta$ contains a line of zeros $\ell_1^\txta=\left\{(\br,\ba,\be)\in\bbR^3\, | \, \ba=\be=0 \right\}$ and a $1$-dimensional centre manifold \\
		$\cE_1^\txta=\left\{(\br,\ba,\be)\in\bbR^3\, | \, \br=\ba=0\right\}$. On the plane $\left\{ \br=0 \right\}$, the centre manifold $\cE_1^\txta$ is unique. The flow along $\cE_1^\txta$ is unstable, that is, it diverges from the origin; while the flow on $\cC_1^\txta$ away from $\ell_1^\txta$ is locally equivalent to that of a saddle.
		\item There is a repelling $2$-dimensional centre manifold $\cC_1^\txtr$. The manifold $\cC_1^\txtr$ contains a line of zeros $\ell_1^\txtr=\left\{(\br,\ba,\be)\in\bbR^3\, | \, \ba=\frac{1}{\alpha_2-\alpha_1},\,\be=0 \right\}$ and a $1$-dimensional centre manifold $\cE_1^\txtr=\left\{(\br,\ba,\be)\in\bbR^3\, | \, \br=0, \ba=\frac{1}{\alpha_2-\alpha_1}+\cO(\be)\right\}$. The flow along $\cE_1^\txtr$ is unstable, that is, it diverges from the equilibrium point $(\ba,\be)=\left(\frac{1}{\alpha_2-\alpha_1},0\right)$; while the flow on $\cC_1^\txtr$ away from $\ell_1^\txtr$ is locally equivalent to that of a saddle.
	\end{enumerate}
\end{proposition}

\begin{proof}
The existence, graph representation and dimension of $\cW_1^\txtc$ is already proven in Proposition \ref{prop1K1}. The existence and dimension of $\cC_1^\txta$ and of $\cC_1^\txtr$	follow from the linearization of \eqref{eq:K1-3d}. The flow on $\cC_1^\txta$ and on $\cC_1^\txtr$ follow from \eqref{eq:K1-3d} by noting that, up to leading order terms, the vector field restricted to either of the centre manifolds is given by
\begin{equation*}
\pushQED{\qed} 
	\begin{split}
		\br' &= -\br\be\\
		\be' &= 2\be^2.\qedhere
	\end{split}
	\popQED
\end{equation*}
\renewcommand{\qedsymbol}{}
\end{proof}

We are now ready to describe the flow of \eqref{eq:K1-3d}. Let $\Pi_1:\Delta_1^\txten\to\Delta_1^\txtex$ be a map defined by the flow of \eqref{eq:K1-3d}. 

	\begin{theorem} The image $\Pi_1(\Delta_1^\txten)$ in $\Delta_1^\txtex$ is of the form
	\begin{equation}
		\Pi_1\left( \begin{matrix}
	  \ba \\[1ex] \delta_1 \\[1ex] \be
		\end{matrix}\right)=\left( \begin{matrix}
			h_{\ba} + O (\ba\be^2) \\[1ex]
			\delta_1\left( \frac{\be}{\mu_1} \right)^{1/2}\\[1ex]
			\mu_1
		\end{matrix} \right),
	\end{equation}
	where the function $h_{\ba}=h_{\ba}(\ba,\delta_1,\be)$ is given by
	\begin{equation}\label{eq:K1ha}
		h_{\ba} = -\frac{\ba}{(\alpha_1-\alpha_2)\ba - ((\alpha_1-\alpha_2)\ba+1)\exp(2T_1)}, \qquad T_1=\frac{1}{2}\left( \frac{1}{\be}-\frac{1}{\mu_1} \right)(1+\cO(\delta_1)).
	\end{equation}

	\end{theorem}
	\begin{proof}
		The proof follows our previous analysis. The term $h_{\ba}$ is obtained from \eqref{eq:K1-sola} and evaluating the transition time \eqref{eq:T1}. The higher order terms $\cO(\ba\be^2)$ follow from \eqref{eq:K1_1p} with $\be>0$ small. For the expression of $h_{\ba}$ it is important to recall Proposition \ref{prop:K1sign}. This means that the initial condition $\ba^*$ in \eqref{eq:K1-sola} has the same sign as $\ba(0)-\bb(0)$, and where $\ba(0),\bb(0)$ are initial conditions of \eqref{eq:K1_0}. 
	\end{proof}

\begin{remark}\leavevmode
\begin{itemize}[leftmargin=*]
	\item If $\alpha_1-\alpha_2=0$ then $h_{\ba}= \ba\exp(-2T_1)$.
	\item If $\ba>\frac{1}{\alpha_1-\alpha_2}$, then the function $h_{\ba}$ is well-defined for any point $(\ba,\delta_1,\be)\in\Delta_1^\txten$. If $\ba\leq\frac{1}{\alpha_1-\alpha_2}$, then $h_{\ba}$ is well-defined only for $T_1<\frac{1}{2}\ln\left( \frac{(\alpha_1-\alpha_2)\ba}{(\alpha_1-\alpha_2)\ba+1} \right)$. In such a case we choose suitably $0<\be<\mu_1\ll\delta_1$ so that the function $h_{\ba}$ is well-defined. 
\end{itemize}
	
\end{remark}

\begin{figure}[htbp]\centering
	\begin{tikzpicture}
		\pgftext{
		\includegraphics[scale=.75]{\main/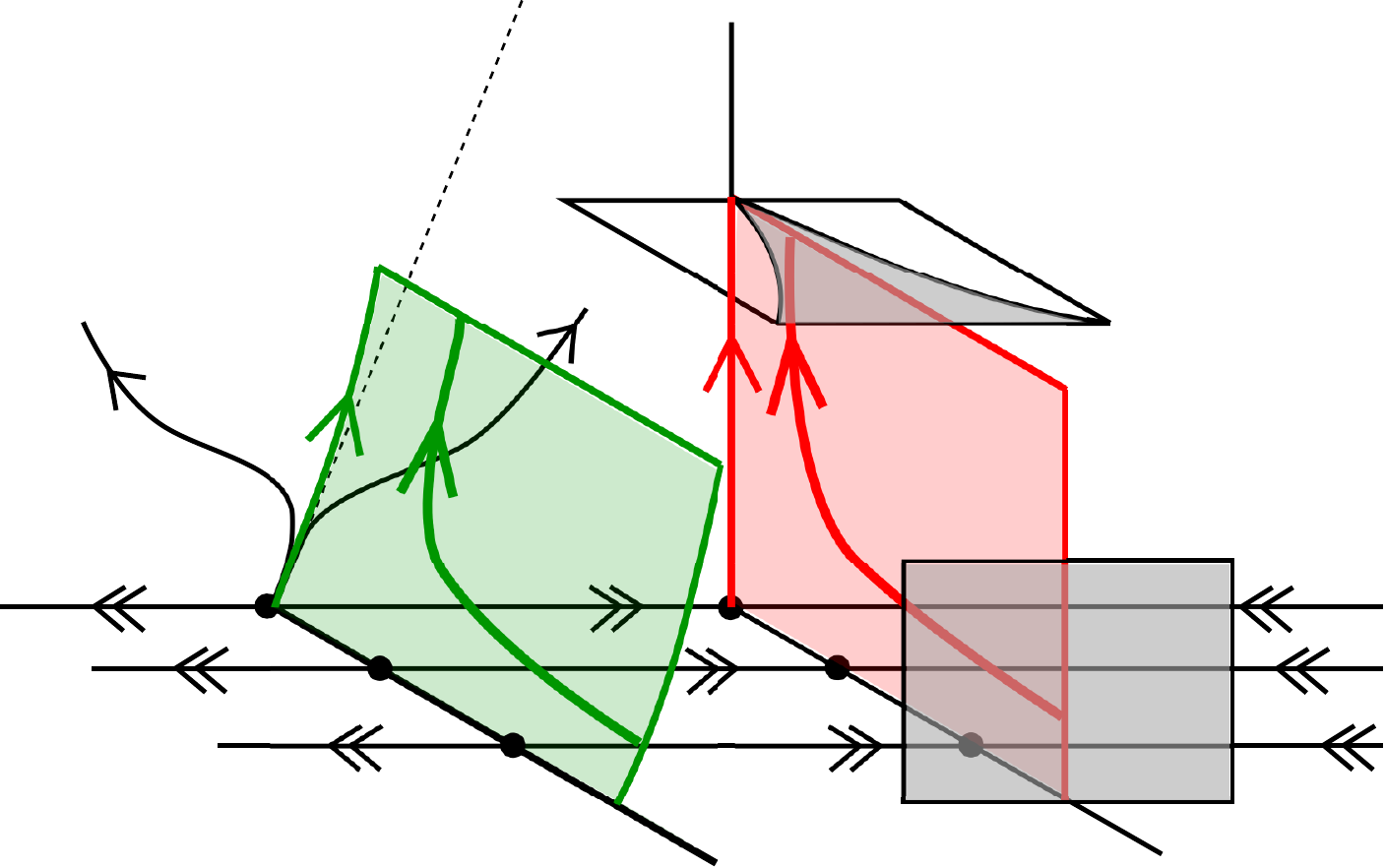}
		}
		\node at (-.75,2.15) {$\Delta_1^{\txtex}$};
		\node at (4.3,-.8) {$\Delta_1^{\txten}$};
		\node at (5.6,-1.35) {$\ba$};
		\node at (0.4,3.35) {$\be$};
		\node at (3.7,-3.45) {$\br$};
		\node[red] at (2.5,0) {$\cC_1^\txta$};
		\node[red] at (0,.8) {$\cE_1^\txta$};
		\node[green!50!black] at (-.1,-.5) {$\cC_1^\txtr$};
		\node[green!50!black] at (-2.8,1) {$\cE_1^\txtr$};
	\end{tikzpicture}
	\caption{Schematic representation of the flow of \eqref{eq:K1_0} restricted to the attracting centre manifold $\cW_1^\txtc$. The wegde-like shape of the image of $\Pi_1(\Delta_1^{\txten})$ (shaded in $\Delta_1^{\txtex}$) is due to the contraction towards $\cC_1^\txta$.}
	\label{fig:K1}
\end{figure}

The analysis in this chart is sketched in Figure \ref{fig:K1}.

\biblio

%% file: subfiles/chartK2.tex
\renewcommand{\br}{r_2}
\renewcommand{\by}{y_2}
\renewcommand{\ba}{a_2}
\renewcommand{\bb}{b_2}
\renewcommand{\bc}{c_2}
\renewcommand{\be}{\ve_2}
\newcommand{\bw}{\bar w}
\newcommand{\bA}{A_2}
\newcommand{\bB}{B_2}
\newcommand{\bC}{C_2}
\newcommand{\bW}{\bar W}

In this chart we study the dynamics of \eqref{eq:orig2} within a small neighbourhood of the origin. The corresponding blow-up map reads as
\begin{equation}
	a=\br\ba, \; b=\br\bb, \; c=\br\bc,\; y = \br\by,\; \ve=\br^2.
\end{equation}
The blown up vector field reads as
\begin{equation}\label{eq:K2-0}
	\begin{split}
		\begin{bmatrix}
			\ba'\\ \bb' \\ \bc'
		\end{bmatrix}&=-( \bar L_0  + \br\bw \bar L_1) \begin{bmatrix}
			\ba\\ \bb \\ \bc
		\end{bmatrix}\\
					\by' &= \br (-1+\cO(\br)),\\
					\br' &=0
	\end{split}
\end{equation}
where
\begin{equation}\label{eq:K2-Ls}
	\begin{split}
		\bar L_0=\begin{bmatrix}
		 \half & \half & -1 \\
		 \half & \half & -1\\
		-1 & -1 & 2
	\end{bmatrix},
	\end{split}\qquad\qquad
	\bar L_1 = \begin{bmatrix}
		1 & -1 & 0\\
		-1 & 1 & 0\\
		0 & 0 & 0
	\end{bmatrix}
\end{equation}
with $\bw = \by+\alpha_1\ba+\alpha_2\bb$. In the rest of this section we omit the equation $\br'=0$ and just keep in mind that $\br$ is a parameter in this chart.

\begin{figure}[htbp]\centering
	\begin{tikzpicture}[baseline=(current bounding box.center)]
			\node[shape=circle,draw=black, inner sep=1pt] (n1) at (0,0) {\small $1$};
			\node[shape=circle,draw=black, inner sep=1pt] (n2) at (2,0) {\small $2$};
			\node[shape=circle,draw=black, inner sep=1pt] (n3) at (1,2) {\small $3$};
			\draw[thick] (n1) -- (n2) node[midway, below] { $-\frac{1}{2}$};
			\draw[thick] (n2) -- (n3) node[midway, right] { $1$};
			\draw[thick] (n1) -- (n3) node[midway, left] { $1$};
	\end{tikzpicture}
	$\qquad+\qquad$
	\begin{tikzpicture}[baseline=(current bounding box.center)]
			\node[shape=circle,draw=black, inner sep=1pt] (n1) at (0,0) {\small $1$};
			\node[shape=circle,draw=black, inner sep=1pt] (n2) at (2,0) {\small $2$};
			\draw[thick] (n1) -- (n2) node[midway, below] { $\br\bw$};
	\end{tikzpicture}
		\caption{Graph representation of \eqref{eq:K2-0}.}
		\label{fig:K2-0}
\end{figure}
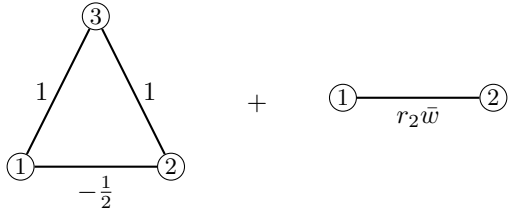

Before proceeding with the analysis, it is again very helpful to study the effect that the blow-up map has on the network's topology. Note that \eqref{eq:K2-0} can be regarded as the model of an $\cO(\br)$ graph preserving perturbation of a static network as shown in Figure \ref{fig:K2-0}.

Roughly speaking, the blow-up separates two types of dynamics: the dynamics of order $\cO(1)$ correspond to a consensus protocol on a degenerate static network. Here by degenerate we mean that the Laplacian of the static network has a kernel of dimension $2$, as can be easily seen in \eqref{eq:K2-0}-\eqref{eq:K2-Ls} with $\br=0$. Next, the dynamics of order $\cO(\br)$ occur in a slower time scale and correspond to the slowly varying edge with weight $\br\bw$.

We proceed with the description of the flow of \eqref{eq:K2-0}.

\begin{proposition}
	For $\br\geq0$ sufficiently small, the equilibrium points of \eqref{eq:K2-0} are given by
	\begin{equation}
		\left\{ (\ba,\bb,\bc,\by,\br)\in\bbR^5\,|\, \br=0,\, (\ba,\bb,\bc)\in\ker(\bar L_0) \right\}.
	\end{equation}
\end{proposition}
\begin{proof}
	Straightforward computations.
\end{proof}

Next we show that \eqref{eq:K2-0} has an attracting $4$-dimensional centre manifold $\cW_2^\txtc$ and a $1$-dimensional stable manifold $\cW_2^s$. These objects, in fact, correspond respectively to $\cW_1^\txtc$ and $\cW_1^\txts$ found in chart $K_1$. In qualitative terms, reduction to $\cW_2^\txtc$ will correspond to representing the behaviour of the third node, with state $\bc$, in terms of the other two nodes.

\begin{proposition}\label{prop:K2-cm} System \eqref{eq:K2-0} has a $4$-dimensional local centre manifold $\cW_2^\txtc$ and a $1$-dimensional local stable manifold $\cW_2^s$ that intersect at $\left\{\br=0\right\}\cap\left\{\bc=\frac{\ba+\bb}{2}\right\}$.
The centre manifold $\cW_2^\txtc$ is given by the graph of $\bc = \frac{\ba+\bb}{2}$ and it hold that $\kappa_{12}(\cW_1^\txtc)=\cW_2^\txtc$.
\end{proposition}
\begin{proof}
	The proof follows the same reasoning (and in fact it is simpler than) the proof of Proposition \ref{prop1K1}. The relation $\kappa_{12}(\cW_1^\txtc)=\cW_2^\txtc$ is straightforward from \eqref{eq:matching}.
\end{proof}

Since the centre manifold $\cW_2^c$ is attracting, and of codimension $1$, the next step is to restrict the dynamics to it. However, the next observation is important (recall Proposition \ref{prop:K1sign}).

\begin{lemma}\label{lemma:signs} The trajectories of \eqref{eq:K2-0} restricted to $\left\{\br=0 \right\}$ have the asymptotic behaviour
\begin{equation}
	\lim_{t_2\to\infty}\left( \ba(t_2),\bb(t_2),\bc(t_2) \right)=\frac{1}{2}\left(  \ba(0)-\bb(0),\bb(0)-\ba(0),0 \right).
\end{equation}
\end{lemma}

As it was the case in chart $K_1$ the previous lemma gives us the relevant sign of $\ba$ on the centre manifold $\cW_2^\txtc$.

The restriction of \eqref{eq:K2-0} to $\cW_2^\txtc$ results on a vector field of order $\cO(\br)$, which can be desingularized as is usual in the blow-up method by dividing by $\br$. By performing the aforementioned steps we obtain
\begin{equation}
	\begin{split}
		\ba' &= \bw(\bb-\ba),\\
		\bb' &= \bw(\ba-\bb),\\
		\by' &= -1+\cO(\br),
	\end{split}
\end{equation}
where we recall that $\bw=\by+\alpha_1\ba+\alpha_2\bb$. From the fact that $a+b+c=\br(\ba+\bb+\bc)=0$ for all $\br\geq0$ and due to the restriction to $\cW_2^\txtc$, that is $\bc=\frac{\ba+\bb}{2}$, we further have that $\ba+\bb=0$. Therefore, the analysis of \eqref{eq:K2-0} is reduced to the analysis of the planar system
\begin{equation}\label{eq:K2-red}
	\begin{split}
		\ba' &= -2(\by+\underbrace{(\alpha_1-\alpha_2)}_{=:\nu\geq0}\ba)\ba,\\
		\by' &= -1+\cO(\br).
	\end{split}
\end{equation}

Note that in the restriction of \eqref{eq:K2-red} to $\left\{ \br=0\right\}$, one has that $\by$ is essentially time in the reverse direction. To describe the flow of \eqref{eq:K2-red}, let $\delta_2>0$ and define the sections
\begin{equation}\label{eq:K2sections1}
	\begin{split}
		\Delta_{2}^\txten &= \left\{ (\br,\ba,\by)\in\bbR^3\, | \, \by=\delta_2  \right\},\\
		\Delta_{2}^\txtex &= \left\{ (\br,\ba,\by)\in\bbR^3\, | \, \by=-\delta_2 \right\}.
	\end{split}
\end{equation}
Accordingly, let $\Pi_{2}:\Delta_{2}^\txten\to\Delta_{2}^\txtex$ be the map defined by the flow of\eqref{eq:K2-red}. We now show the following.

\begin{proposition}\label{prop:K21} Consider  \eqref{eq:K2-red}. Then the following hold.
	\begin{enumerate}[leftmargin=*]
		\item There exists a trajectory $\gamma_c$ given by
		\begin{equation}
			\gamma_c(t_2)=\left( \ba(t_2),\by(t_2) \right)=(0,-t_2)\, \qquad t_2\in\bbR.
		\end{equation}
		No other trajectory of \eqref{eq:K2-red} converges to the $\by$-axis as $t_2\to\pm\infty$.
		\item There exist orbits $\gamma_2^\txtr$ and $\gamma_2^\txta$ that are defined, respectively, in the quadrants $\left\{ \ba<0,\by>0\right\}$ and $\left\{ \ba>0,\by<0\right\}$ and are given by
		\begin{equation}\label{eq:K2gammas}
		\begin{split}
			\gamma_2^\txtr &= \left\{ (\ba,\by)\in\bbR^2\,|\, \ba=-\frac{1}{2\nu D_+(\by)}, \, \ba<0 \right\},\\
			\gamma_2^\txta &= \left\{ (\ba,\by)\in\bbR^2\,|\, \ba=-\frac{1}{2\nu D_+(\by)}, \, \ba>0 \right\},
		\end{split}
		\end{equation}
		where $D_+(\by)$ stands for the Dawson function \cite[pp. 219 and 235]{abramowitz1972handbook}. Furthermore, since $\by$ is essentially time, the trajectory $\gamma_2^j$, $j=\txtr,\txta$, has asymptotic expansions
		\begin{equation}
			\begin{aligned}
				\gamma_2^j &= -\frac{1}{2\nu\by} + \cO(\by^{-3}), &&\qquad \by\to 0\\[1ex]
				\gamma_2^j &= -\frac{1}{\nu}\by + \frac{1}{2\nu\by} + \cO(\by^{-3}), && \qquad |\by|\to\infty.
			\end{aligned}
		\end{equation}
		All trajectories of \eqref{eq:K2-red} with initial condition $\ba^*>0$ and $\by^*>0$ are asymptotic to $\gamma_2^\txta$ as $t_2\to\infty$.
		\item The transition map $\Pi_{2}:\Delta_{2}^\txten\to\Delta_{2}^\txtex$ is well-defined if and only if 
		\begin{equation}\label{eq:K2cond}
			\ba|_{\Delta_{2}^\txten}>-\frac{1}{4D_+(\delta_2)},
		\end{equation}
		and given by
		\begin{equation}
			\Pi_{2}\left( \begin{matrix}
				\br\\ \ba \\ \delta_2
			\end{matrix}\right) = \left( \begin{matrix}
				\br\\\frac{\ba}{1+4\ba\nu D_+(\delta_2)} +\cO(\br)\\ -\delta_2
			\end{matrix}\right).
		\end{equation}
		\item If $\ba|_{\Delta_{2}^\txten}\leq-\frac{1}{4D_+(\delta_2)}$ then the corresponding orbit has asymptote $\by=\Gamma_2$ which is implicitly given by
		\begin{equation}\label{eq:K2-Gamma2}
			\exp(\Gamma_2^2)D_+(\Gamma_2) = \exp(\delta_2^2)\left( \frac{1}{2\ba^*\nu}+D(\delta_2) \right). 
		\end{equation}
		\end{enumerate}

\end{proposition}
\begin{proof}
The first item follows from the invariance of $\ba=0$ and linear analysis along the $\by$-axis. For the second item, no distinction between the orbits is needed, one only needs to check that $\gamma_2^j$ satisfies \eqref{eq:K2-red}, for which~\cite{abramowitz1972handbook} $D_+'(\by)=1-2\by D_+(\by)$ is useful. More precisely, from the expression $\gamma_2^j=\left\{\ba=-\frac{1}{2\nu D_+(\by)} \right\}$, one has 
\begin{equation}
	\begin{split}
		\ba' &= \frac{D_+'(\by)\by'}{2\nu D_+^2(\by)}=\frac{2\by D_+(\by)-1}{2\nu D_+^2(\by)} = \frac{1}{\nu D_+(\by)}\left( \by - \frac{1}{2D_+(\by)}\right)=-2(\by+\nu\ba)\ba.
	\end{split}
\end{equation}

Next, the asymptotic expansions for $\gamma_2^j$ follow directly from~\cite{abramowitz1972handbook}, where one finds
\begin{equation}
	\begin{aligned}
		D_+(\xi) &= \xi+\cO(\xi^3), \qquad &&\xi\to0\\[1ex]
		D_+(\xi) &= \frac{1}{2\xi}+\cO(\xi^{-3}), \qquad &&\xi\to\infty.
	\end{aligned}
\end{equation}

The fact that $\gamma_2^\txta$ attracts all trajectories with the given initial conditions follows from: i) $\ba=0$ is invariant, ii) in the limit $|\by|\to\infty$ the curve $\gamma_2^\txta$ is asymptotic to $\by+\nu\ba=0$, and iii) the set $\left\{ \by+\nu\ba=0\right\}$ is attracting in the quadrant $\ba>0$, $\by<0$.

For the transition map we have that \eqref{eq:K2-red} has an explicit solution given by
\begin{equation}\label{eq:K2sol}
	\ba(\by) = \frac{\ba^*}{\exp((\by^*)^2-\by^2)(1+2\ba^*\nu D_+(\by^*))-2\ba^*\nu D_+(\by)},
\end{equation}
where $(\ba^*,\by^*)$ denotes an initial condition. Thus, for the map $\Pi_{2}$ to be well-defined we need to ensure that the denominator in \eqref{eq:K2sol} does not vanish. Let us substitute $(\ba^*,\by^*)=(\ba^*,\delta_2)$ with $\delta_2>0$, and compute $\ba(-\delta_2)$. For this it is useful to recall that $D_+$ is an odd function. So we get
\begin{equation}
	\ba(-\delta_2) = \frac{\ba^*}{1+4\ba^*\nu D_+(\delta_2)},
\end{equation}
which indeed leads to \eqref{eq:K2cond} and the form of $\Pi_{2}$ also follows. Finally, the expression of the asymptote $\Gamma_2$ is obtained by solving the denominator of \eqref{eq:K2sol} equal to $0$ and with initial condition $(\ba^*,\by^*)=(\ba^*,\delta_2)$.
\end{proof}

\begin{remark}\leavevmode
\begin{itemize}[leftmargin=*]
	\item In particular, it follows from the third item of Proposition \ref{prop:K21} that the map $\Pi_2(\br,\ba,\delta_2)$ is well defined for all $\ba\geq0$.
	\item For $\delta_2>0$ sufficiently large and $\ba^*$ sufficiently small one has that $\ba(-\delta_2)\approx \ba^*$.
	\item If $\alpha_1-\alpha_2=0$ then $\Pi(\ba,\delta_2)=(\ba,-\delta_2)$.

\end{itemize}
	
\end{remark}

We now relate the curves $\gamma_2$ and $\gamma_c$ with centre manifolds found in chart $K_1$.
\begin{proposition}\label{prop:K1ident}
	The curves $\gamma_2^\txtr$ and $\gamma_c$ correspond, respectively, to the centre manifolds $\cE_1^\txtr$ and $\cE_1^\txta$ of chart $K_1$.
\end{proposition}
\begin{proof}
	We detail the relation between $\gamma_2$ and $\cE_1^\txtr$, the correspondence between $\gamma_c$ and $\cE_1^\txta$ is trivial since they are given by $\left\{\ba=0\right\}$ and $\left\{a_1=0\right\}$ respectively. We can transform $\gamma_2$ into the coordinates of chart $K_1$ via the map $\kappa_{21}$, which gives
	\begin{equation}\label{eq:K2g1}
		\ve_1^{-1/2}a_1=-\frac{1}{2\nu D_+(\ve_1^{-1/2})}.
	\end{equation}
	Taking the limit $\ve_1\to0$ in \eqref{eq:K2g1} one gets $a_1=-\frac{1}{\nu}=\frac{1}{\alpha_2-\alpha_1}$. Thus the claim follows from the analysis performed in chart $K_1$ particularly for $r_1=0$.
\end{proof}

\begin{remark}
	The trajectory $\gamma_c$ corresponds to a singular maximal canard of \eqref{eq:orig2}, while $\gamma_2^\txtr$ and $\gamma_2^\txta$ correspond to the manifolds $\cM_0^\txtr$ and $\cM_0^\txta$. Accordingly,  $\cO(\br)$-small perturbation of such orbits correspond to $\cN_\ve$, $\cM_\ve^\txtr$ and $\cM_\ve^\txta$ for $\ve>0$ sufficiently small.
\end{remark}

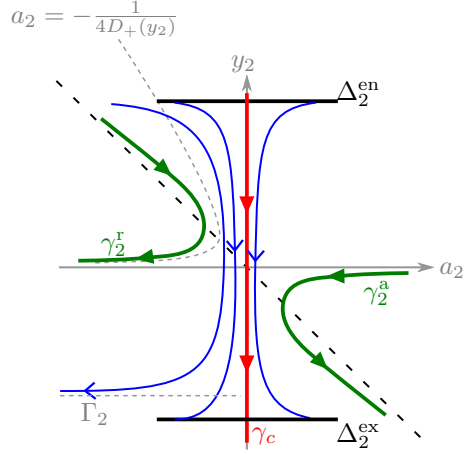
\begin{figure}[htbp]\centering
	\begin{tikzpicture}
		\node at (0,0){
			\includegraphics{\main/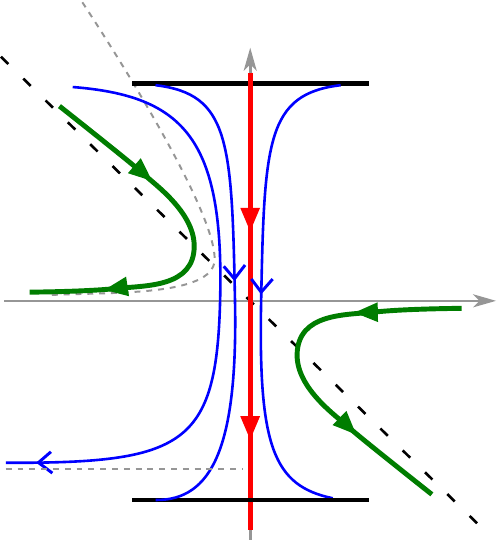}
		};
		\node[color=gray] at (2.75,-.35) {$\ba$};
		\node[color=gray] at (-.025,2.4) {$\by$};
		\node at (1.5,2.0){$\Delta_2^\txten$};
		\node at (1.5,-2.55){$\Delta_2^\txtex$};
		\node[red] at (.25,-2.55){$\gamma_c$};
		\node[green!50!black] at (1.75,-.65){$\gamma_2^\txta$};
		\node[green!50!black] at (-1.75,0){$\gamma_2^\txtr$};
		\node[gray] at (-2,-2.25){$\Gamma_2$};
		\node[gray] at (-2,3){$\ba=-\frac{1}{4D_+(\by)}$};
	\end{tikzpicture}
	\caption{Flow of \eqref{eq:K2-0} along $\cW_2^\txtc$. This flow is equivalent to that of \eqref{eq:orig2} within a small neighborhood of the origin and for $\ve>0$ sufficiently small. We observe that trajectories starting at $\Delta_2^\txten$ are first attracted to the invariant set $\left\{ \ba=0 \right\}$, which in terms of the network means consensus. Then, once the trajectories pass through the origin, they are repelled from consensus. All trajectories with initial condition $\ba^*>0$ are eventually attracted towards $\gamma_2^\txta$, which in terms of the original coordinates corresponds to the clustering manifold. }
	\label{fig:K2}
\end{figure}

The analysis performed in this chart is sketched in Figure \ref{fig:K2}.

\biblio

%% file: subfiles/chartK3.tex
\renewcommand{\br}{r_3}
\renewcommand{\by}{y_3}
\renewcommand{\ba}{a_3}
\renewcommand{\bb}{b_3}
\renewcommand{\bc}{c_3}
\renewcommand{\be}{\ve_3}
\newcommand{\bw}{\bar w}
\newcommand{\bA}{A_3}
\newcommand{\bB}{B_3}
\newcommand{\bC}{C_3}
\newcommand{\bW}{\bar W}

The analysis in this chart is similar to that in chart $K_1$ performed in Section \ref{sec:K1}. Therefore, we shall only point out the main information required from this chart, and omit the proofs.

In this chart the blow up map is given by
\begin{equation}
	a=\br\ba, \; b=\br\bb, \; c=\br\bc,\; y = -\br,\; \ve=\br^2\be.
\end{equation}

We then obtain the blown up vector field
\begin{equation}\label{eq:K3_0}
	\begin{split}
		\begin{bmatrix}
			\ba' \\ \bb' \\ \bc'
		\end{bmatrix} &=
		-\begin{bmatrix}
			\half & \half & -1 \\
			\half & \half & -1\\
			-1 & -1 & 2
		\end{bmatrix}\begin{bmatrix}
			\ba\\ \bb \\ \bc
		\end{bmatrix} + \br f_3(\ba,\bb,\bc,\br,\be)\\
		\br' &= -\br^2\be(-1+\cO(\br))\\
		\be' &= 2\br\be^2(-1+\cO(\br)),
	\end{split}
\end{equation}
where $f_3(\ba,\bb,\bc,\br,\be)$ reads as
\begin{equation}\label{eq:f3}
	f_3 = \left( (-1+\alpha_1\ba+\alpha_2\bb) \begin{bmatrix}
		-1 & 1 & 0 \\
		1 & -1 & 0\\
		0 & 0 & 0
	\end{bmatrix}  + \be(-1+\cO(\br)) \begin{bmatrix}
		1 & 0 & 0\\
		0 & 1 & 0\\
		0 & 0 & 1
	\end{bmatrix} \right)\begin{bmatrix}
		\ba\\ \bb \\ \bc
	\end{bmatrix}.
\end{equation}

The flow of  \eqref{eq:K3_0} is described as follows.

	\begin{proposition} The following claims hold for  \eqref{eq:K3_0}.

	\begin{enumerate}[leftmargin=*]
		\item  There exist a $1$-dimensional local stable manifold $\cW_3^\txts$ and a $4$-dimensional local centre-stable manifold  $\cW_3^\txtc$. The centre manifold is given by the graph of $\bc=\frac{\ba+\bb}{2}$.
	\end{enumerate}
	Restricted to $\cW_3^\txtc$ one has $\bb=-\ba$, $\bc=0$, and:
	\begin{enumerate}[leftmargin=*]
		\setcounter{enumi}{1}
		\item There is a repelling $2$-dimensional centre manifold $\cC_3^\txtr$. The manifold $\cC_3^\txtr$ contains a line of zeros $\ell_3^\txtr=\left\{(\br,\ba,\be)\in\bbR^3\, | \, \ba=\be=0 \right\}$ and a $1$-dimensional centre manifold\\
		$\cE_3^\txtr=\left\{(\br,\ba,\be)\in\bbR^3\, | \, \br=\ba=0\right\}$. On the plane $\left\{ \br=0 \right\}$, the centre manifold $\cE_3^\txtr$ is unique. The flow along $\cE_3^\txtr$ is stable, that is, it converges to the origin; while the flow on $\cC_3^\txtr$ away from $\ell_3^\txtr$ is locally equivalent to that of a saddle.
		\item There is an attracting $2$-dimensional centre manifold $\cC_3^\txta$. The manifold $\cC_3^\txta$ contains a line of zeros $\ell_3^\txta=\left\{(\br,\ba,\be)\in\bbR^3\, | \, \ba=\frac{1}{\alpha_1-\alpha_2},\,\be=0 \right\}$ and a $1$-dimensional centre manifold $\cE_3^\txta=\left\{(\br,\ba,\be)\in\bbR^3\, | \, \br=0, \ba=\frac{1}{\alpha_1-\alpha_2}+\cO(\be)\right\}$. The flow along $\cE_3^\txta$ is stable, that is, it converges to the equilibrium point $(\ba,\be)=\left(\frac{1}{\alpha_1-\alpha_2},0\right)$; while the flow on $\cC_3^\txta$ away from $\ell_3^\txta$ is locally equivalent to that of a saddle.
	
	\end{enumerate}

	Define the sections
\begin{equation}\label{eq:K3sections}
	\begin{split}
		\Delta_3^\txten &= \left\{ (\ba,\br,\be)\in\bbR^3\, | \, \br<\delta_3,\, \be=\mu_3 \right\}\\
		\Delta_3^\txtex &= \left\{ (\ba,\br,\be)\in\bbR^3\, | \, \br=\delta_3,\, \be<\mu_3 \right\},
	\end{split}
\end{equation}
where $\delta_3>0$,  and $\mu_3>0$ is sufficiently small. Let $\Pi_3:\Delta_3^\txten\to\Delta_3^\txtex$ denote the map induced by the flow of \eqref{eq:K3_0} restricted to $\cW_3^c$. Then $\Pi_3$ has the form
	\begin{equation}\label{eq:K3P3}
		\Pi_3\left( \begin{matrix}
	  \ba\\ \br \\ \mu_3
		\end{matrix}\right)=\left( \begin{matrix}
			h_{\ba} +\cO(\ba\mu_3^2)\\
			\delta_3\\
			\mu_3\left( \frac{r_3}{\delta_3}\right)^2
		\end{matrix} \right),
	\end{equation}
	where the function $h_{\ba}=h_{\ba}(\ba,\delta_3,\be)$ reads as
	\begin{equation}
		h_{\ba} = \frac{\ba\exp(2T_3)}{(\alpha_1-\alpha_2)\ba(\exp(2T_3)-1)+1},\qquad T_3=\frac{1}{2\mu_3}\left( \left( \frac{\delta_3}{r_3}\right)^2-1 \right).
	\end{equation}
		\end{proposition}

\begin{remark}\leavevmode
\begin{itemize}[leftmargin=*]
	\item If $\ba\geq0$, then the function $h_{\ba}$ is well-defined for any point $(\ba,\br,\mu_3)\in\Delta_3^\txten$. If $\ba<0$, the function $h_{\ba}$ is well-defined only for $T_3<\frac{1}{2}\ln\left( 1+\frac{1}{(\alpha_1-\alpha_2)|\ba|} \right)$. In such a case, we choose suitably $0<\br<\delta_3$ so that the function $h_{\ba}$ is well-defined. 
	\item For $\ba>0$, and $T_3>0$ sufficiently large, one has $h_{\ba}\approx\frac{1}{\alpha_1-\alpha_2}$.
	\item If $\alpha_1-\alpha_2=0$ then $h_{\ba}\approx \ba\exp(2T_3)$.
	\item $\kappa_{23}(\gamma_2^\txtc)=\cE_3^\txtr$ and $\kappa_{23}(\gamma_2^\txta)=\cE_3^\txta$.
\end{itemize}
\end{remark}

\begin{figure}[htbp]\centering
	\begin{tikzpicture}
		\pgftext{
		\includegraphics[scale=.75]{\main/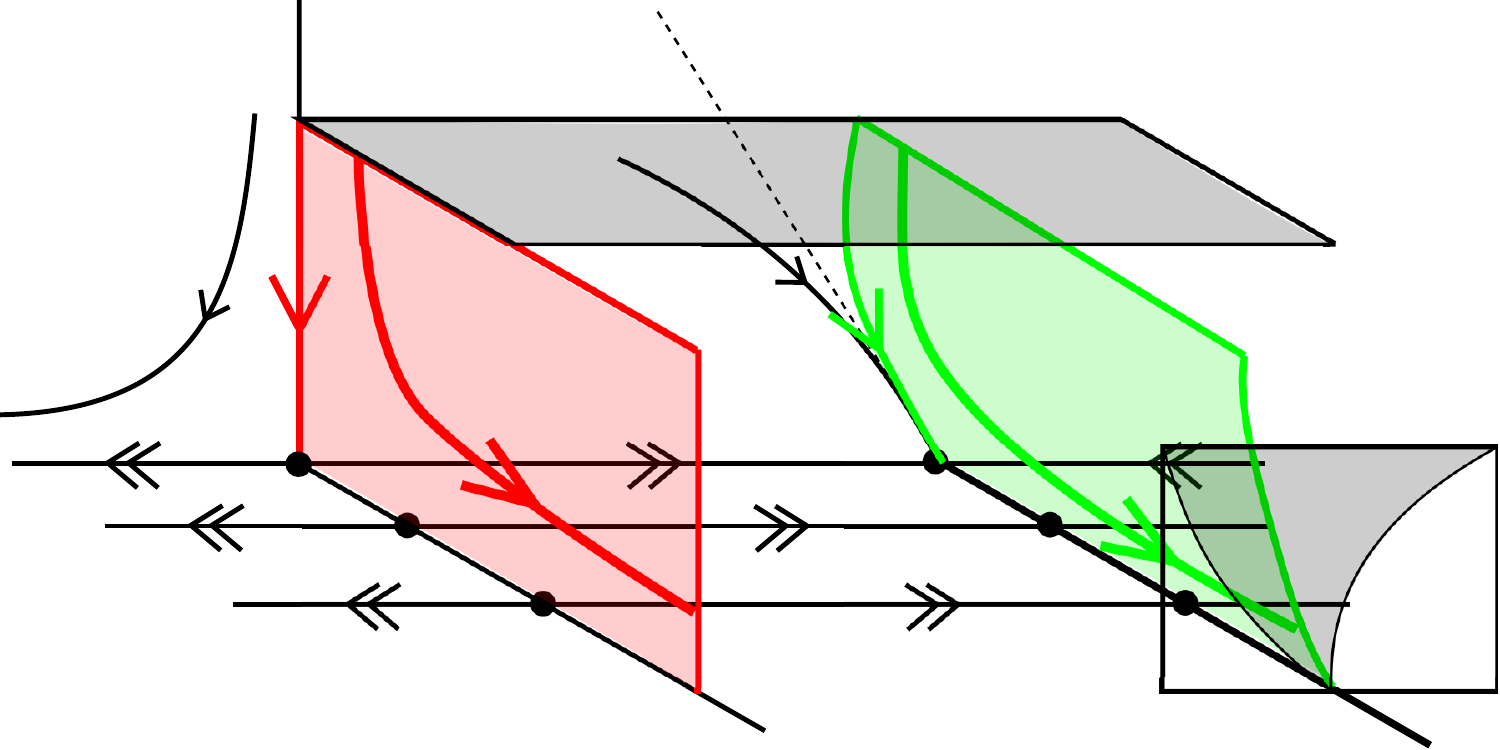}
		}
		\node at (2,2.2) {$\Delta_3^{\txten}$};
		\node at (3.8,-2.8) {$\Delta_3^{\txtex}$};  
		\node at (-6.15,-.7) {$-\ba$};
		\node at (-3.4,3) {$\be$};
		\node at (0.4,-2.9) {$\br$};
		\node[green!80!black] at (2.5,.5) {$\cC_3^\txta$};
		\node[red] at (-1.5,.5) {$\cC_3^\txtr$};

		\node[green!80!black] at (.8,-.2) {$\cE_3^\txta$};
		\node[red] at (-3.75,-.2) {$\cE_3^\txtr$};
	\end{tikzpicture}
	\caption{Schematic representation of the flow of \eqref{eq:K3_0} restricted to the attracting centre manifold $\cW_3^\txtc$. The wegde-like shape of the image of $\Pi_3(\Delta_3^{\txten})$ (shaded in $\Delta_1^{\txtex}$) is due to the contraction towards $\cC_3^\txta$.}
	\label{fig:K3}
\end{figure}

The flow in this chart is as depicted in Figure \ref{fig:K3}.

\biblio

%% file: subfiles/transition.tex
In this section we prove items (T4)-(T6) of Theorem \ref{thm:main1}. First of all note that if we choose $\delta_1=\delta_3=\delta$, then the sections $\Delta_1^\txten$ and $\Delta_3^\txtex$ are precisely the sections $\Sigma^\txten|_\cA$ and  $\Sigma^\txtex|_\cA$ in the blow-up coordinates. Moreover, the set $\cA$ corresponds, in each chart, to the centre manifold $\cW_1^\txta$, $\cW_2^\txta$, and $\cW_3^\txta$ respectively. Thus it will suffice to consider the transition map $\bar\Pi:\Delta_1^\txten\to\Delta_3^\txtex$ in the blow-up space (or equivalently $\Pi|_\cA$).

 The map $\bar\Pi(\Delta_1^\txten)$ is then given as
\begin{equation}
	\bar\Pi(\Delta_1^\txten) = \Pi_3\circ\kappa_{23}\circ\Pi_2\circ\kappa_{12}\circ\Pi_1(\Delta_1^\txten),
\end{equation}
where the maps $\Pi_1$, $\Pi_2$, and $\Pi_3$ are given in Sections \ref{sec:K1}, \ref{sec:K2}, and \ref{sec:K3} respectively, and where the maps $\kappa_{12}$ and $\kappa_{23}$ are defined in Lemma \ref{lemma:matching}. We compute $\bar\Pi(\Delta_1^\txten)$ as follows. For brevity we disregard the higher order terms in the chart maps.
\begin{enumerate}[leftmargin=*]
	\item We start from $\Delta_1^\txten=(a_1,\delta,\ve)$ and compute $\Pi_1(\Delta_1^\txten)=(h_{a_1},\delta \ve_1^{1/2}\mu^{-1/2},\mu)$, where $h_{a_1}$ is as in \eqref{eq:K1ha} and we let $\mu_1=\mu$. 
	\item Next we compute $\kappa_{12}\circ\Pi_1(\Delta_1^\txten)$ from \eqref{eq:matching}, obtaining
	\begin{equation}
		\kappa_{12}\circ\Pi_1(\Delta_1^\txten)=\left(\delta\ve^{1/2},\mu^{-1/2}h_{a_1},\mu^{-1/2} \right).
	\end{equation}
	By defining $\mu^{1/2}=\delta_2$ we have from \eqref{eq:K2sections1} that $\kappa_{12}\circ\Pi_1(\Delta_1^\txten)\subset\Delta_2^\txten$. 
	\item Next we can compute $\Pi_2\circ\kappa_{12}\circ\Pi_1(\Delta_1^\txten)$ using Proposition \ref{prop:K21}. We get
	\begin{equation}
		\Pi_2\circ\kappa_{12}\circ\Pi_1(\Delta_1^\txten) = \left(  \delta \ve_1^{1/2}, \underbrace{\frac{\mu^{-1/2}h_{a_1}}{1+4\mu^{-1/2}h_{a_1}(\alpha_1-\alpha_2)D_+(\mu^{-1/2})} }_{h_{a_2}},-\mu^{-1/2} \right).
	\end{equation}
	\item Next we compute $\kappa_{23}\circ\Pi_2\circ\kappa_{12}\circ\Pi_1(\Delta_1^\txten)$ again using \eqref{eq:matching}, obtaining
	\begin{equation}
		\kappa_{23}\circ\Pi_2\circ\kappa_{12}\circ\Pi_1(\Delta_1^\txten)=\left( \mu^{1/2}h_{a_2},\delta\ve_1^{1/2}\mu
		^{-1/2},\mu \right).
	\end{equation}
	By defining $\mu_3=\mu$ we have from \eqref{eq:K3sections} that $\kappa_{23}\circ\Pi_2\circ\kappa_{12}\circ\Pi_1(\Delta_1^\txten)\subset\Delta_3^\txten$.
	\item Finally we compute $\Pi_3\circ\kappa_{23}\circ\Pi_2\circ\kappa_{12}\circ\Pi_1(\Delta_1^\txten)$ from \eqref{eq:K3P3}, obtaining
	\begin{equation}
		\Pi_3\circ\kappa_{23}\circ\Pi_2\circ\kappa_{12}\circ\Pi_1(\Delta_1^\txten) = \left( h_{a_3},\delta, \ve_1\right),
	\end{equation}
	where $h_{a_3}= \frac{ \mu^{1/2}h_{a_2}\exp(2T_3)}{(\alpha_1-\alpha_2) \mu^{1/2}h_{a_2}(\exp(2T_3)-1)+1}$ reads, after substitutions, as
	\begin{equation}
		h_{a_3}=\frac{a_1\exp(2 T)}{ (\alpha_1-\alpha_2)a_1\left( 2\exp(2 T)-2+4\mu^{-1/2}D_+(\mu^{-1/2}) \right)+\exp(2T)}
	\end{equation}
	where $T=T_1=\frac{1}{2}\left( \frac{1}{\ve_1}-\frac{1}{\mu} \right)(1+\cO(\delta))$. We see that, as it was already evident in each chart, the function $h_{a_3}$ is well-defined for $a_1\geq0$, while for $a_1<0$ the function $h_{a_3}$ is well-defined only for finite time $T$. We thus assume that, in either case, we choose appropriate constants $\delta>0$, and $\mu>0$ sufficiently small, such that $h_{a_3}$ is well-defined.

\end{enumerate}

Note that, restricted to $\cW_j^\txtc$, all the sets $\left\{a_j=0\right\}$, $j=1,2,3$, are invariant along the blow-up space. In other words, if $(\bar{a},\bar{b},\bar{c},\bar{y},\bar\ve,r)$ denote the global blow-up coordinates as in \eqref{eq:bu}, then we have that $\left\{ \ba=0\right\}$ is invariant on the set $\left\{ \bar c=\frac{\bar a+\bar b}{2}\right\}$. We now proceed with proving each item (T4)-(T6) of Theorem \ref{thm:main1}.

\begin{itemize}[leftmargin=*]
	\item[(T4)] Indeed, we have that $\bar\Pi(0,\delta,\ve_1)=(0,\delta,\ve_1)$. Moreover, it follows from item 1 of Proposition \ref{prop:K21} and Section \ref{sec:K3} that the curve $\gamma_c$ connects the point $(a_1,r_1,\ve_1)=(0,0,0)$ in chart $K_1$ with the point $(a_3,r_3,\ve_3)=(0,0,0)$ in chart $K_3$. Therefore, the center manifolds $\cC_1^\txta$ and $\cC_3^\txtr$ are also connected in the central chart $K_2$ via the map $\Pi_2$ for $r_2\geq0$ sufficiently small. It is then clear that, since we can identify $\cC_1^\txta$ with $\cN_\ve^\txta$ and $\cC_1^\txtr$ with $\cN_\ve^\txtr$, the manifolds $\cN_\ve^\txta$ and $\cN_1^\txtr$ are also connected for $\ve\geq0$ sufficiently small. The stability of $\cN_\ve$ follows from $\cC_1^\txta$ being attracting in chart $K_1$ and $\cC_3^\txtr$ being repelling in chart $K_3$.
	\item[(T5)] If $\alpha_1-\alpha_2=0$ we have that $h_{a_3}(a_1,\delta_1,\ve_1)=a_1$. Since $\Delta_3^\txtex$ is sufficiently away from the origin, the claim for $a\neq0$ follows from Fenichel's theory and the stability properties of $\cN_0$.
	\item[(T6)] The expression of $\cM_\ve^\txta$ and $\cM_\ve^\txtr$ are obtained from blowing down $\gamma_2^\txta$ and $\gamma_2^\txtr$ given in \eqref{eq:K2gammas}. Accordingly the fact that $\cM_\ve^\txta$ attracts all trajectories with initial condition $a_1|_{\Delta_1^\txten}>0$ follows equivalent arguments as those for the second item if Proposition \ref{prop:K21}. On the contrary, when we have  $a_1|_{\Delta_1^\txten}<0$ we see from the expression of $h_{a_3}$ that the trajectories become unbounded in finite time $T$, see also the remark at the end of Section \ref{sec:K3}. The proof is finalized by recalling the relationship between the signs of initial conditions in $\Sigma^\txten$ and the corresponding sign of $a_j$ in $\cW_j^\txtc$ given at the beginning of the proof.
\end{itemize}

\biblio

%% file: subfiles/generalizations.tex
Here we develop a couple of generalizations for the results presented in Section \ref{sec:triangle}. The first one is concerned with triangle motifs with one dynamic weight while the other two weights are fixed and positive, but not necessarily equal. The second generalization deals with consensus protocols defined on arbitrary graphs, where just one weight is dynamic.

\subsection{The nonsymetric triangle motif}\label{sec:nonsym}
\begin{figure}[htbp]\centering
\begin{tikzpicture}
		\node[shape=circle,draw=black, inner sep=1pt] (n1) at (0,0) {\small $1$};
		\node[shape=circle,draw=black, inner sep=1pt] (n2) at (2,0) {\small $2$};
		\node[shape=circle,draw=black, inner sep=1pt] (n3) at (1,2) {\small $3$};
		\draw[thick] (n1) -- (n2) node[midway, below] { $w$};
		\draw[thick] (n2) -- (n3) node[midway, right] { $w_{23}$};
		\draw[thick] (n1) -- (n3) node[midway, left] { $w_{13}$};
		%\node at (5,1) {$w_{ij}=w_{ji}$};
	\end{tikzpicture}
	\caption{Triangle motif: a network of three nodes connected on a cycle. The weights $w_{13}$ and $w_{23}$ are fixed and positive.}
	\label{fig:triangle_ns}
\end{figure}
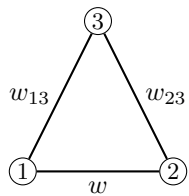

In Section \ref{sec:triangle} we studied a triangle motif with the fixed weights equal to $1$. \emph{In this section we show that the non-symmetric case is topologically equivalent to the symmetric one}.  Let us start by considering the fast-slow system
\begin{equation}\label{eq:sfs2}
	\begin{split}
		x' &= -L(x,y,\ve)x\\
		y' &= \ve g(x,y,\ve),
	\end{split}\qquad\qquad
	L=\begin{bmatrix}
		w + w_{13} & -w & -w_{13}\\
		-w & w+w_{23} & -w_{23}\\
		-w_{13} & -w_{23} & w_{13}+w_{23}
	\end{bmatrix},
\end{equation}
where $w=w^*+y+\alpha_1x_1+\alpha_2x_2$, the weights $w_{13}$ and $w_{23}$ are fixed and positive, and $w^*\in\bbR$ is such that $\dim\ker L|_{\left\{w=w^*\right\}}=2$, see more details below. We recall that, from the arguments at the beginning of Section \ref{sec:main}, we may also assume that $x_1+x_2+x_3=0$. System \eqref{eq:sfs2} corresponds to the network shown in Figure \ref{fig:triangle_ns}.

\begin{remark}
	If $w_{13}=w_{23}=\tilde w>0$, one can show, for example using the exact same transformation $T$ of Lemma \ref{lemma:T}, that the eigenvalues of $L$ are $\left\{ 0,3\tilde w, 2 w+ \tilde w \right\}$. Thus, the analysis in this case is completely equivalent to the one already performed in Section \ref{sec:triangle}. The only difference would be the rate of convergence towards the set $\left\{x_3=\frac{x_1+x_2}{2}\right\}$. Therefore, in this section we rather assume $w_{13}\neq w_{23}$.
\end{remark}

It is straightforward to show that the spectrum of $L$ is given by
\begin{equation}
	\spec L = \left\{ 0,\lambda_2,\lambda_3\right\} = \left\{ 0,  w+w_{13}+w_{23} \pm\sqrt{w^2+w_{13}^2+w_{23}^2 - ww_{13} - ww_{23}- w_{13}w_{23}} \right\}.
\end{equation}
We note the following:
\begin{itemize}
	\item If $ww_{13} + ww_{23}+ w_{13}w_{23}>0$ then $\lambda_2>0$, $\lambda_3>0$,
	\item If $ww_{13} + ww_{23}+ w_{13}w_{23}=0$ then $\lambda_2>0$, $\lambda_3=0$,
	\item If $ww_{13} + ww_{23}+ w_{13}w_{23}<0$ then $\lambda_2>0$, $\lambda_3<0$.
\end{itemize}

This means that, as in the symmetric case, the Laplacian matrix has always a positive eigenvalue $\lambda_2(w)$, and another $\lambda_3(w)$ whose sign depends on $w$.

\begin{proposition}\label{prop:eqs_ns} Let $\ve=0$. Then the equilibrium points of \eqref{eq:sfs2} are given by 
\begin{itemize}[leftmargin=*]
	\item $x_1=x_2=x_3=0$, that is consensus.
	\item $w=-\frac{w_{13}w_{23}}{w_{13}+w_{23}}=:w^*$ and $x_3=\frac{w_{13}x_1+w_{23}x_2}{w_{13}+w_{23}}$, which from the invariant $x_1+x_2+x_3=0$ results in $x_2=-\frac{2w_{13}+w_{23}}{w_{13}+2w_{23}}x_1$. That is (generically) clustering.
	\item The equilibrium point $p=(x_1,x_2,x_3)$ defined by $p=\left\{ w=w^*\, | \, x_1=x_2= x_3=0 \right\}$ is non-hyperbolic.
\end{itemize}
\end{proposition}
\begin{proof}
	Straightforward computations.
\end{proof}

We emphasize that, just as in the symmetric case, the previous proposition provides the characterization of the critical manifold of \eqref{eq:sfs2} given as the union of ``the consensus'' and the ``the clustering'' manifolds. 

\begin{proposition}
	The layer equation of \eqref{eq:sfs2} has a transcritical bifurcation at $p$.
\end{proposition}
\begin{proof}
	We shall compute the eigenvalues associated to the equilibrium points of Proposition \ref{prop:eqs_ns}. We start by substituting $x_3=\frac{w_{13}x_1+w_{23}x_2}{w_{13}+w_{23}}$, and afterwards $x_2=-\frac{2w_{13}+w_{23}}{w_{13}+2w_{23}}x_1$. Naturally, both relations hold at all the equilibrium points. In particular, the equilibrium points along the consensus manifold are given by $x_1=0$, while those along the clustering manifold are given by $x_1=\frac{-(w_{13}+2w_{23})y}{ \alpha_1(w_{13}+2w_{23})-\alpha_2(2w_{13}+w_{23}) }$. By doing so we rewrite \eqref{eq:sfs2} in the limit $\ve\to0$ as
	\begin{equation}\label{eq:aux-red}
	\begin{split}
		x_1' &=  \frac{-3x_1(w_{13}+w_{13})( (w_{13}+2w_{23})y+x_1(\alpha_1(w_{13}+2w_{23})-\alpha_2(2w_{13}+w_{23})) )}{(w_{13}+2w_{23})^2},
	\end{split}
	\end{equation}
	where $y$ is considered as a parameter. Let $x_1'=\lambda_Nx_1$ denote the linearization of \eqref{eq:aux-red} at $x_1=0$, and $x_1'=\lambda_Mx_1$ denote the linearization of \eqref{eq:aux-red} at $x_1=\frac{-(w_{13}+2w_{23})y}{ \alpha_1(w_{13}+2w_{23})-\alpha_2(2w_{13}+w_{23}) }$. We have
	\begin{equation}
		\begin{split}
			\lambda_N &= -3\frac{w_{13}+w_{23}}{w_{13}+2w_{23}}y=-\lambda_M.
		\end{split}
	\end{equation}
	It is then clear that the eigenvalues along the consensus and the clustering manifolds ($\lambda_N$ and $\lambda_M$ respectively) have opposite signs and that there is an exchange in their signs at $y=0$. Finally we note that $y=0$ precisely corresponds to the point $p$, completing the proof.
\end{proof}

With the previous analysis we have shown that, at the singular level, the dynamics of the symmetric and the non-symmetric graphs are topologically equivalent. The reasons for this are: a) the uniform positive eigenvalue of the Laplacian matrix $L$ and b) that the dynamics of both systems are organized by a transcritical singularity corresponding to the intersection of the consensus and of the clustering manifolds. It remains to show that the passage through the transcritical singularity is also equivalent in both cases. We shall show this in the central chart $K_2$. We recall the the blow-up in the central chart is given by
\renewcommand{\ba}{a_2}
\renewcommand{\bb}{b_2}
\renewcommand{\bc}{c_2}
\renewcommand{\by}{y_2}
\renewcommand{\br}{r_2}
\begin{equation}
	x_1=\br\ba, \; x_2=\br\bb, \;x_3=\br\bc, \; y=\br\by, \; \ve=\br^2.
\end{equation}
Accordingly, the blown up vector field reads as
\begin{equation}\label{eq:ns-vf}
	\begin{split}
		\begin{bmatrix}
			\ba'\\ \bb'\\ \bc'
		\end{bmatrix} &=-\left( \begin{bmatrix}
			w^*+w_{13} & -w^* & -w_{13}\\
			-w^* & w^*+w_{23} & w_{23} \\
			-w_{13} & -w_{23} & w_{13} +w_{23}
		\end{bmatrix} +\br(\by+\alpha_1\ba+\alpha_2\bb)\begin{bmatrix}
			1 & -1 & 0 \\ -1 & 1 & 0\\ 0 & 0 & 0
		\end{bmatrix}\right)\begin{bmatrix}
			\ba\\ \bb\\ \bc
		\end{bmatrix}\\
		\by &= -\br(1+\cO(\br))
	\end{split}
\end{equation}
For $\br=0$ we get the linear system
\begin{equation}
	\begin{bmatrix}
			\ba'\\ \bb'\\ \bc'
		\end{bmatrix} =-\underbrace{\begin{bmatrix}
			w^*+w_{13} & -w^* & -w_{13}\\
			-w^* & w^*+w_{23} & -w_{23} \\
			-w_{13} & -w_{23} & w_{13} +w_{23}
		\end{bmatrix}}_{L_0}\begin{bmatrix}
			\ba\\ \bb\\ \bc
		\end{bmatrix}
\end{equation}

Using similar arguments as in Section \ref{sec:triangle} we can show that $\cW_2^\txtc = \left\{ \bc = \frac{w_{13}\ba+w_{23}\bb}{w_{13}+w_{23}} \right\}$ is an attracting centre manifold. Restriction to $\cW_2^\txtc$ and division of the vector field by $\br$ results in
\begin{equation}\label{eq:K2-ns-red}
	\begin{split}
		\ba' &= (\by+\alpha_1\ba+\alpha_2\bb)(\bb-\ba)\\
		\bb' &= (\by+\alpha_1\ba+\alpha_2\bb)(\ba-\ba)\\
		\by' &= -1.
	\end{split}
\end{equation}
Note that \eqref{eq:K2-ns-red} is the exact same equation \eqref{eq:K2-red} that we obtained in the symmetric case. This means that in both cases (symmetric and non-symmetric), the restriction to the centre manifold $\cW_2^\txtc$ coincide. With the above analysis we have shown the following:
\begin{proposition}
	The fast-slow system \eqref{eq:sfs2} with $w_{13}>0$ and $w_{23}>0$ is topologically equivalent to the case $w_{13}=w_{23}=1$.
\end{proposition}

In qualitative terms, the only difference between the symmetric and the non-symmetric cases is the convergence rate towards the invariant set $\cA$. Once trajectories have converged to $\cA$, the dynamics are organized \emph{by the same} transcritical singularity.

\subsection{Arbitrary graphs}\label{sec:arbitrary}

A natural question that arises is whether the analysis we have performed in Section \ref{sec:triangle}, particularly in Section \ref{sec:blowup}, has any relevance for consensus dynamics on arbitrary weighted graphs. In this section we argue that indeed, given some natural assumptions, generic consensus dynamics with one slowly varying weight behave essentially as the triangle motif.

Let us consider an undirected weighted graph $\cG=\left\{ \cV,\cE,\cW \right\}$. We assume that there are no self-loops, i.e. $w_{ii}=0$ for all $i\in1,\ldots,m$ and that there is at most one edge connecting two nodes. Next, without loss of generality, let us assume that $w_{12}=w\in\bbR$ is dynamic, while the rest of the weights $w_{ij}$, $(i,j)\in\cE$ and $(i,j)\neq(1,2)$, are fixed and positive. Denote by $L(w)$ the Laplacian matrix corresponding to the graphs we have defined so far. As in the analysis of the triangle motif we assume that the dynamic weight $w$ depends smoothly on the nodes it connects (namely $(x_1,x_2)$) and on an extrinsic slow variable $y\in\bbR$. Thus we consider the fast-slow system
\begin{equation}\label{eq:gen-consensus}
	\begin{split}
		x'&=-L(w)x,\\
		y' &= \ve g(x,y,\ve),
	\end{split}\qquad w=w(x_1,x_2,y,\ve),
\end{equation}
where $x\in\bbR^m$, $y\in\bbR$ and $L$ is a Laplacian matrix. Recall from the rescaling we performed in Section \ref{sec:blowup}, that it suffices to consider trajectories such that $\sum_{i=0}^m x_i(0)=0$, and to set $w=w^*+y+a_1 x_1+\alpha_2 x_2 +\cO(2)$, where by $\cO(2)$ we denote monomials of degree at least $2$. Moreover $w^*$ is a particular value at which $\dim(\ker L(w^*))>1$. 

\begin{lemma} Consider the consensus protocol \eqref{eq:gen-consensus} as defined above. Then, the following hold.
\begin{enumerate}[leftmargin=*]
	\item  Generically, if $\dim\ker(L(w))>1$, then $\dim\ker(L(w))=2$.
	\item Suppose that $L(w)$ is analytic in $w$, and that  $\dim\ker(L(w))=2$ at points $(x,y)$ defined by $w=w^*$. Then the consensus protocol \eqref{eq:gen-consensus} undergoes a (singular) transcritical bifurcation at $\left\{ w=w^* \right\}\cap\left\{ x=0\right\}$.
\end{enumerate}
	
\end{lemma}
\begin{proof}
	Recall that $0$ is a trivial eigenvalue with eigenvector $\bm 1_n$. Therefore, $\dim\ker L(w)$ increases its dimension whenever an eigenvalue crosses $0$. Next, we note that smooth one parameter families of symmetric matrices have \emph{generically} simple eigenvalues \cite{jonker1993partition}. This means that, generically, only one of the nontrivial eigenvalues of $L(w)$ can vanish for a certain value of $w$. Thus, the first item in the lemma follows. 

	Next, assuming analytic dependence of $L(w)$ on $w$ and according to \cite[II.6.2]{kato2013perturbation}, the matrix $L(w)$ can be orthogonally diagonalized as $L(w)=Q^\top(w)D(w)Q(w)$, where $Q(w)$ and $D(w)$ are also analytic in $w$. This means that the layer equation of \eqref{eq:gen-consensus} is conjugate to
	\begin{equation}
	\begin{split}
		z_1'&=0\\
		z_j'&=\lambda_j(\tilde w)z_j\\
		y'&=0,
	\end{split}\qquad j=2,\ldots,m,
	\end{equation}
where $\tilde w = \tilde w(z,y)= w^*+y+A(z)+\cO(2)$, where $A$ is a linear function with $A(0)=0$. We know from our previous arguments that there is a $k\in[2,3,\ldots,m]$ such that $\lambda_k(w^*)=0$ while $\lambda_j(w^*)\neq0$ for all $j\neq k$. Thus, we note that $z_k$ undergoes an exchange of stability if $\lambda_k$ crosses transversally the origin. But note that $\frac{\partial \lambda_k}{\partial y}(0)=\frac{\partial \lambda_k}{\partial \tilde w}\frac{\partial \tilde w}{\partial y}(0)$. Since we have that $\lambda_k$ depends analytically on $\tilde w$, we expect that $\frac{\partial \lambda_k}{\partial \tilde w}(0)$ is, generically, nonzero. Thus, we conclude that in a generic setting, whenever eigenvalues cross the origin, they do it transversally and thus a typical exchange of stability (transcritical bifurcation) occurs.
\end{proof}

\begin{remark}
	We note that singular pitchfork bifurcations can also occur, but this requires extra conditions on $w$, namely $w=y+\alpha_1x_1^2+\alpha_2 x_2^2$. In that case one can follow a similar analysis as that performed here; see also~\cite{krupa2001extendingtrans}.
\end{remark}

\biblio

%% file: subfiles/conclusions.tex
We have studied a class of adaptive networks under a linear average consensus protocol. The main property of the networks we studied here, is that one of the edges is slowly varying and can take values over the whole reals. The fact that the dynamic weight can be negative implies that the Laplacian of the graph may not be positive semi-definite, as is the case of non-negative weights. We have shown that the dynamics of the network are organized by a transcritical singularity. Interestingly, the network structure induces a generic maximal canard (unlike the fast-slow transcritical singularity of~\cite{krupa2001extendingtrans}). Moreover, we have shown that the blow-up method preserves the network structure. That is, on each of the blow-up charts we have found a related network but whose analysis is simpler compared to the full fast-slow one. Overall we have provided a case study in which center manifold reduction and tools of geometric singular perturbation theory, in particular the blow-up method, have been successfully used to describe the dynamics of a class of two time scale networks with a dynamic weight.

As we have considered one of the simplest network communication protocols, we conjecture that similar and more complicated problems are to be encountered when one studies general complex networks. For example, already assuming two slow weights presents new mathematical challenges: on one hand the critical manifold is potentially more complicated leading to more complex singular dynamics; on the other hand, it does not necessarily hold that the nonzero eigenvalues of the Laplacian are simple, and that they cross zero transversally. These two obstacles imply that, probably, one needs to desingularize higher codimension singularities and that generic results would be harder to obtain. Furthermore, although in this paper we have studied continuous-time and continuous state-space network dynamics, we expect that analogous multiscale phenomena can be found in discrete systems. Similar relevance can be expected in more general adaptive network scenarios including directed networks or those incorporating stochastic phenomena. Beyond our work here, there are currently relatively few works linking techniques between multiple time scale dynamical systems and network science~\cite{AshwinCreaserTsaneva,CappellettiWiuf,KuehnNetworks,KuehnNetworks1}. Yet, we conjecture that the inclusion of multiscale dynamics into network science will have far reaching consequences and will have a high impact in better modeling, analysis and understanding complex phenomena.

%% file: subfiles/num.tex
In this section we briefly discuss a numerical issue that may appear when simulating a fast-slow consensus network. This numerical issue is related to the presence of a maximal canard and yields trajectories not diverging when they should. Accordingly, we present an algorithm that overcomes, to some extent, such an issue. Afterwards, since the main part of the analysis of this paper concerns the triangle motif of Section \ref{sec:triangle}, we present some numerical examples of larger networks showcasing a transition through a transcritical singularity. In this way, we also present numerical evidence of the generalizations presented in Section \ref{sec:arbitrary}. At the end of this appendix, motivated by the fact that a dynamic weight in the consensus protocol opens-up the possibility for more complicated dynamics rather than simple average consensus, we numerically investigate the existence of periodic solutions. Although we do not rigorously study such a problem here, the geometric understanding of the way the dynamics are organized allows us to conjecture situations in which periodic solutions indeed appear. 

\subsection{An issue with numerical integration}

Due to the equivalences shown in Section \ref{sec:gen}, we shall restrict ourselves in this section to the triangle motif with $w_{13}=w_{23}=1$. For convenience we recall that the model reads as
\begin{equation}\label{eq:orig-num}
	\begin{split}
		\begin{bmatrix}
			a'\\ b'\\ c'
		\end{bmatrix} &=-\underbrace{\begin{bmatrix}
			w + 1 & -w & -1\\
			-w & w+1 & -1\\
			-1 & -1 & 2
		\end{bmatrix}}_{L(a,b,y)}\begin{bmatrix}
			a \\ b\\ c\\
		\end{bmatrix}\\
		y' &=\ve (-1+\cO(a,b,c,y,\ve))
	\end{split}
\end{equation}
where $w=-\frac{1}{2} + y + \alpha_1 a + \alpha_2 b$ and $a(t)+b(t)+c(t)=0$ for all $t\geq0$. We know that $L(a,b,y)$ has eigenvalues $\lambda_1=0$, $\lambda_2=3$, $\lambda_3=\lambda_3(a,b,y)$. Thus, roughly speaking, the numerical issue we present below has to do with which non-trivial eigenvalue dominates at $t=0$. Recall that the smallest non-trivial eigenvalue of the Laplacian, called in some instances ``spectral gap'' or ``algebraic connectivity'', bounds the rate of convergence towards consensus.

Let us discretize \eqref{eq:orig-num} using the \emph{forward Euler} approximation method $x'\approx\frac{x_{k+1}-x_k}{\Delta t}$, where $x_k=x(k\Delta t)$ with $k\in\mathbb N_0$. This discretization preserves the invariant $a(t)+b(t)+c(t)=0$ for all $t\geq0$,  resulting in $a_k+b_k+c_k=0$ for all $k\geq0$. Taking into account that $c_k=-a_k-b_k$, and disregarding the higher order terms in $y'$, we write the discretized version of \eqref{eq:orig-num} as
\begin{equation}\label{eq:num2}
	\begin{split}
		\begin{bmatrix}
			a_{k+1} \\ b_{k+1}
		\end{bmatrix} &= \underbrace{\begin{bmatrix}
			1-\Delta t(w_k+2) & \Delta t(w_k-1) \\ \Delta t(w_k-1) & 1-\Delta t(w_k+2)
		\end{bmatrix}}_{A_k(a_k,b_k,y_k)}\begin{bmatrix}
			a_k \\ b_k
		\end{bmatrix}\\{}
		y_{k+1} &= y_k-\ve\Delta t,
	\end{split}
\end{equation}
where $w_k=-\frac{1}{2} + y_k + \alpha_1 a_k + \alpha_2 b_k$. The matrix $A_k$ has spectrum $\spec A_k=\left\{\lambda_2, \lambda_3(k) \right\}=\left\{ 1-3\Delta t, 1-(2w_k+1)\Delta t \right\}$. From Section \ref{sec:prel_a}, we further know that the eigenvalue $\lambda_2$ is related to the convergence rate towards $\left\{c_k=\frac{a_k+b_k}{2}\right\}$, which together with the invariant $a_k+b_k+c_k=0$ is equivalent to the convergence towards $\left\{a_k=-b_k\right\}$. On the other hand, the eigenvalue $\lambda_3(k)$ is related to the stability of the consensus manifold, that is $\left\{a_k=b_k\right\}$. Since the aforementioned convergences are both exponential, and for meaningful\footnote{Those resulting in a transition from consensus to clustering, recall Theorem \ref{thm:main1}.} initial conditions both eigenvalues are within the unit circle for $k=0$, the relationship between the two eigenvalues plays an important role in the numerical integration.

Suppose that $|\lambda_3(0)|>M|\lambda_2|$, for some $M>0$. Then, up to truncation and computer precision errors, the relation $a_k=b_k$ dominates over $a_k=-b_k$ for some $k>K>0$. Substitution of $a_k=b_k$ in \eqref{eq:num2} implies that the governing difference equation is
\begin{equation}\label{eq:case1}
	a_{k+1} = (1-3\Delta t)a_k,
\end{equation}
where we observe that the weight $w_k$ does not play a role any more, and that $a_k=0$ is invariant. This means that, again disregarding numerical errors and approximations due to computer precision, $a_k=0$ is in this case a unique equilibrium point, \emph{independent} of the value of the weight $w_k$. In essence, this means that a computer algorithm may not recognize the instability of the consensus manifold; we note that this mechanism of ``over-stabilization'' induced by conserved quantities is different from recently discovered extended stabilization by the Euler method near transcritical and pitchfork singularities~\cite{ArcidiaconoEngelKuehn,engel2018discretized}.

On the contrary case, if $|\lambda_3(0)|<M|\lambda_2|$ for some $M>0$, then the relation $a_k=-b_k$ dominates and the governing difference equation is
\begin{equation}\label{eq:case2}
	a_{k+1} = (1-(2w_k+1)\Delta t)a_k,
\end{equation}
where the transition through a transcritical singularity, depending on the value of $w_k$, is present. Therefore, as expected from our analysis, there is a typical exchange of stability through a transcritical singularity for some negative values of $w_k$.

\begin{figure}[htbp]

\includegraphics{\main/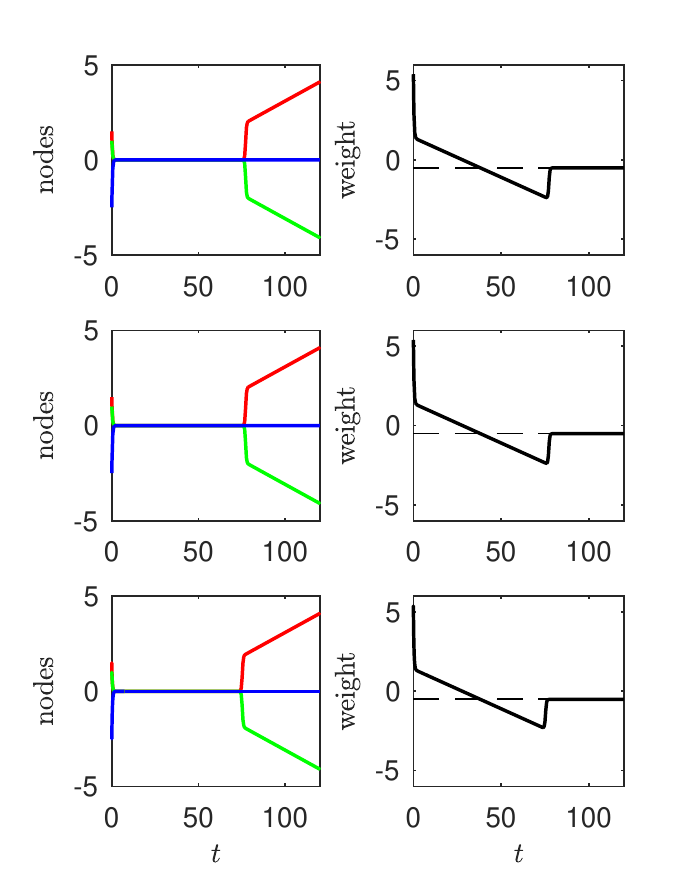}\hfill
\includegraphics{\main/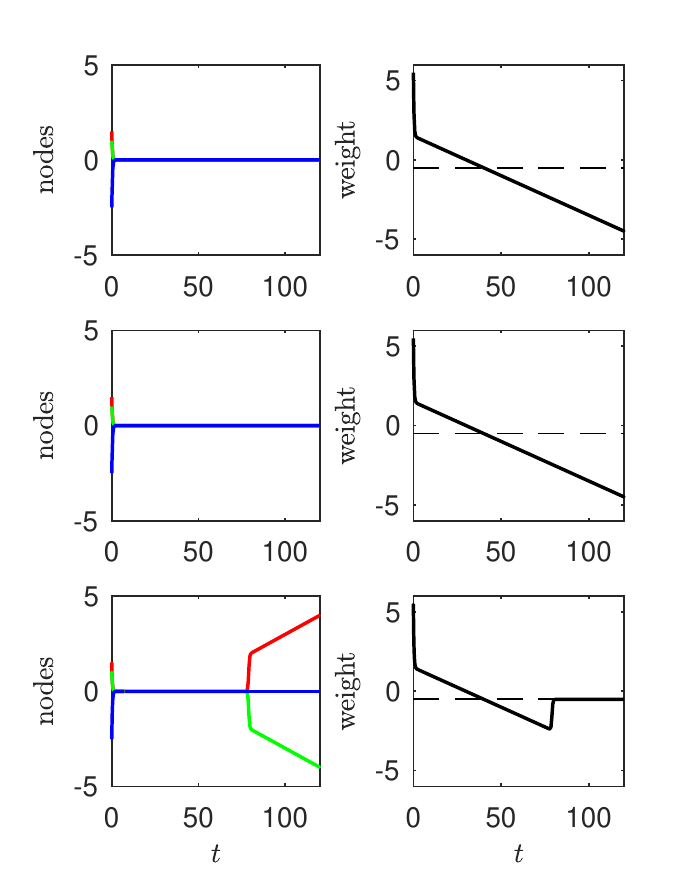}

\caption{Numerical simulation of \eqref{eq:orig-num} realized in Matlab 2017b. The values of the parameters are $(\alpha_1,\alpha_2,\ve)=(2,1,0.05)$ and the initial conditions of the nodes are $(a(0),b(0),c(0))=(1.5,1,-2.5)$. The plots on the left are for $y(0)=1.9$, while the plots on the right for $y(0)=2$. In both figures we show: a) Integration via forward Euler on the first row; b) Integration via {\tt{ode15s}} on the second row; and c) Integration via the proposed algorithm in the third row, where we set $e_\cA=1\times10^{-9}$. The Euler integrations are performed with $\Delta t=1\times 10^{-3}$ while for the {\tt{ode15s}} we set the {\tt{MaxStep}} to $1\times 10^{-3}$. We observe that for $y(0)=1.9$ the three algorithms provide similar outcome, namely the trajectories first approach consensus and then, when the weight is sufficiently negative, they transition towards clustering. The delay on the transition towards clustering is due to the trajectories being exponentially close to the maximal canard. However, for the case $y(0)=2$, only our proposed method provides a result aligned with the qualitative analysis of Section \ref{sec:triangle}. Note that the simulation for $y(0)=2$ falls under the case where after some time the dynamics are governed by \eqref{eq:case1}, explaining the observation of the trajectories remaining at (or close to) consensus even though the weight is largely negative. We recall that a negative weight means that the consensus manifold is unstable, yet trayectories remain close to it. In conclusion, our conjecture is that, in general, numerical algorithms for fast-slow consensus networks in the context of this paper work well when initial conditions are near the invariant set $\cA$ and $y(0)$ is not too large. Otherwise, numerical integration methods may require careful set-up in order to guarantee correct simulations.}
\label{fig:comparison}
\end{figure}

With the above simplified analysis we have argued that numerical integration of \eqref{eq:orig-num} may not be correct when $|2w(0)+1|>3M$ for some $M>0$. To prevent the aforementioned issue, one may adapt the analysis presented in Section \ref{sec:blowup} as shown in Algorithm \ref{alg1} and Figure \ref{fig:comparison}.

\begin{algorithm}
\KwData{$t_0$ (initial time), $\Delta t$ (time-step), $t_f$ (final time), $(a_0,b_0,c_0,y_0,\ve)$ (initial conditions and $\ve$), $e_\cA$ (a distance bound from the set $\cA=\left\{c_k=\frac{a_k+b_k}{2}\right\}$)}
\Begin{
	\For{$k=0,\ldots,\left[ \frac{t_f-t_0}{\Delta t} \right]$ }{

		\begin{equation*}
			\begin{split}
				\begin{bmatrix}
					a_{k+1} \\
					b_{k+1} \\
					c_{k+1}
				\end{bmatrix} &= (I -\Delta t L(a_k,b_k,y_k))\begin{bmatrix}
					a_{k} \\
					b_{k} \\
					c_{k}
				\end{bmatrix}\\
				y_{k+1} &= y_k-\ve\Delta t
			\end{split}
		\end{equation*}

		\If{$\left| c_k-\frac{a_k+b_k}{2} \right|\leq e_\cA$}{
			Stop and exit {\bfseries{for}} loop saving $(a_k,b_k,c_k,y_k,k)$.
		}		
	}
	\For{$j=k,\ldots,\left[ \frac{t_f-t_0}{\Delta t} \right]$}{
	\begin{equation*}
		\begin{split}
			a_{j+1} &= a_j-2\Delta t(y_j+(\alpha_1-\alpha_2)a_j)a_j\\
			b_{j+1} &= -a_{j+1}\\
			c_{j+1} &= 0\\
			y_{j+1} &= y_j-\ve\Delta t
		\end{split}
	\end{equation*}
	}
}
\caption{Pseudo-code used to simulate \eqref{eq:orig-num}. The first for-loop numerically integrates \eqref{eq:orig2} via \eqref{eq:orig-num}. This for-loop stops when the trajectories reach a small neighbourhood of the globally attracting set $\cA$. Such a threshold is set by $e_\cA$ in the algorithm. Afterwards we switch to a simplified system which is obtained as a restriction to $\cA$. Note the resemblance between the system in the second for-loop and \eqref{eq:K2-red}. Naturally, the actual numerical integration method is arbitrary, here for simplicity we have implemented forward Euler, yet providing reliable simulation results, see Figure \ref{fig:comparison}.}
\label{alg1}
\end{algorithm}

\subsection{Ring, complete and star networks}

\begin{figure}
\begin{subfigure}{0.5\textwidth}
	\includegraphics{\main/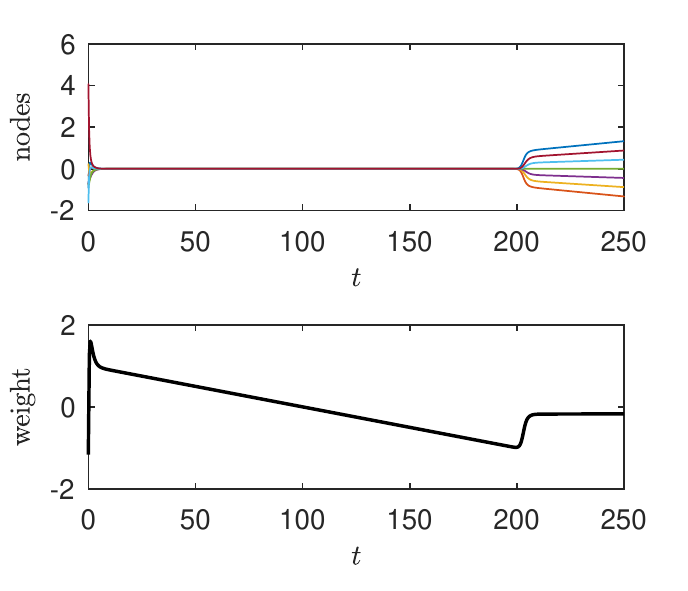}
	\caption{Ring graph with $7$ nodes}
\end{subfigure}%
\begin{subfigure}{0.5\textwidth}
	\includegraphics{\main/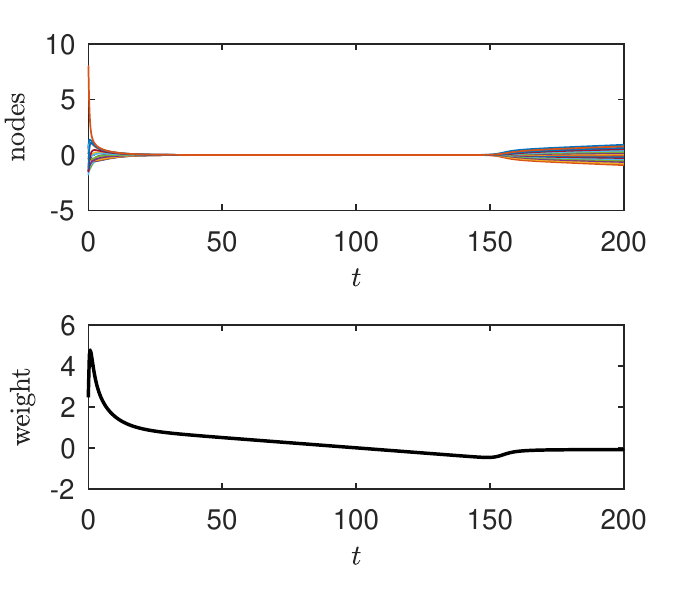}
	\caption{Ring graph with $16$ nodes}
\end{subfigure}\\[1ex]
\begin{subfigure}{0.5\textwidth}
	\includegraphics{\main/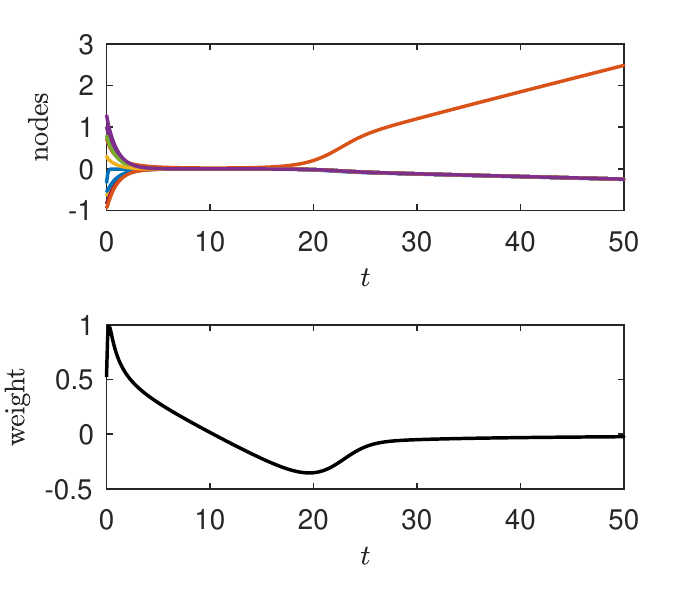}
	\caption{Star graph with $11$ nodes}
\end{subfigure}%
\begin{subfigure}{0.5\textwidth}
	\includegraphics{\main/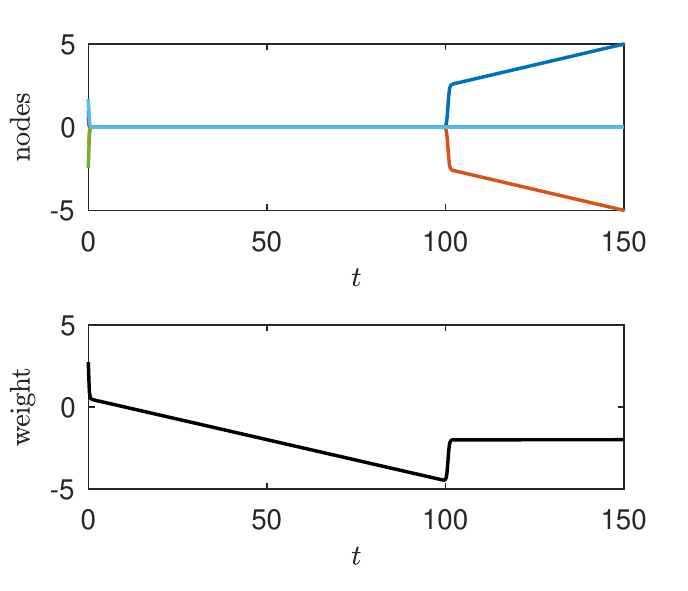}
	\caption{Complete graph with $6$ nodes}
\end{subfigure}
\caption{Simulation of a fast-slow consensus protocol \eqref{eq:gen-consensus} under several graph topologies. For all these simulations we consider the dynamic weight $w=y+2x_1+x_2$ on the edge $\left\{ 1,2 \right\}$ and all other weights fixed to $1$. We also set $\ve=0.05$ and $y(0)=0.5$. 
We observe in all cases that the trajectories first converge towards the origin (average consensus), and then, when $w$ reaches some negative value, the trajectories leave a small vicinity of the consensus manifold and approach a clustering state. The actual convergence rate in each case is given by the so-called spectral gap, or the smallest non-trivial eigenvalue.
}
\label{fig:gen-graphs}
\end{figure}

To complement the analysis performed in Section \ref{sec:gen}, we present here a couple of numerical simulations of fast-slow consensus protocols defined on ring, complete and star networks with more than $3$ nodes. In all the simulations shown in Figure \ref{fig:gen-graphs} we keep the dynamic weight $w=y+\alpha_1x_1+\alpha_2x_2$ between nodes $x_1$ and $x_2$ while all other edges are fixed to $1$. We note that we produce simulations for initial conditions that lead to exchange between consensus and clustering. We observe that in all cases of Figure \ref{fig:gen-graphs}, we have a qualitative behaviour similar to the one analyzed in the triangle motif, thus validating the arguments of Section \ref{sec:gen}.

\subsection{Periodic orbits on the triangle motif}

Commonly, and historically, one considers consensus protocols seeking consensus among the nodes. This behaviour is the only possible one if the weights of the graph are non-negative. However, as we have seen in the main part of this paper, introducing negative weights enriches the dynamics. In this section we want to provide numerical evidence on the existence of periodic orbits in consensus networks. We do not attempt to give a full treatment of the problem. Rather, we present two instances that are constructed from the geometric description we developed in Section \ref{sec:triangle}.

Let us then consider the triangle motif \eqref{eq:orig-num} and let us consider initial conditions with $y(0)>0$. In order to produce periodic orbits one must introduce a return mechanism that allows trajectories to return to $y>0$ after they have crossed the origin. Since we already know that the dynamics are organized as sketched in Figure \ref{fig:crit0}, we propose two return mechanisms.

\begin{enumerate}[leftmargin=*]
	\item The first natural way to introduce a return mechanism is to add a drift on the slow variable $y$ that acts away from the consensus manifold $\cN_\ve$. For this strategy to result in periodic orbits, one further requires that the clustering manifold $\cM_0$ is aligned with the fast-foliation, that is $\alpha_1-\alpha_2=0$. For this example we them propose
	\begin{equation}\label{eq:per1}
		\begin{split}
		\begin{bmatrix}
			a'\\ b'\\ c'
		\end{bmatrix} &=-\begin{bmatrix}
			w + 1 & -w & -1\\
			-w & w+1 & -1\\
			-1 & -1 & 2
		\end{bmatrix}\begin{bmatrix}
			a \\ b\\ c\\
		\end{bmatrix},\qquad w=-\frac{1}{2} + y + \alpha(a+b),\\
		y' &=\ve (-1 +k_1a_1^2 )
	\end{split}
	\end{equation}

We emphasize that in this case, the existence of periodic solutions does not follow from the singular limit. However, we know from our geometric analysis, especially in the chart $K_2$, that this model indeed has periodic solutions. We show in Figure \ref{fig:sim-periodic1} a simulation of \eqref{eq:per1}.

\begin{figure}[htbp]\centering
	\includegraphics{\main/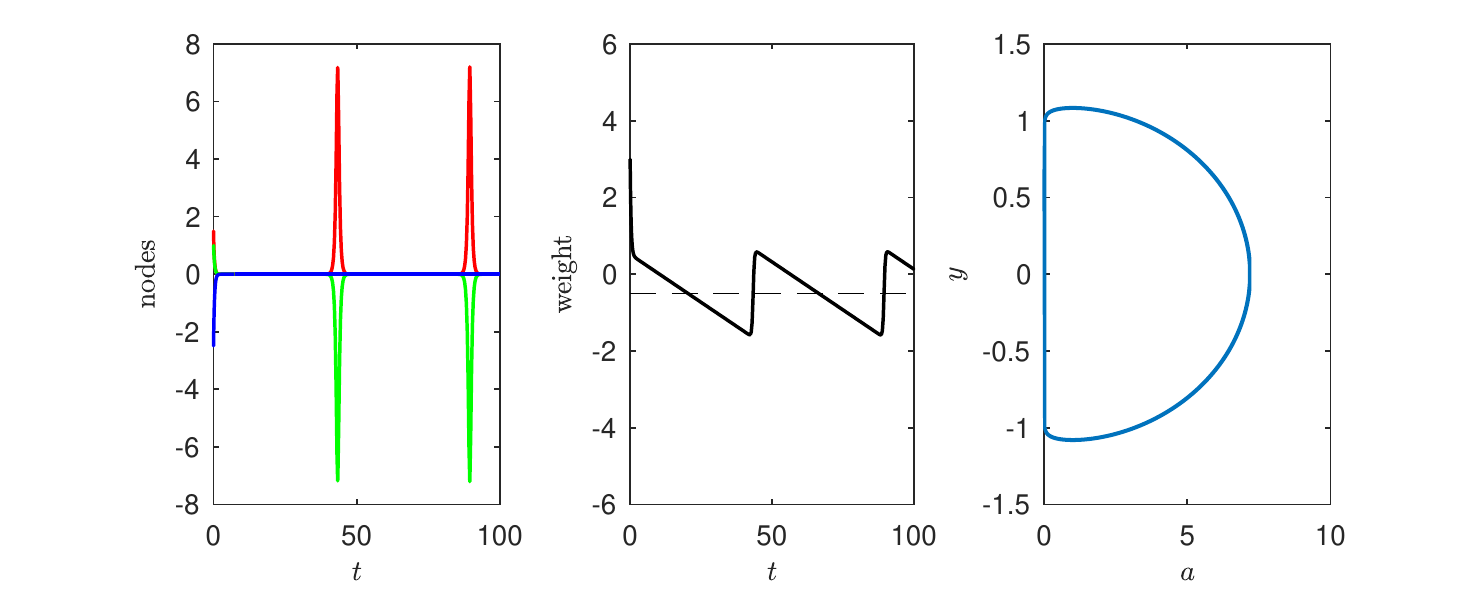}
	\caption{Simulation of \eqref{eq:per1} with $(y(0),\alpha,k_1,\ve)=(0.5,1,1,0.05)$. From left to right, the first two plots are the time series for the nodes and the weight, while the third plot depicts the corresponding phase portrait.}
	\label{fig:sim-periodic1}
\end{figure}

\item Another way to produce a return mechanism is to add higher order terms to the dynamic weight. Based on the singular limit sketched in Figure \ref{fig:crit0}, the idea is to achieve a singular limit as depicted in Figure \ref{fig:sing-return2}.

\begin{figure}[htbp]\centering

	\begin{tikzpicture}
		\node at (0,0){
		\includegraphics{\main/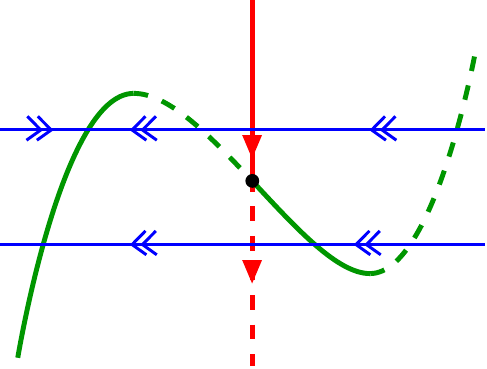}
		};
		\node[red] at (0,2.1) {$\cN_0$};
		\node[green!50!black] at (2.2,1.6) {$\cM_0$};
		\draw[->] (3,-2)--++(0.5,0) node[right]{$a$};
		\draw[->] (3,-2)--++(0,0.5) node[above]{$y$};
	\end{tikzpicture}
	\caption{Singular limit for a return mechanism via higher order terms in the weight $w$.}
	\label{fig:sing-return2}
\end{figure}

We note that near the origin, the singular limit is just as in Figure \ref{fig:crit0}. Furthermore, by inspecting Figure \ref{fig:sing-return2} we observe that returning on the left side of $\cM_0$ is ``easier'' than on the right, because, to produce a return mechanism on the right side we would actually require to generate a canard cycle. So, for the purposes of this example we shall concentrate on the case where the dynamics with $\ba<0$ exhibit periodic solutions. Following Figure \ref{fig:sing-return2} we propose the weight
\begin{equation}
	w=-\frac{1}{2} + \alpha_1 a + \alpha_2 b + \alpha_3 a^3 + \alpha_4 b^3,
\end{equation}
with $\alpha_1-\alpha_2>0$ and $\alpha_3-\alpha_4<0$. Next, we must introduce a new slow vector field so that trajectories can travel back along the left branch of the cubic critical manifold. Naturally there are many ways to achieve this. For example:
\begin{equation}\label{eq:per2}
		\begin{split}
		\begin{bmatrix}
			a'\\ b'\\ c'
		\end{bmatrix} &=-\begin{bmatrix}
			w + 1 & -w & -1\\
			-w & w+1 & -1\\
			-1 & -1 & 2
		\end{bmatrix}\begin{bmatrix}
			a \\ b\\ c\\
		\end{bmatrix},\qquad w=-\frac{1}{2} + \alpha_1 a + \alpha_2 b + \alpha_3 a^3 + \alpha_4 b^3,\\
		y' &=\ve (-1 + \beta_1 a),
	\end{split}
	\end{equation}
for some $\beta<0$. Figure \ref{fig:sim-periodic2} shows a corresponding simulation.

\begin{figure}[htbp]\centering
	\includegraphics{\main/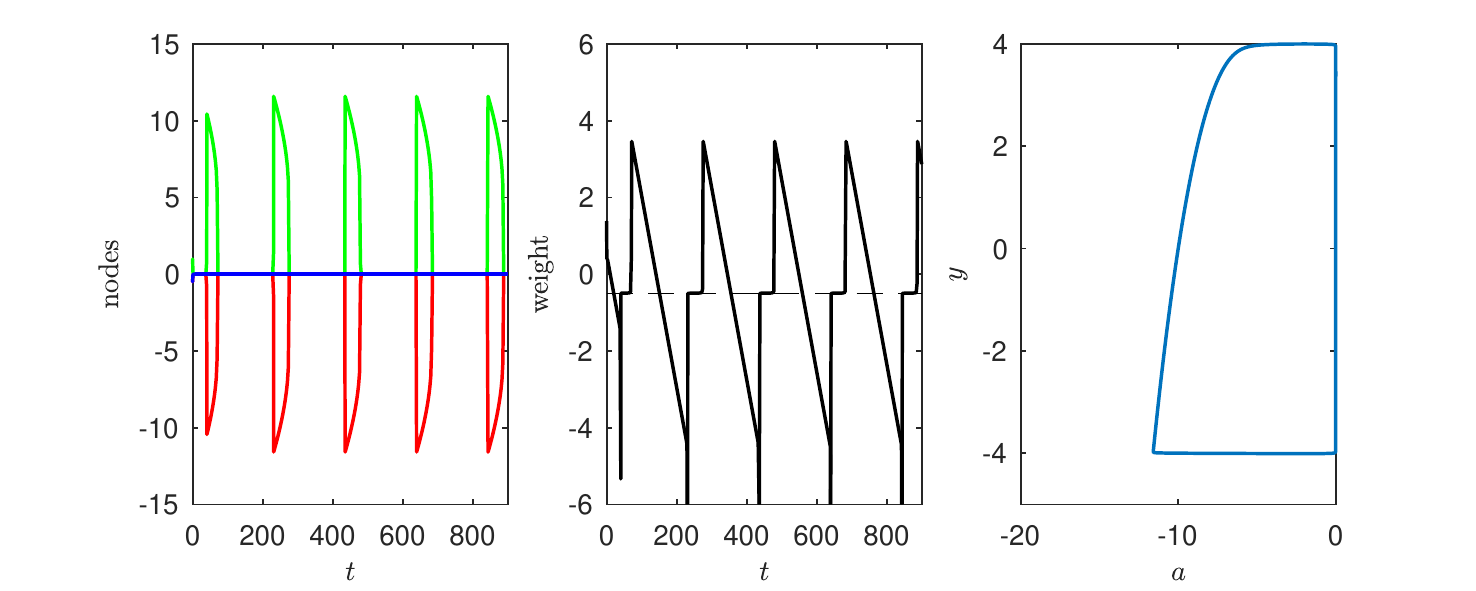}
	\caption{An example of periodic orbits on the triangle motif due to a cubic weight. We have used parameters $(\alpha_1,\alpha_2,\alpha_3,\alpha_4,\beta,\ve)=(2,1,1,1.01,-.5,.05)$.}
	\label{fig:sim-periodic2}
\end{figure}

\end{enumerate}

\biblio

%% file: main.bbl
\begin{thebibliography}{10}

\bibitem{abramowitz1972handbook}
M.~Abramowitz and I.~A. Stegun.
\newblock {\em Handbook of Mathematical Functions}.
\newblock Dover, 1972.

\bibitem{albert2002statistical}
R.~Albert and A.-L. Barab{\'a}si.
\newblock Statistical mechanics of complex networks.
\newblock {\em Rev. Mod. Phys.}, 74(1):47, 2002.

\bibitem{altafini2013consensus}
C.~Altafini.
\newblock Consensus problems on networks with antagonistic interactions.
\newblock {\em {IEEE} Trans. Autom. Control}, 58(4):935--946, 2013.

\bibitem{alves2007unveiling}
N.~A. Alves.
\newblock Unveiling community structures in weighted networks.
\newblock {\em Phys. Rev. E}, 76(3):036101, 2007.

\bibitem{aral2009distinguishing}
S.~Aral, L.~Muchnik, and A.~Sundararajan.
\newblock Distinguishing influence-based contagion from homophily-driven
  diffusion in dynamic networks.
\newblock {\em Proc. Natl. Acad. Sci.}, 106(51):21544--21549, 2009.

\bibitem{ArcidiaconoEngelKuehn}
L.~Arcidiacono, M.~Engel, and C.~Kuehn.
\newblock Discretized fast-slow systems near pitchfork singularities.
\newblock {\em arXiv:1902.06512}, pages 1--29, 2019.

\bibitem{AshwinCreaserTsaneva}
P.~Ashwin, J.~Creaser, and K.~Tsaneva-Atanasova.
\newblock Fast and slow domino regimes in transient network dynamics.
\newblock {\em Phys. Rev. E}, 96(5):052309, 2017.

\bibitem{awad2018time}
A.~Awad, A.~Chapman, E.~Schoof, A.~Narang-Siddarth, and M.~Mesbahi.
\newblock Time-scale separation in networks: State-dependent graphs and
  consensus tracking.
\newblock {\em {IEEE} Trans. Control Netw. Syst.}, 2018.

\bibitem{barabasi2016network}
A.-L. Barab{\'a}si.
\newblock {\em Network science}.
\newblock CUP, 2016.

\bibitem{barrat2004architecture}
A.~Barrat, M.~Barthelemy, R.~Pastor-Satorras, and A.~Vespignani.
\newblock The architecture of complex weighted networks.
\newblock {\em Proc. Natl. Acad. Sci.}, 101(11):3747--3752, 2004.

\bibitem{BarratBarthelemyVespignani}
A.~Barrat, M.~Barth\'{e}lemy, and A.~Vespignani.
\newblock {\em Dynamical Processes on Complex Networks}.
\newblock CUP, 2008.

\bibitem{boccaletti2006complex}
S.~Boccaletti, V.~Latora, Y.~Moreno, M.~Chavez, and D-U. Hwang.
\newblock Complex networks: Structure and dynamics.
\newblock {\em Phys. Rep.}, 424(4-5):175--308, 2006.

\bibitem{bronski2014spectral}
J.~C. Bronski and L.~DeVille.
\newblock Spectral theory for dynamics on graphs containing attractive and
  repulsive interactions.
\newblock {\em SIAM J. App. Math.}, 74(1):83--105, 2014.

\bibitem{CappellettiWiuf}
D.~Cappelletti and C.~Wiuf.
\newblock Uniform approximation of solutions by elimination of intermediate
  species in deterministic reaction networks.
\newblock {\em SIAM J. Appl. Dyn. Syst.}, 16(4):2259--2286, 2017.

\bibitem{casteigts2012time}
A.~Casteigts, P.~Flocchini, W.~Quattrociocchi, and N.~Santoro.
\newblock Time-varying graphs and dynamic networks.
\newblock {\em Int. J. Parallel Emergent Distrib. Syst.}, 27(5):387--408, 2012.

\bibitem{chen2016semidefiniteness}
W.~Chen, D.~Wang, J.~Liu, T.~Ba{\c{s}}ar, K.~H. Johansson, and L.~Qiu.
\newblock On semidefiniteness of signed laplacians with application to
  microgrids.
\newblock {\em IFAC-PapersOnLine}, 49(22):97--102, 2016.

\bibitem{chen2013programmable}
Y-J. Chen, N.~Dalchau, N.~Srinivas, A.~Phillips, L.~Cardelli, D.~Soloveichik,
  and G.~Seelig.
\newblock Programmable chemical controllers made from {DNA}.
\newblock {\em Nat. Nanotechnol.}, 8(10):755, 2013.

\bibitem{de2015planar}
P.~De~Maesschalck.
\newblock Planar canards with transcritical intersections.
\newblock {\em Acta App. Math.}, 137(1):159--184, 2015.

\bibitem{dumortier1996canard}
F.~Dumortier and R.~Roussarie.
\newblock {\em Canard cycles and center manifolds}, volume 577.
\newblock Amer. Math. Soc., 1996.

\bibitem{eckhaus2011asymptotic}
W.~Eckhaus.
\newblock {\em Asymptotic analysis of singular perturbations}, volume~9.
\newblock Elsevier, 2011.

\bibitem{eckhaus2011matched}
W.~Eckhaus.
\newblock {\em Matched asymptotic expansions and singular perturbations},
  volume~6.
\newblock Elsevier, 2011.

\bibitem{engel2018discretized}
M.~Engel and C.~Kuehn.
\newblock {Discretized Fast-Slow Systems near Transcritical Singularities}.
\newblock {\em arXiv preprint arXiv:1806.06561}, pages 1--, 2018.

\bibitem{fax2003information}
J.~A. {Fax} and R.~M. {Murray}.
\newblock Information flow and cooperative control of vehicle formations.
\newblock {\em {IEEE} Trans. Autom. Control}, 49(9):1465--1476, Sep. 2004.

\bibitem{fenichel1979geometric}
N.~Fenichel.
\newblock Geometric singular perturbation theory for ordinary differential
  equations.
\newblock {\em J. Differ. Equ.}, 31(1):53--98, 1979.

\bibitem{Field}
M.~Field.
\newblock Combinatorial dynamics.
\newblock {\em Dynam. Syst.}, 19(3):217--243, 2004.

\bibitem{GolubitskyStewart2}
M.~Golubitsky and I.~Stewart.
\newblock Coordinate changes for network dynamics.
\newblock {\em Dynam. Syst.}, 32(1):80--116, 2017.

\bibitem{GrossSayama}
T.~Gross and H.~Sayama, editors.
\newblock {\em Adaptive Networks: Theory, Models and Applications}.
\newblock Springer, 2009.

\bibitem{holland2004using}
B.~R. Holland, K.~T. Huber, V.~Moulton, and P.~J. Lockhart.
\newblock Using consensus networks to visualize contradictory evidence for
  species phylogeny.
\newblock {\em Mol. Biol. Evol.}, 21(7):1459--1461, 2004.

\bibitem{jadbabaie2003coordination}
A.~{Jadbabaie}, , and A.~S. {Morse}.
\newblock Coordination of groups of mobile autonomous agents using nearest
  neighbor rules.
\newblock {\em IEEE Trans. Autom. Control}, 48(6):988--1001, June 2003.

\bibitem{jardon2019survey}
H.~Jard\'on-Kojakhmetov and C.~Kuehn.
\newblock A survey on the blow-up method for fast-slow systems.
\newblock {\em arXiv preprint arXiv:1901.01402}, pages 1--, 2019.

\bibitem{jones1995geometric}
C.~K. R.~T. Jones.
\newblock Geometric singular perturbation theory.
\newblock In {\em Dynamical Systems}, pages 44--118. Springer, 1995.

\bibitem{jonker1993partition}
P.~Jonker, G.~Still, and F.~Twilt.
\newblock {\em On the partition of real symmetric matrices according to the
  multiplicities of their eigenvalues}.
\newblock Fac. of Applied Math., University of Twente, 1993.

\bibitem{Kaper}
T.J. Kaper.
\newblock An introduction to geometric methods and dynamical systems theory for
  singular perturbation problems.
\newblock In J.~Cronin and R.E. O'Malley, editors, {\em Analyzing multiscale
  phenomena using singular perturbation methods}, pages 85--131. Springer,
  1999.

\bibitem{kato2013perturbation}
T.~Kato.
\newblock {\em Perturbation theory for linear operators}, volume 132.
\newblock Springer Science \& Business Media, 2013.

\bibitem{knyazev2017signed}
A.~V. Knyazev.
\newblock Signed laplacian for spectral clustering revisited.
\newblock {\em arXiv preprint arXiv:1701.01394}, 1, 2017.

\bibitem{krupa2001extending}
M.~Krupa and P.~Szmolyan.
\newblock Extending geometric singular perturbation theory to nonhyperbolic
  points---fold and canard points in two dimensions.
\newblock {\em SIAM J. Math. Anal.}, 33(2):286--314, 2001.

\bibitem{krupa2001extendingtrans}
M.~Krupa and P.~Szmolyan.
\newblock Extending slow manifolds near transcritical and pitchfork
  singularities.
\newblock {\em Nonlinearity}, 14(6):1473, 2001.

\bibitem{KuehnNetworks}
C.~Kuehn.
\newblock Time-scale and noise optimality in self-organized critical adaptive
  networks.
\newblock {\em Phys. Rev. E}, 85(2):026103, 2012.

\bibitem{kuehn2015multiple}
C.~Kuehn.
\newblock {\em Multiple time scale dynamics}, volume 191.
\newblock Springer, 2015.

\bibitem{KuehnNetworks1}
C.~Kuehn.
\newblock Multiscale dynamics of an adaptive catalytic network model.
\newblock {\em Math. Model. Nat. Pheno.}, pages 1--, 2019.
\newblock accepted / to appear.

\bibitem{kuehn2015multiscale}
C.~Kuehn and P.~Szmolyan.
\newblock Multiscale geometry of the {O}lsen model and non-classical relaxation
  oscillations.
\newblock {\em J. Nonlinear Sci.}, 25(3):583--629, 2015.

\bibitem{kunkel2006differential}
P.~Kunkel and V.~Mehrmann.
\newblock {\em Differential-algebraic equations: analysis and numerical
  solution}, volume~2.
\newblock European Mathematical Society, 2006.

\bibitem{lynch1996distributed}
N.~A. Lynch.
\newblock {\em Distributed algorithms}.
\newblock Elsevier, 1996.

\bibitem{merris1994laplacian}
R.~Merris.
\newblock Laplacian matrices of graphs: a survey.
\newblock {\em Linear Algebra Appl.}, 197:143--176, 1994.

\bibitem{mesbahi2005state}
M.~Mesbahi.
\newblock On state-dependent dynamic graphs and their controllability
  properties.
\newblock {\em {IEEE} Trans. Autom. Control}, 50(3):387--392, 2005.

\bibitem{mesbahi2010graph}
M.~Mesbahi and M.~Egerstedt.
\newblock {\em Graph theoretic methods in multiagent networks}, volume~33.
\newblock Princeton University Press, 2010.

\bibitem{milo2002network}
R.~Milo, S.~Shen-Orr, S.~Itzkovitz, N.~Kashtan, D.~Chklovskii, and U.~Alon.
\newblock Network motifs: simple building blocks of complex networks.
\newblock {\em Science}, 298(5594):824--827, 2002.

\bibitem{mohar1991laplacian}
B.~Mohar.
\newblock The laplacian spectrum of graphs.
\newblock In {\em Graph Theory, Combinatorics, and Applications}, pages
  871--898. Wiley, 1991.

\bibitem{moreau2004stability}
L.~Moreau.
\newblock Stability of continuous-time distributed consensus algorithms.
\newblock In {\em Proceedings of the 43rd IEEE Conference on Decision and
  Control}, volume~4, pages 3998--4003, 2004.

\bibitem{moreau2005stability}
L.~Moreau.
\newblock Stability of multiagent systems with time-dependent communication
  links.
\newblock {\em {IEEE} Trans. Autom. Control}, 50(2):169--182, 2005.

\bibitem{moreno1934shall}
J.~L. Moreno.
\newblock {\em Who shall survive?: A new approach to the problem of human
  interrelations.}
\newblock Nerv. \& Ment. Dis. Publ. Co., 1934.

\bibitem{NijholtRinkSanders}
E.~Nijholt, B.~Rink, and J.~Sanders.
\newblock Center manifolds of coupled cell networks.
\newblock {\em SIAM J. Math. Anal.}, 49(5):4117--4148, 2017.

\bibitem{olfati2005distributed}
R.~Olfati-Saber.
\newblock Distributed kalman filter with embedded consensus filters.
\newblock In {\em Proceedings of the 44th IEEE Conference on Decision and
  Control}, pages 8179--8184, 2005.

\bibitem{olfati2007consensus}
R.~Olfati-Saber, J.~A. Fax, and R.~M. Murray.
\newblock Consensus and cooperation in networked multi-agent systems.
\newblock {\em Proc. {IEEE}}, 95(1):215--233, 2007.

\bibitem{olfati2004consensus}
R.~Olfati-Saber and R.~M. Murray.
\newblock Consensus problems in networks of agents with switching topology and
  time-delays.
\newblock {\em {IEEE} Trans. Autom. Control}, 49(9):1520--1533, 2004.

\bibitem{OMalley1991}
R.~E. O'Malley, Jr.
\newblock {\em Singular perturbation methods for ordinary differential
  equations}, volume~89 of {\em Applied Mathematical Sciences}.
\newblock Springer-Verlag, New York, 1991.

\bibitem{pan2016laplacian}
L.~Pan, H.~Shao, and M.~Mesbahi.
\newblock Laplacian dynamics on signed networks.
\newblock In {\em Proceedings of the 55th Conference on Decision and Control},
  pages 891--896, 2016.

\bibitem{pastor2001epidemic}
R.~Pastor-Satorras and A.~Vespignani.
\newblock Epidemic dynamics and endemic states in complex networks.
\newblock {\em Phys. Rev. E}, 63(6):066117, 2001.

\bibitem{proskurnikov2013average}
A.~V. Proskurnikov.
\newblock Average consensus in networks with nonlinearly delayed couplings and
  switching topology.
\newblock {\em Automatica}, 49(9):2928--2932, 2013.

\bibitem{proskurnikov2016opinion}
A.~V. Proskurnikov, A.~S. Matveev, and M.~Cao.
\newblock Opinion dynamics in social networks with hostile camps: Consensus vs.
  polarization.
\newblock {\em {IEEE} Trans. Autom. Control}, 61(6):1524--1536, 2016.

\bibitem{proskurnikov2017tutorial}
A.~V. Proskurnikov and R.~Tempo.
\newblock A tutorial on modeling and analysis of dynamic social networks. part
  i.
\newblock {\em Annual Reviews in Control}, 43:65--79, 2017.

\bibitem{proskurnikov2018tutorial}
A.~V. Proskurnikov and R.~Tempo.
\newblock A tutorial on modeling and analysis of dynamic social networks. part
  ii.
\newblock {\em Annual Reviews in Control}, 2018.

\bibitem{ren2008distributed}
W.~Ren and R.~W. Beard.
\newblock {\em Distributed consensus in multi-vehicle cooperative control}.
\newblock Springer, 2008.

\bibitem{ren2005survey}
W.~Ren, R.~W. Beard, and E.~M. Atkins.
\newblock A survey of consensus problems in multi-agent coordination.
\newblock In {\em Proceedings of the American Control Conference, 2005.}, pages
  1859--1864.

\bibitem{ren2007information}
W.~Ren, R.~W. Beard, and E.~M. Atkins.
\newblock Information consensus in multivehicle cooperative control.
\newblock {\em {IEEE} Control Syst. Mag.}, 27(2):71--82, 2007.

\bibitem{saber2003consensus}
R.~O. {Saber} and R.~M. {Murray}.
\newblock Consensus protocols for networks of dynamic agents.
\newblock In {\em Proceedings of the American Control Conference, 2003.},
  volume~2, pages 951--956.

\bibitem{schweitzer2009economic}
F.~Schweitzer, G.~Fagiolo, D.~Sornette, F.~Vega-Redondo, A.~Vespignani, and
  D.~R. White.
\newblock Economic networks: The new challenges.
\newblock {\em Science}, 325(5939):422--425, 2009.

\bibitem{strogatz2001exploring}
S.~H. Strogatz.
\newblock Exploring complex networks.
\newblock {\em Nature}, 410(6825):268, 2001.

\bibitem{takens1976constrained}
F.~Takens.
\newblock Constrained equations; a study of implicit differential equations and
  their discontinuous solutions.
\newblock In {\em Structural stability, the theory of catastrophes, and
  applications in the sciences}, pages 143--234. Springer, 1976.

\bibitem{tanner2007flocking}
H.~G. Tanner, A.~Jadbabaie, and G.~J. Pappas.
\newblock Flocking in fixed and switching networks.
\newblock {\em {IEEE} Trans. Autom. Control}, 52(5):863--868, 2007.

\bibitem{thomas1977majority}
R.~H. Thomas.
\newblock A majority consensus approach to concurrency control for multiple
  copy databases.
\newblock {\em ACM Trans. Database Syst.}, 4(2):180--209, jun 1979.

\bibitem{tikhonov1952systems}
A.~N. Tikhonov.
\newblock Systems of differential equations containing small parameters in the
  derivatives.
\newblock {\em Mat. Sb.}, 73(3):575--586, 1952.

\bibitem{vanderHofstad}
R.~van~den Hofstad.
\newblock {\em Random Graphs and Complex Networks: Volume 1}.
\newblock CUP, 2016.

\bibitem{verhulst2005methods}
F.~Verhulst.
\newblock {\em Methods and applications of singular perturbations: boundary
  layers and multiple timescale dynamics}, volume~50.
\newblock Springer Science \& Business Media, 2005.

\bibitem{xia2004analyzing}
Y.~Xia, H.~Yu, R.~Jansen, M.~Seringhaus, S.~Baxter, D.~Greenbaum, H.~Zhao, and
  M.~Gerstein.
\newblock Analyzing cellular biochemistry in terms of molecular networks.
\newblock {\em Annu. Rev. Biochem.}, 73(1):1051--1087, 2004.

\bibitem{xiao2007distributed}
L.~Xiao, S.~Boyd, and S-J. Kim.
\newblock Distributed average consensus with least-mean-square deviation.
\newblock {\em J. Parallel Distr. Com.}, 67(1):33--46, 2007.

\bibitem{xie2011social}
J.~Xie, S.~Sreenivasan, G.~Korniss, W.~Zhang, C.~Lim, and B.~K. Szymanski.
\newblock Social consensus through the influence of committed minorities.
\newblock {\em Phys. Rev. E}, 84(1):011130, 2011.

\bibitem{zelazo2014definiteness}
D.~Zelazo and M.~B{\"u}rger.
\newblock On the definiteness of the weighted laplacian and its connection to
  effective resistance.
\newblock In {\em 53rd IEEE Conference on Decision and Control}, pages
  2895--2900, 2014.

\bibitem{zelazo2017robustness}
D.~Zelazo and M.~B{\"u}rger.
\newblock On the robustness of uncertain consensus networks.
\newblock {\em {IEEE} Trans. Control Netw. Syst.}, 4(2):170--178, 2017.

\end{thebibliography}


\newcommand{\etalchar}[1]{$^{#1}$}
\begin{thebibliography}{ACS{\etalchar{+}}18}

\bibitem[ACS{\etalchar{+}}18]{awad2018time}
Armand Awad, Airlie Chapman, Eric Schoof, Anshu Narang-Siddarth, and Mehran
  Mesbahi.
\newblock Time-scale separation in networks: State-dependent graphs and
  consensus tracking.
\newblock {\em IEEE Transactions on Control of Network Systems}, 2018.

\end{thebibliography}


\begin{thebibliography}{JKK19}

\bibitem[AS72]{abramowitz1972handbook}
M.~Abramowitz and I.A. Stegun.
\newblock {\em Handbook of Mathematical Functions}.
\newblock Dover, 1972.

\bibitem[DR96]{dumortier1996canard}
F.~Dumortier and R.~Roussarie.
\newblock {\em Canard cycles and center manifolds}, volume 577.
\newblock American Mathematical Society, 1996.

\bibitem[JKK19]{jardon2019survey}
H.~Jard\'on-Kojakhmetov and C.~Kuehn.
\newblock A survey on the blow-up method for fast-slow systems.
\newblock {\em arXiv preprint arXiv:1901.01402}, 2019.

\bibitem[KS01]{krupa2001extending}
M.~Krupa and P.~Szmolyan.
\newblock Extending geometric singular perturbation theory to nonhyperbolic
  points---fold and canard points in two dimensions.
\newblock {\em SIAM journal on mathematical analysis}, 33(2):286--314, 2001.

\bibitem[Kue15]{kuehn2015multiple}
C.~Kuehn.
\newblock {\em Multiple time scale dynamics}, volume 191.
\newblock Springer, 2015.

\end{thebibliography}


\begin{thebibliography}{}

\end{thebibliography}


\begin{thebibliography}{AS72}

\bibitem[AS72]{abramowitz1972handbook}
M.~Abramowitz and I.A. Stegun.
\newblock {\em Handbook of Mathematical Functions}.
\newblock Dover, 1972.

\end{thebibliography}


\begin{thebibliography}{JST93}

\bibitem[JST93]{jonker1993partition}
Peter Jonker, Georg Still, and Frank Twilt.
\newblock {\em On the partition of real symmetric matrices according to the
  multiplicities of their eigenvalues}.
\newblock Fac. of Applied Math., University of Twente, 1993.

\bibitem[Kat13]{kato2013perturbation}
Tosio Kato.
\newblock {\em Perturbation theory for linear operators}, volume 132.
\newblock Springer Science \& Business Media, 2013.

\bibitem[KS01]{krupa2001extendingtrans}
M.~Krupa and P.~Szmolyan.
\newblock Extending slow manifolds near transcritical and pitchfork
  singularities.
\newblock {\em Nonlinearity}, 14(6):1473, 2001.

\end{thebibliography}


\newcommand{\etalchar}[1]{$^{#1}$}
\begin{thebibliography}{BBPSV04}

\bibitem[AB02]{albert2002statistical}
R{\'e}ka Albert and Albert-L{\'a}szl{\'o} Barab{\'a}si.
\newblock Statistical mechanics of complex networks.
\newblock {\em Reviews of modern physics}, 74(1):47, 2002.

\bibitem[ACS{\etalchar{+}}18]{awad2018time}
Armand Awad, Airlie Chapman, Eric Schoof, Anshu Narang-Siddarth, and Mehran
  Mesbahi.
\newblock Time-scale separation in networks: State-dependent graphs and
  consensus tracking.
\newblock {\em IEEE Transactions on Control of Network Systems}, 2018.

\bibitem[Alt13]{altafini2013consensus}
Claudio Altafini.
\newblock Consensus problems on networks with antagonistic interactions.
\newblock {\em IEEE Transactions on Automatic Control}, 58(4):935--946, 2013.

\bibitem[Alv07]{alves2007unveiling}
Nelson~A Alves.
\newblock Unveiling community structures in weighted networks.
\newblock {\em Physical Review E}, 76(3):036101, 2007.

\bibitem[AMS09]{aral2009distinguishing}
Sinan Aral, Lev Muchnik, and Arun Sundararajan.
\newblock Distinguishing influence-based contagion from homophily-driven
  diffusion in dynamic networks.
\newblock {\em Proceedings of the National Academy of Sciences},
  106(51):21544--21549, 2009.

\bibitem[B{\etalchar{+}}16]{barabasi2016network}
Albert-L{\'a}szl{\'o} Barab{\'a}si et~al.
\newblock {\em Network science}.
\newblock Cambridge university press, 2016.

\bibitem[BBPSV04]{barrat2004architecture}
Alain Barrat, Marc Barthelemy, Romualdo Pastor-Satorras, and Alessandro
  Vespignani.
\newblock The architecture of complex weighted networks.
\newblock {\em Proceedings of the national academy of sciences},
  101(11):3747--3752, 2004.

\bibitem[BD14]{bronski2014spectral}
Jared~C Bronski and Lee DeVille.
\newblock Spectral theory for dynamics on graphs containing attractive and
  repulsive interactions.
\newblock {\em SIAM Journal on Applied Mathematics}, 74(1):83--105, 2014.

\bibitem[BLM{\etalchar{+}}06]{boccaletti2006complex}
Stefano Boccaletti, Vito Latora, Yamir Moreno, Martin Chavez, and D-U Hwang.
\newblock Complex networks: Structure and dynamics.
\newblock {\em Physics reports}, 424(4-5):175--308, 2006.

\bibitem[CDS{\etalchar{+}}13]{chen2013programmable}
Yuan-Jyue Chen, Neil Dalchau, Niranjan Srinivas, Andrew Phillips, Luca
  Cardelli, David Soloveichik, and Georg Seelig.
\newblock Programmable chemical controllers made from dna.
\newblock {\em Nature nanotechnology}, 8(10):755, 2013.

\bibitem[CFQS12]{casteigts2012time}
Arnaud Casteigts, Paola Flocchini, Walter Quattrociocchi, and Nicola Santoro.
\newblock Time-varying graphs and dynamic networks.
\newblock {\em International Journal of Parallel, Emergent and Distributed
  Systems}, 27(5):387--408, 2012.

\bibitem[CWL{\etalchar{+}}16]{chen2016semidefiniteness}
Wei Chen, Dan Wang, Ji~Liu, Tamer Ba{\c{s}}ar, Karl~H Johansson, and Li~Qiu.
\newblock On semidefiniteness of signed laplacians with application to
  microgrids.
\newblock {\em IFAC-PapersOnLine}, 49(22):97--102, 2016.

\bibitem[DR96]{dumortier1996canard}
F.~Dumortier and R.~Roussarie.
\newblock {\em Canard cycles and center manifolds}, volume 577.
\newblock American Mathematical Society, 1996.

\bibitem[Eck11a]{eckhaus2011asymptotic}
W.~Eckhaus.
\newblock {\em Asymptotic analysis of singular perturbations}, volume~9.
\newblock Elsevier, 2011.

\bibitem[Eck11b]{eckhaus2011matched}
W.~Eckhaus.
\newblock {\em Matched asymptotic expansions and singular perturbations},
  volume~6.
\newblock Elsevier, 2011.

\bibitem[Fen79]{fenichel1979geometric}
N.~Fenichel.
\newblock Geometric singular perturbation theory for ordinary differential
  equations.
\newblock {\em Journal of differential equations}, 31(1):53--98, 1979.

\bibitem[FM03]{fax2003information}
J~Alexander Fax and Richard~M Murray.
\newblock Information flow and cooperative control of vehicle formations.
\newblock 2003.

\bibitem[HHML04]{holland2004using}
Barbara~R Holland, Katharina~T Huber, Vincent Moulton, and Peter~J Lockhart.
\newblock Using consensus networks to visualize contradictory evidence for
  species phylogeny.
\newblock {\em Molecular biology and evolution}, 21(7):1459--1461, 2004.

\bibitem[JKK19]{jardon2019survey}
H.~Jard\'on-Kojakhmetov and C.~Kuehn.
\newblock A survey on the blow-up method for fast-slow systems.
\newblock {\em arXiv preprint arXiv:1901.01402}, 2019.

\bibitem[JLM03]{jadbabaie2003coordination}
Ali Jadbabaie, Jie Lin, and A~Stephen Morse.
\newblock Coordination of groups of mobile autonomous agents using nearest
  neighbor rules.
\newblock {\em Departmental Papers (ESE)}, page~29, 2003.

\bibitem[Jon95]{jones1995geometric}
C.~K. R.~T. Jones.
\newblock Geometric singular perturbation theory.
\newblock In {\em Dynamical Systems}, pages 44--118. Springer, 1995.

\bibitem[KM06]{kunkel2006differential}
P.~Kunkel and V.~Mehrmann.
\newblock {\em Differential-algebraic equations: analysis and numerical
  solution}, volume~2.
\newblock European Mathematical Society, 2006.

\bibitem[Kny17]{knyazev2017signed}
Andrew~V Knyazev.
\newblock Signed laplacian for spectral clustering revisited.
\newblock {\em arXiv preprint arXiv:1701.01394}, 1, 2017.

\bibitem[KS01]{krupa2001extending}
M.~Krupa and P.~Szmolyan.
\newblock Extending geometric singular perturbation theory to nonhyperbolic
  points---fold and canard points in two dimensions.
\newblock {\em SIAM journal on mathematical analysis}, 33(2):286--314, 2001.

\bibitem[KS15]{kuehn2015multiscale}
C.~Kuehn and P.~Szmolyan.
\newblock Multiscale geometry of the {O}lsen model and non-classical relaxation
  oscillations.
\newblock {\em Journal of Nonlinear Science}, 25(3):583--629, 2015.

\bibitem[Kue15]{kuehn2015multiple}
C.~Kuehn.
\newblock {\em Multiple time scale dynamics}, volume 191.
\newblock Springer, 2015.

\bibitem[Lyn96]{lynch1996distributed}
Nancy~A Lynch.
\newblock {\em Distributed algorithms}.
\newblock Elsevier, 1996.

\bibitem[MACO91]{mohar1991laplacian}
Bojan Mohar, Y~Alavi, G~Chartrand, and OR~Oellermann.
\newblock The laplacian spectrum of graphs.
\newblock {\em Graph theory, combinatorics, and applications}, 2(871-898):12,
  1991.

\bibitem[ME10]{mesbahi2010graph}
Mehran Mesbahi and Magnus Egerstedt.
\newblock {\em Graph theoretic methods in multiagent networks}, volume~33.
\newblock Princeton University Press, 2010.

\bibitem[Mer94]{merris1994laplacian}
Russell Merris.
\newblock Laplacian matrices of graphs: a survey.
\newblock {\em Linear algebra and its applications}, 197:143--176, 1994.

\bibitem[Mes05]{mesbahi2005state}
Mehran Mesbahi.
\newblock On state-dependent dynamic graphs and their controllability
  properties.
\newblock {\em IEEE Transactions on Automatic Control}, 50(3):387--392, 2005.

\bibitem[Mor34]{moreno1934shall}
Jacob~Levy Moreno.
\newblock Who shall survive?: A new approach to the problem of human
  interrelations.
\newblock 1934.

\bibitem[Mor04]{moreau2004stability}
Luc Moreau.
\newblock Stability of continuous-time distributed consensus algorithms.
\newblock In {\em 2004 43rd IEEE conference on decision and control (CDC)(IEEE
  Cat. No. 04CH37601)}, volume~4, pages 3998--4003. IEEE, 2004.

\bibitem[Mor05]{moreau2005stability}
Luc Moreau.
\newblock Stability of multiagent systems with time-dependent communication
  links.
\newblock {\em IEEE Transactions on automatic control}, 50(2):169--182, 2005.

\bibitem[O'M91]{OMalley1991}
R.~E. O'Malley, Jr.
\newblock {\em Singular perturbation methods for ordinary differential
  equations}, volume~89 of {\em Applied Mathematical Sciences}.
\newblock Springer-Verlag, New York, 1991.

\bibitem[OS05]{olfati2005distributed}
Reza Olfati-Saber.
\newblock Distributed kalman filter with embedded consensus filters.
\newblock In {\em Proceedings of the 44th IEEE Conference on Decision and
  Control}, pages 8179--8184. IEEE, 2005.

\bibitem[OSFM07]{olfati2007consensus}
Reza Olfati-Saber, J~Alex Fax, and Richard~M Murray.
\newblock Consensus and cooperation in networked multi-agent systems.
\newblock {\em Proceedings of the IEEE}, 95(1):215--233, 2007.

\bibitem[OSM04]{olfati2004consensus}
Reza Olfati-Saber and Richard~M Murray.
\newblock Consensus problems in networks of agents with switching topology and
  time-delays.
\newblock {\em IEEE Transactions on automatic control}, 49(9):1520--1533, 2004.

\bibitem[PMC16]{proskurnikov2016opinion}
Anton~V Proskurnikov, Alexey~S Matveev, and Ming Cao.
\newblock Opinion dynamics in social networks with hostile camps: Consensus vs.
  polarization.
\newblock {\em IEEE Transactions on Automatic Control}, 61(6):1524--1536, 2016.

\bibitem[Pro13]{proskurnikov2013average}
Anton~V Proskurnikov.
\newblock Average consensus in networks with nonlinearly delayed couplings and
  switching topology.
\newblock {\em Automatica}, 49(9):2928--2932, 2013.

\bibitem[PSM16]{pan2016laplacian}
Lulu Pan, Haibin Shao, and Mehran Mesbahi.
\newblock Laplacian dynamics on signed networks.
\newblock In {\em 2016 IEEE 55th Conference on Decision and Control (CDC)},
  pages 891--896. IEEE, 2016.

\bibitem[PSV01]{pastor2001epidemic}
Romualdo Pastor-Satorras and Alessandro Vespignani.
\newblock Epidemic dynamics and endemic states in complex networks.
\newblock {\em Physical Review E}, 63(6):066117, 2001.

\bibitem[PT17]{proskurnikov2017tutorial}
Anton~V Proskurnikov and Roberto Tempo.
\newblock A tutorial on modeling and analysis of dynamic social networks. part
  i.
\newblock {\em Annual Reviews in Control}, 43:65--79, 2017.

\bibitem[PT18]{proskurnikov2018tutorial}
Anton~V Proskurnikov and Roberto Tempo.
\newblock A tutorial on modeling and analysis of dynamic social networks. part
  ii.
\newblock {\em Annual Reviews in Control}, 2018.

\bibitem[RB08]{ren2008distributed}
Wei Ren and Randal~W Beard.
\newblock {\em Distributed consensus in multi-vehicle cooperative control}.
\newblock Springer, 2008.

\bibitem[RBA05]{ren2005survey}
Wei Ren, Randal~W Beard, and Ella~M Atkins.
\newblock A survey of consensus problems in multi-agent coordination.
\newblock In {\em Proceedings of the 2005, American Control Conference, 2005.},
  pages 1859--1864. IEEE, 2005.

\bibitem[RBA07]{ren2007information}
Wei Ren, Randal~W Beard, and Ella~M Atkins.
\newblock Information consensus in multivehicle cooperative control.
\newblock {\em IEEE Control systems magazine}, 27(2):71--82, 2007.

\bibitem[SFS{\etalchar{+}}09]{schweitzer2009economic}
Frank Schweitzer, Giorgio Fagiolo, Didier Sornette, Fernando Vega-Redondo,
  Alessandro Vespignani, and Douglas~R White.
\newblock Economic networks: The new challenges.
\newblock {\em science}, 325(5939):422--425, 2009.

\bibitem[SM03]{saber2003consensus}
Reza~Olfati Saber and Richard~M Murray.
\newblock Consensus protocols for networks of dynamic agents.
\newblock 2003.

\bibitem[Str01]{strogatz2001exploring}
Steven~H Strogatz.
\newblock Exploring complex networks.
\newblock {\em nature}, 410(6825):268, 2001.

\bibitem[Tak76]{takens1976constrained}
F.~Takens.
\newblock Constrained equations; a study of implicit differential equations and
  their discontinuous solutions.
\newblock In {\em Structural stability, the theory of catastrophes, and
  applications in the sciences}, pages 143--234. Springer, 1976.

\bibitem[Tho77]{thomas1977majority}
Robert~H Thomas.
\newblock A majority consensus approach to concurrency control for multiple
  copy data bases.
\newblock Technical report, BOLT BERANEK AND NEWMAN INC CAMBRIDGE MA, 1977.

\bibitem[Tik52]{tikhonov1952systems}
A.~N. Tikhonov.
\newblock Systems of differential equations containing small parameters in the
  derivatives.
\newblock {\em Matematicheskii sbornik}, 73(3):575--586, 1952.

\bibitem[TJP07]{tanner2007flocking}
Herbert~G Tanner, Ali Jadbabaie, and George~J Pappas.
\newblock Flocking in fixed and switching networks.
\newblock {\em IEEE Transactions on Automatic control}, 52(5):863--868, 2007.

\bibitem[VB90]{vasil1990asymptotic}
AB~Vasil’Eva and VF~Butuzov.
\newblock Asymptotic methods in the theory of singular perturbations, 1990.

\bibitem[Ver05]{verhulst2005methods}
F.~Verhulst.
\newblock {\em Methods and applications of singular perturbations: boundary
  layers and multiple timescale dynamics}, volume~50.
\newblock Springer Science \& Business Media, 2005.

\bibitem[XBK07]{xiao2007distributed}
Lin Xiao, Stephen Boyd, and Seung-Jean Kim.
\newblock Distributed average consensus with least-mean-square deviation.
\newblock {\em Journal of parallel and distributed computing}, 67(1):33--46,
  2007.

\bibitem[XSK{\etalchar{+}}11]{xie2011social}
Jierui Xie, Sameet Sreenivasan, Gyorgy Korniss, Weituo Zhang, Chjan Lim, and
  Boleslaw~K Szymanski.
\newblock Social consensus through the influence of committed minorities.
\newblock {\em Physical Review E}, 84(1):011130, 2011.

\bibitem[XYJ{\etalchar{+}}04]{xia2004analyzing}
Yu~Xia, Haiyuan Yu, Ronald Jansen, Michael Seringhaus, Sarah Baxter, Dov
  Greenbaum, Hongyu Zhao, and Mark Gerstein.
\newblock Analyzing cellular biochemistry in terms of molecular networks.
\newblock {\em Annual review of biochemistry}, 73(1):1051--1087, 2004.

\bibitem[ZB14]{zelazo2014definiteness}
Daniel Zelazo and Mathias B{\"u}rger.
\newblock On the definiteness of the weighted laplacian and its connection to
  effective resistance.
\newblock In {\em 53rd IEEE Conference on Decision and Control}, pages
  2895--2900. IEEE, 2014.

\bibitem[ZB17]{zelazo2017robustness}
Daniel Zelazo and Mathias B{\"u}rger.
\newblock On the robustness of uncertain consensus networks.
\newblock {\em IEEE Transactions on Control of Network Systems}, 4(2):170--178,
  2017.

\end{thebibliography}


\newcommand{\etalchar}[1]{$^{#1}$}
\begin{thebibliography}{BBPSV04}

\bibitem[ACS{\etalchar{+}}18]{awad2018time}
A.~Awad, A.~Chapman, E.~Schoof, A.~Narang-Siddarth, and M.~Mesbahi.
\newblock Time-scale separation in networks: State-dependent graphs and
  consensus tracking.
\newblock {\em {IEEE} Trans. Control Netw. Syst.}, 2018.

\bibitem[Alt13]{altafini2013consensus}
C.~Altafini.
\newblock Consensus problems on networks with antagonistic interactions.
\newblock {\em {IEEE} Trans. Autom. Control}, 58(4):935--946, 2013.

\bibitem[Alv07]{alves2007unveiling}
N.~A. Alves.
\newblock Unveiling community structures in weighted networks.
\newblock {\em Phys. Rev. E}, 76(3):036101, 2007.

\bibitem[BBPSV04]{barrat2004architecture}
A.~Barrat, M.~Barthelemy, R.~Pastor-Satorras, and A.~Vespignani.
\newblock The architecture of complex weighted networks.
\newblock {\em Proc. Natl. Acad. Sci.}, 101(11):3747--3752, 2004.

\bibitem[BD14]{bronski2014spectral}
J.~C. Bronski and L.~DeVille.
\newblock Spectral theory for dynamics on graphs containing attractive and
  repulsive interactions.
\newblock {\em SIAM J. App. Math.}, 74(1):83--105, 2014.

\bibitem[CDS{\etalchar{+}}13]{chen2013programmable}
Y-J. Chen, N.~Dalchau, N.~Srinivas, A.~Phillips, L.~Cardelli, D.~Soloveichik,
  and G.~Seelig.
\newblock Programmable chemical controllers made from {DNA}.
\newblock {\em Nat. Nanotechnol.}, 8(10):755, 2013.

\bibitem[CFQS12]{casteigts2012time}
A.~Casteigts, P.~Flocchini, W.~Quattrociocchi, and N.~Santoro.
\newblock Time-varying graphs and dynamic networks.
\newblock {\em Int. J. Parallel Emergent Distrib. Syst.}, 27(5):387--408, 2012.

\bibitem[CWL{\etalchar{+}}16]{chen2016semidefiniteness}
W.~Chen, D.~Wang, J.~Liu, T.~Ba{\c{s}}ar, K.~H. Johansson, and L.~Qiu.
\newblock On semidefiniteness of signed laplacians with application to
  microgrids.
\newblock {\em IFAC-PapersOnLine}, 49(22):97--102, 2016.

\bibitem[DR96]{dumortier1996canard}
F.~Dumortier and R.~Roussarie.
\newblock {\em Canard cycles and center manifolds}, volume 577.
\newblock Amer. Math. Soc., 1996.

\bibitem[Eck11a]{eckhaus2011asymptotic}
W.~Eckhaus.
\newblock {\em Asymptotic analysis of singular perturbations}, volume~9.
\newblock Elsevier, 2011.

\bibitem[Eck11b]{eckhaus2011matched}
W.~Eckhaus.
\newblock {\em Matched asymptotic expansions and singular perturbations},
  volume~6.
\newblock Elsevier, 2011.

\bibitem[Fen79]{fenichel1979geometric}
N.~Fenichel.
\newblock Geometric singular perturbation theory for ordinary differential
  equations.
\newblock {\em J. Differ. Equ.}, 31(1):53--98, 1979.

\bibitem[Fie04]{Field}
M.~Field.
\newblock Combinatorial dynamics.
\newblock {\em Dynam. Syst.}, 19(3):217--243, 2004.

\bibitem[FM04]{fax2003information}
J.~A. {Fax} and R.~M. {Murray}.
\newblock Information flow and cooperative control of vehicle formations.
\newblock {\em {IEEE} Trans. Autom. Control}, 49(9):1465--1476, Sep. 2004.

\bibitem[GS17]{GolubitskyStewart2}
M.~Golubitsky and I.~Stewart.
\newblock Coordinate changes for network dynamics.
\newblock {\em Dynam. Syst.}, 32(1):80--116, 2017.

\bibitem[HHML04]{holland2004using}
B.~R. Holland, K.~T. Huber, V.~Moulton, and P.~J. Lockhart.
\newblock Using consensus networks to visualize contradictory evidence for
  species phylogeny.
\newblock {\em Mol. Biol. Evol.}, 21(7):1459--1461, 2004.

\bibitem[JKK19]{jardon2019survey}
H.~Jard\'on-Kojakhmetov and C.~Kuehn.
\newblock A survey on the blow-up method for fast-slow systems.
\newblock {\em arXiv preprint arXiv:1901.01402}, pages 1--, 2019.

\bibitem[JM03]{jadbabaie2003coordination}
A.~{Jadbabaie}, , and A.~S. {Morse}.
\newblock Coordination of groups of mobile autonomous agents using nearest
  neighbor rules.
\newblock {\em IEEE Trans. Autom. Control}, 48(6):988--1001, June 2003.

\bibitem[Jon95]{jones1995geometric}
C.~K. R.~T. Jones.
\newblock Geometric singular perturbation theory.
\newblock In {\em Dynamical Systems}, pages 44--118. Springer, 1995.

\bibitem[KM06]{kunkel2006differential}
P.~Kunkel and V.~Mehrmann.
\newblock {\em Differential-algebraic equations: analysis and numerical
  solution}, volume~2.
\newblock European Mathematical Society, 2006.

\bibitem[Kny17]{knyazev2017signed}
A.~V. Knyazev.
\newblock Signed laplacian for spectral clustering revisited.
\newblock {\em arXiv preprint arXiv:1701.01394}, 1, 2017.

\bibitem[KS01]{krupa2001extending}
M.~Krupa and P.~Szmolyan.
\newblock Extending geometric singular perturbation theory to nonhyperbolic
  points---fold and canard points in two dimensions.
\newblock {\em SIAM J. Math. Anal.}, 33(2):286--314, 2001.

\bibitem[KS15]{kuehn2015multiscale}
C.~Kuehn and P.~Szmolyan.
\newblock Multiscale geometry of the {O}lsen model and non-classical relaxation
  oscillations.
\newblock {\em J. Nonlinear Sci.}, 25(3):583--629, 2015.

\bibitem[Kue15]{kuehn2015multiple}
C.~Kuehn.
\newblock {\em Multiple time scale dynamics}, volume 191.
\newblock Springer, 2015.

\bibitem[Lyn96]{lynch1996distributed}
N.~A. Lynch.
\newblock {\em Distributed algorithms}.
\newblock Elsevier, 1996.

\bibitem[Mer94]{merris1994laplacian}
R.~Merris.
\newblock Laplacian matrices of graphs: a survey.
\newblock {\em Linear Algebra Appl.}, 197:143--176, 1994.

\bibitem[Mes05]{mesbahi2005state}
M.~Mesbahi.
\newblock On state-dependent dynamic graphs and their controllability
  properties.
\newblock {\em {IEEE} Trans. Autom. Control}, 50(3):387--392, 2005.

\bibitem[Moh91]{mohar1991laplacian}
B.~Mohar.
\newblock The laplacian spectrum of graphs.
\newblock In {\em Graph Theory, Combinatorics, and Applications}, pages
  871--898. Wiley, 1991.

\bibitem[Mor04]{moreau2004stability}
L.~Moreau.
\newblock Stability of continuous-time distributed consensus algorithms.
\newblock In {\em Proceedings of the 43rd IEEE Conference on Decision and
  Control}, volume~4, pages 3998--4003. IEEE, 2004.

\bibitem[Mor05]{moreau2005stability}
L.~Moreau.
\newblock Stability of multiagent systems with time-dependent communication
  links.
\newblock {\em {IEEE} Trans. Autom. Control}, 50(2):169--182, 2005.

\bibitem[NRS17]{NijholtRinkSanders}
E.~Nijholt, B.~Rink, and J.~Sanders.
\newblock Center manifolds of coupled cell networks.
\newblock {\em SIAM J. Math. Anal.}, 49(5):4117--4148, 2017.

\bibitem[O'M91]{OMalley1991}
R.~E. O'Malley, Jr.
\newblock {\em Singular perturbation methods for ordinary differential
  equations}, volume~89 of {\em Applied Mathematical Sciences}.
\newblock Springer-Verlag, New York, 1991.

\bibitem[OS05]{olfati2005distributed}
R.~Olfati-Saber.
\newblock Distributed kalman filter with embedded consensus filters.
\newblock In {\em Proceedings of the 44th IEEE Conference on Decision and
  Control}, pages 8179--8184. IEEE, 2005.

\bibitem[OSFM07]{olfati2007consensus}
R.~Olfati-Saber, J.~A. Fax, and R.~M. Murray.
\newblock Consensus and cooperation in networked multi-agent systems.
\newblock {\em Proc. {IEEE}}, 95(1):215--233, 2007.

\bibitem[OSM04]{olfati2004consensus}
R.~Olfati-Saber and R.~M. Murray.
\newblock Consensus problems in networks of agents with switching topology and
  time-delays.
\newblock {\em {IEEE} Trans. Autom. Control}, 49(9):1520--1533, 2004.

\bibitem[PMC16]{proskurnikov2016opinion}
A.~V. Proskurnikov, A.~S. Matveev, and M.~Cao.
\newblock Opinion dynamics in social networks with hostile camps: Consensus vs.
  polarization.
\newblock {\em {IEEE} Trans. Autom. Control}, 61(6):1524--1536, 2016.

\bibitem[Pro13]{proskurnikov2013average}
A.~V. Proskurnikov.
\newblock Average consensus in networks with nonlinearly delayed couplings and
  switching topology.
\newblock {\em Automatica}, 49(9):2928--2932, 2013.

\bibitem[PSM16]{pan2016laplacian}
L.~Pan, H.~Shao, and M.~Mesbahi.
\newblock Laplacian dynamics on signed networks.
\newblock In {\em Proceedings of the 55th Conference on Decision and Control},
  pages 891--896. IEEE, 2016.

\bibitem[RB08]{ren2008distributed}
W.~Ren and R.~W. Beard.
\newblock {\em Distributed consensus in multi-vehicle cooperative control}.
\newblock Springer, 2008.

\bibitem[RBA07]{ren2007information}
W.~Ren, R.~W. Beard, and E.~M. Atkins.
\newblock Information consensus in multivehicle cooperative control.
\newblock {\em {IEEE} Control Syst. Mag.}, 27(2):71--82, 2007.

\bibitem[SM03]{saber2003consensus}
R.~O. {Saber} and R.~M. {Murray}.
\newblock Consensus protocols for networks of dynamic agents.
\newblock In {\em Proceedings of the American Control Conference, 2003.},
  volume~2, pages 951--956, June 2003.

\bibitem[Tak76]{takens1976constrained}
F.~Takens.
\newblock Constrained equations; a study of implicit differential equations and
  their discontinuous solutions.
\newblock In {\em Structural stability, the theory of catastrophes, and
  applications in the sciences}, pages 143--234. Springer, 1976.

\bibitem[Tho79]{thomas1977majority}
R.~H. Thomas.
\newblock A majority consensus approach to concurrency control for multiple
  copy databases.
\newblock {\em ACM Trans. Database Syst.}, 4(2):180--209, jun 1979.

\bibitem[Tik52]{tikhonov1952systems}
A.~N. Tikhonov.
\newblock Systems of differential equations containing small parameters in the
  derivatives.
\newblock {\em Mat. Sb.}, 73(3):575--586, 1952.

\bibitem[TJP07]{tanner2007flocking}
H.~G. Tanner, A.~Jadbabaie, and G.~J. Pappas.
\newblock Flocking in fixed and switching networks.
\newblock {\em {IEEE} Trans. Autom. Control}, 52(5):863--868, 2007.

\bibitem[Ver05]{verhulst2005methods}
F.~Verhulst.
\newblock {\em Methods and applications of singular perturbations: boundary
  layers and multiple timescale dynamics}, volume~50.
\newblock Springer Science \& Business Media, 2005.

\bibitem[XBK07]{xiao2007distributed}
L.~Xiao, S.~Boyd, and S-J. Kim.
\newblock Distributed average consensus with least-mean-square deviation.
\newblock {\em J. Parallel Distr. Com.}, 67(1):33--46, 2007.

\bibitem[XSK{\etalchar{+}}11]{xie2011social}
J.~Xie, S.~Sreenivasan, G.~Korniss, W.~Zhang, C.~Lim, and B.~K. Szymanski.
\newblock Social consensus through the influence of committed minorities.
\newblock {\em Phys. Rev. E}, 84(1):011130, 2011.

\bibitem[ZB14]{zelazo2014definiteness}
D.~Zelazo and M.~B{\"u}rger.
\newblock On the definiteness of the weighted laplacian and its connection to
  effective resistance.
\newblock In {\em 53rd IEEE Conference on Decision and Control}, pages
  2895--2900. IEEE, 2014.

\bibitem[ZB17]{zelazo2017robustness}
D.~Zelazo and M.~B{\"u}rger.
\newblock On the robustness of uncertain consensus networks.
\newblock {\em {IEEE} Trans. Control Netw. Syst.}, 4(2):170--178, 2017.

\end{thebibliography}


\newcommand{\etalchar}[1]{$^{#1}$}
\begin{thebibliography}{MSOI{\etalchar{+}}02}

\bibitem[AS72]{abramowitz1972handbook}
M.~Abramowitz and I.A. Stegun.
\newblock {\em Handbook of Mathematical Functions}.
\newblock Dover, 1972.

\bibitem[DM15]{de2015planar}
P.~De~Maesschalck.
\newblock Planar canards with transcritical intersections.
\newblock {\em Acta Applicandae Mathematicae}, 137(1):159--184, 2015.

\bibitem[KS01]{krupa2001extendingtrans}
M.~Krupa and P.~Szmolyan.
\newblock Extending slow manifolds near transcritical and pitchfork
  singularities.
\newblock {\em Nonlinearity}, 14(6):1473, 2001.

\bibitem[MSOI{\etalchar{+}}02]{milo2002network}
Ron Milo, Shai Shen-Orr, Shalev Itzkovitz, Nadav Kashtan, Dmitri Chklovskii,
  and Uri Alon.
\newblock Network motifs: simple building blocks of complex networks.
\newblock {\em Science}, 298(5594):824--827, 2002.

\end{thebibliography}
